\documentclass[]{article}
\usepackage{cite}
\usepackage{ebrahim}
\usepackage{dsfont}
\usepackage{comment}
\usepackage{tikz}
\usepackage{mathtools}
\usetikzlibrary{tikzmark}
\usetikzlibrary{cd}

\excludecomment{excludeThis}

\excludecomment{thesisOnly} 
\includecomment{paperOnly}  

\newcommand{\spec}[1]{{\operatorname{spec}(#1)}}
\newcommand{\grspec}[1]{{\operatorname{gr-spec}(#1)}}
\newcommand{\sspec}[2]{{\operatorname{{#1}-spec}(#2)}}

\title{The Prime Spectrum and Representation Theory of the $2\times 2$ Reflection Equation Algebra\footnote{
This work will form part of the author's PhD thesis at UC Santa Barbara. 
It 
was partially supported by NSF grant DMS-1601184.
}}
\author{Ebrahim Ebrahim}

\begin{document}
\maketitle

\def \A {\mathcal{A}}
\def \GL {\operatorname{GL}}
\def \Mat {\operatorname{M}}
\def \O {\mathcal{O}}

\abstract{
The theory of generalized Weyl algebras is used to study the $2\times 2$ reflection equation algebra
$\A=\A_q(\Mat_2)$ in the case that $q$ is not a root of unity,
where the $R$-matrix used to define $\A{}$ is the 
standard one of type $A$.
Simple finite dimensional $\A{}$-modules are classified,
finite dimensional weight modules are shown to be semisimple,
$\operatorname{Aut}(\A)$ is computed,
and the prime spectrum of $\A{}$ is computed along with its Zariski topology.
Finally, it is shown that $\A$ satisfies the Dixmier-Moeglin equivalence.
}

\setlength{\parskip}{.5cm}
\parindent0ex 

\section{Introduction}

Throughout, $k$ is a field and $q\in k^\times$ is not a root of unity.

Let $n$ be a positive integer.
Consider the action by right conjugation of the algebraic group $\GL_n(k)$ of invertible $n\times n$ matrices
on the space $\Mat_n(k)$ of all $n\times n$ matrices:
$$
M\xmapsto{g\in \GL_n(k)}g^{-1}Mg.
$$
At the level of coordinate rings, the action map becomes an algebra homomorphism
\begin{equation}\label{comodule_alg}
\O(\Mat_n)\rightarrow \O(\Mat_n)\otimes \O(\GL_n),
\end{equation}
where we have dropped mention of the base field $k$ to simplify notation.
This gives $\O(\Mat_n)$ the structure of a comodule-algebra over the Hopf algebra $\O(\GL_n)$.
We shall consider what happens when this picture is carried into a quantum algebra setting.
The construction of \cite{FRT} yields a noncommutative deformation $\O_q(\Mat_n)$ of $\O(\Mat_n)$,
using the the $R$-matrix
\begin{equation}\label{rmatrix}
R_{jl}^{ik}=\begin{cases}
q&\text{if $i=j=k=l$}\\
1&\text{if $i=j$, $k=l$, and $i\neq k$}\\
q-q^{-1}&\text{if $i>j$, $i=l$, and $j= k$}\\
0&\text{otherwise}.
\end{cases}
\end{equation}
More precisely, the $k$-algebra $\O_q(\Mat_n)$ has a presentation with $n^2$ generators $\{t^i_j\bbel 0\leq i,j\leq n\}$
and the relations
$$
R^{ik}_{ab}t^a_jt^b_l = R^{ab}_{jl}t^k_bt^i_a,
$$
where repeated indices are summed over.
Further, $\O_q(\Mat_n)$ is a bialgebra in a way that matches the comultiplication on $\O(\Mat_n)$ induced by matrix multiplication in $\Mat_n$:
\begin{align*}
&\Delta(t^i_j)=t^i_k\otimes t^k_j\\
&\epsilon(t^{i}_j)=\delta^i_j.
\end{align*}
Inverting a suitable central determinant-like element in $\O_q(\Mat_n)$ yields a noncommutative deformation $\O_q(\GL_n)$
of $\O(\GL_n)$; see \cite[I.2.4]{KGnB} for example. 

One may attempt to mimic the map of (\ref{comodule_alg}) for $\O_q(\Mat_n)$ with the hope of making this algebra
a comodule-algebra over $\O_q(\GL_n)$,
\begin{equation}\label{not_comodule_alg}
\begin{array}{ccc}
\O_q(\Mat_n) & \rightarrow & \O_q(\Mat_n)\otimes \O_q(\GL_n)\\
t^i_j &\mapsto & t^k_l \otimes S(t^i_k)t^l_j,
\end{array}
\end{equation}
but such a prescription yields only a coaction map and not an algebra homomorphism.
The remedy is to replace the $\O_q(\GL_n)$-comodule $\O_q(\Mat_n)$ by a different
noncommutative deformation of $\O(\Mat_n)$.
The needed construction is provided by the transmutation theory of Majid, presented in \cite{MaBook};
it is a $k$-algebra $\A_q(\Mat_n)$ with $n^2$ generators $\{u^i_j\bbel 0\leq i,j\leq n\}$ and the relations
\begin{equation}\label{re}
R^{li}_{mn}R^{pm}_{qr}u^k_lu^n_p
\ =\ 
R^{ki}_{ml}R^{nm}_{qp}u^l_nu^p_r ,
\end{equation}
where the $R$-matrix is still (\ref{rmatrix}), the same one used to build $\O_q(\Mat_n)$.
Replacing (\ref{not_comodule_alg}) with
\begin{equation}\label{yes_comodule_alg}
\begin{array}{ccc}
\A_q(\Mat_n) & \rightarrow & \A_q(\Mat_n)\otimes \O_q(\GL_n)\\
u^i_j &\mapsto & u^k_l \otimes S(t^i_k)t^l_j
\end{array}
\end{equation}
does give an algebra homomorphism, making $\A_q(\Mat_n)$ a comodule-algebra over $\O_q(\GL_n)$
and providing a more suitable ``quantization'' of (\ref{comodule_alg}).
The algebra $\A_q(\Mat_n)$ is referred to as a \emph{braided matrix algebra} by Majid,
and as a \emph{reflection equation algebra} elsewhere in the literature.


We shall focus on the case $n=2$, the $2\times 2$ reflection equation algebra,
denoted throughout by $\A:=\A_q(\Mat_2)$.
It is generated by $u_{ij}$ for $i,j\in\{1,2\}$ with the relations
given in (\ref{re}), which simplify to:
\def\n{\\[0.3em]}
\begin{equation}\label{rels2}\begin{array}{ll}
u_{11}u_{22}=u_{22}u_{11} & \n
u_{11}u_{12}=u_{12}(u_{11}+(q^{-2}-1)u_{22}) \hspace{2em} & 
u_{21}u_{11}=(u_{11}+(q^{-2}-1)u_{22})u_{21} \hspace{2em} \n
u_{22}u_{12}=q^2u_{12}u_{22} & 
u_{21}u_{22}=q^2u_{22}u_{21} \n
u_{21}u_{12}-u_{12}u_{21}=(q^{-2}-1)u_{22}(u_{22}-u_{11}). &
\end{array}\end{equation}
Observe that $u_{12}$ and $u_{21}$ normalize the subalgebra generated by $u_{11}$ and $u_{22}$,
and they do so via inverse automorphisms of that subalgebra.
This suggests
that $\A{}$ is a generalized Weyl algebra,
a fact this paper is devoted to exploiting.



\paragraph{Brief History}
The ``reflection equation'' (\ref{re}) was first introduced by Cherednik in his study \cite{cherednik} of 
factorizable scattering on a half-line,
and reflection equation algebras later emerged from Majid's transmutation theory in \cite{Ma}.
In \cite{KSk}, Kulish and Sklyanin prove several things about $\A=\A_q(\Mat_2)$. They show that
$\A$ has a $k$-basis consisting of monomials in the generators $u_{ij}$.
They compute the center of $\A$.
They find a determinant-like element of $\A$ and they show that inverting $u_{22}$
and setting the determinant-like element equal to $1$ yields $U_q(\mathfrak{sl}_2)$,
and they note that this can be used to pull back representations of $U_q(\mathfrak{sl}_2)$ to representations of $\A$.
(We shall see in this paper that all irreducible representations
that are not annihilated by $u_{22}$
arise in this way.)
Domokos and Lenagan address $\A_q(\Mat_n)$ for general $n$ in \cite{DL}.
They show that $\A_q(\Mat_n)$ is a noetherian domain, and that it has a $k$-basis consisting of monomials in the generators $u_{ij}$.

\paragraph{Paper Outline}
Section \ref{sec:gwa} builds up the needed background and notation regarding generalized Weyl algebras (GWAs),
most of which is a collection of results from \cite{bav1}, \cite{bav_simplemod}, \cite{bav_bilateral}, and \cite{drozd}.
A description of homogeneous ideals of GWAs is given in section \ref{sec:gwa:ideals},
and localization is explored in section \ref{sec:gwa:loc}.
Section \ref{sec:gwa:gkdim} addresses GK dimension by
transporting the arguments of \cite{LMO} for skew Laurent rings into the GWA setting.
Section \ref{sec:gwa:rep} explores a certain aspect of the finite dimensional representation theory of GWAs,
focusing on the setting that will apply to the $2\times 2$ reflection equation algebra when $q$ is not a root of unity.

Section \ref{sec:rea} applies GWA theory to our $2\times 2$ reflection equation algebra $\A{}$.
Normal elements are identified and then used to compute the automorphism group.
Sections \ref{sec:aq2:simples} and \ref{sec:aq2:weight} contain a classification of finite dimensional simple $\A$-modules
and an identification of a large class of semisimple $\A$-modules.
Finally, the prime spectrum of $\A{}$ is fully worked out in section \ref{sec:rea:spec},
and some consequences are explored.

\def\G{\mathcal{G}}

\paragraph{Notation}
All rings are rings with $1$, and they are not necessarily commutative.
Given a ring $R$ and an automorphism $\sigma$ of $R$, we use $R[x;\sigma]$ to denote 
the skew polynomial ring and $R[x^\pm;\sigma]$ to denote the skew Laurent ring.
Our convention for the twisting is such that $xr=\sigma(r)x$ for $r\in R$.
If there is further twisting by a $\sigma$-derivation $\delta$, then
the notation becomes $R[x;\sigma,\delta]$.
We will also use $R((x^\pm;\sigma))$ to denote the skew Laurent series ring.
Given a subset $\G$ of a ring $R$, we indicate by $\gen{\G}$ the two-sided ideal of $R$ generated by $\G$.
When there is some ambiguity as to the ring in which ideal generation takes place,
we resolve it by using a subscript $\gen{\G}_R$ or by writing $\gen{\G}\triangleleft R$.

\paragraph{Acknowledgements}
The author would like to thank Ken Goodearl for his advice and support throughout this work.

\section{Generalized Weyl Algebras}\label{sec:gwa}

\newcommand{\sline}[2]{{\sigma^{[#1,#2]}(z)}}
\newcommand{\slinezn}[2]{{\sigma^{[#1,#2]}(z_n)}}

Generalized Weyl algebras, henceforth known as GWAs, were introduced by Bavula in \cite{bav_original}.
Examples include the ordinary Weyl algebra and the classical and quantized universal enveloping algebras of $\mathfrak{sl}_2$.

We shall define GWAs by presenting them as rings over a given base ring.
A ring $S$ \emph{over} a ring $R$, also known as an \emph{$R$-ring}, is simply a ring homomorphism $R\rightarrow S$.
A morphism $S\rightarrow S'$ of rings over $R$ is a ring homomorphism such that
$$ \begin{tikzcd}
& R \arrow[ld] \arrow[rd] & \\
S\arrow[rr] & & S' 
\end{tikzcd}$$
commutes.
Given any set $\mathcal{X}$, one can show that a \emph{free} $R$-ring on $\mathcal{X}$ exists.
This provides meaning to the notion of a presentation of a ring over $R$;
it can be thought of as a ring over $R$ satisfying a universal property described in terms of the relations.

\defn{
Let $R$ be a ring, $\sigma$ an automorphism of $R$, and $z$ an element of the center of $R$.
The GWA based on this data is the ring over $R$ generated by $x$ and $y$ subject to the relations
\begin{equation}\label{gwa_relns}
\begin{array}{l@{\hspace{3em}}l@{\hspace{3em}}l}
 yx=z & xy=\sigma(z) \\
 x r=\sigma(r)x &y r=\sigma^{-1}(r)y & \forall\ r\in R.
\end{array}
\end{equation}
}
We denote this construction by
$$
R[x,y;\sigma,z]
$$
and we adapt some useful notation from \cite{bav1} as follows. Define
$$
v_n =
\begin{cases}
x^n & n\geq 0\\
y^{(-n)} & n\leq 0
\end{cases}
$$
for $n\in\N$, and define
\begin{equation*}
\sline{j}{k}=
{\prod_{l={j}}^{k}\sigma^{l}(z)} 
\end{equation*}
for integers $j\leq k$.
We take a product over an empty index set to be $1$.
Define the following special elements of $Z(R)$:
\newcommand{\zz}[2]{[\hspace{-1.2pt}[ #1, #2 ]\hspace{-1.2pt}] }
\begin{equation} \label{eqn:zz}
\zz{n}{m}=
\begin{cases}
\sline{n+m+1}{n} & n>0,\ m<0,\ |n|\geq |m| \\
\sline{1}{n} & n>0,\ m<0,\ |n|\leq |m| \\
\sline{n+1}{n+m} & n<0,\ m>0,\ |n|\geq |m| \\
\sline{n+1}{0} & n<0,\ m>0,\ |n|\leq |m| \\
1 & \textrm{for other $n,m\in\Z$.}
\end{cases}
\end{equation}
Now we have $v_nv_m=\zz{n}{m}v_{n+m}$ for $n,m\in\Z$.

\subsection{Basic Properties}

This section lays down some basic ring-theoretic properties of GWAs,
ones we will need to reference in later sections.
The following two propositions are easy observations:

\prop{\label{prop:gwa_laurent}
There is an $R$-ring homomorphism $\phi:R[x,y;\sigma,z]\rightarrow R[x^{\pm};\sigma]$
sending $x$ to $x$ and $y$ to $zx^{-1}$.
There is also an $R$-ring homomorphism $\phi':R[x,y;\sigma,z]\rightarrow R[x^{\pm};\sigma]$
sending $x$ to $xz$ and $y$ to $x^{-1}$.}
\begin{thesisOnly}
 \begin{pf}
 Observe that the
 skew Laurent ring $R[x^{\pm};\sigma]$
 is an $R$-ring with
 $x$ and $zx^{-1}$ satisfying the defining relations for the GWA,
 and similarly for $xz$ and $x^{-1}$.
 \end{pf}
 \done
\end{thesisOnly}

\prop{\label{prop:gwa_op}
 $R[x,y;\sigma, z]$ has the alternative expression $R[y,x;\sigma^{-1},\sigma(z)]$,
 and $R[x,y;\sigma, z]^\text{op}$ can be expressed as $R^\text{op}[x,y;\sigma^{-1},\sigma(z)]$.
}
\begin{thesisOnly}
 \begin{pf}
 Check, carefully, that the GWA relations hold where needed.
 \end{pf}
 \done
\end{thesisOnly}

Proposition \ref{prop:gwa_op} roughly means that if we prove something (ring-theoretic) about $x$,
then we get a $y$ version of the result by swapping $x$ and $y$,
replacing $\sigma$ with $\sigma^{-1}$,
and replacing $z$ with $\sigma(z)$.
And if we prove something left-handed,
then we get a right-handed version of the result by
replacing $\sigma$ with $\sigma^{-1}$
and replacing $z$ with $\sigma(z)$.

\prop{\label{thm:gwa_model}
Consider a GWA $R[x,y;\sigma, z]$.
\begin{enumerate}
\item \label{thm:gwa_model:1}
It is a free left (right) $R$-module on $\{v_i\bbel i\in\Z\}$.
\item \label{thm:gwa_model:2}
It has a $\Z$-grading with homogeneous components $Rv_n=v_nR$:
$$
R[x,y;\sigma,z]=\bigoplus_{{n}\in\Z} Rv_n
$$
\item \label{thm:gwa_model:3}
It contains a copy of the ring $R$ as the subring of degree zero elements.
The subring generated by $R$ and $x$
is a skew polynomial ring $R[x;\sigma]$, and
the subring generated by $R$ and $y$
is $R[y;\sigma^{-1}]$.
\item \label{thm:gwa_model:5}
It is left (right) noetherian if $R$ is left (right) noetherian.
\end{enumerate}
}
\begin{pf}
See \cite[Lemma II.3.1.6]{rosen} for a proof of assertion \ref{thm:gwa_model:1}.
Assertions \ref{thm:gwa_model:2} and \ref{thm:gwa_model:3} are then easily shown.
\begin{paperOnly}
Assertion \ref{thm:gwa_model:5} was proven in \cite[Proposition 1.3]{bav1}.
\end{paperOnly}
\begin{excludeThis}
We shall prove these properties by defining a model ring for which they clearly hold,
and then showing that the model is isomorphic to the GWA.
As usual with this sort of proof,
most of the work goes into showing that the model ring is a ring.
Let $W$ be a free left $R$-module with basis $\{b_i\bbel i\in\Z\}$.
Define a multiplication on $W$ by
$$
\left(\sum_{i\in\Z}c_ib_i\right)\cdot
\left(\sum_{j\in\Z}d_jb_j\right)
=\sum_{k\in\Z}
 \left( \sum_{i+j=k} c_i\ \sigma^i(d_j)\ \zz{i}{j} \right)
 b_{k}.
$$
Before proving associativity, we need the identity
\begin{equation}\label{eqn:assoc_id}
\zz{a}{b}\zz{a+b}{c}=\sigma^a(\zz{b}{c})\zz{a}{b+c}
\end{equation}
for $a,b,c\in\Z$.
Instead of checking through the many cases, we may take advantage of
the known associativity of skew Laurent rings.
Let $\phi,\phi'$ be as in Proposition \ref{prop:gwa_laurent},
and denote the respective images of $v_n$ under $\phi,\phi'$
by $u_n,u_n'\in R[x^\pm;\sigma]$.
One has
$$\begin{array}{lll}
(u_au_b)u_c &= \zz{a}{b} u_{a+b} u_c &= \zz{a}{b}\zz{a+b}{c}u_{a+b+c} \\
\hspace{5mm}\rotatebox{90}{=} & \\
u_a(u_bu_c) &= u_a \zz{b}{c} u_{b+c} &= \sigma^a(\zz{b}{c}) \zz{a}{b+c} u_{a+b+c}
\end{array}$$
for $a,b,c\in\Z$, and similarly with primes.
If $a+b+c\geq 0$ then $u_{a+b+c}=x^{a+b+c}$ can be cancelled to give the desired identity.
If $a+b+c\leq 0$ then use the primed version and cancel $u_{a+b+c}'=x^{a+b+c}$ to get the identity.
Finally, we  use (\ref{eqn:assoc_id}) to show that the multiplication on $W$
is associative:
\begin{align*}
\left(
\left(\sum_{i\in\Z}c_ib_i\right)
\left(\sum_{j\in\Z}d_jb_j\right)
\right)
\left(\sum_{k\in\Z}e_kb_k\right)
&=
\sum_{m\in\Z}\ 
\sum_{i,j\in\Z}
c_i \sigma^i(d_{j-i})\sigma^{j}(e_{m-j})\zz{i}{j-i}\zz{j}{m-j}
b_m \\
&=
\sum_{m\in\Z}\ 
\sum_{k,l\in\Z}
c_k \sigma^k(d_{l})\sigma^{k+l}(e_{m-(k+l)})\zz{k}{l}\zz{k+l}{m-(k+l)}
b_m \\
&=
\sum_{m\in\Z}\ 
\sum_{k,l\in\Z}
c_k \sigma^k(d_{l})\sigma^{k+l}(e_{m-(k+l)})\sigma^k(\zz{l}{m-(k+l)})\zz{k}{m-k}
b_m \\
&=
\left(\sum_{i\in\Z}c_ib_i\right)
\left(
\left(\sum_{j\in\Z}d_jb_j\right)
\left(\sum_{k\in\Z}e_kb_k\right)
\right).
\end{align*}
It is trivial to check that $b_0$ serves as a multiplicative identity.
Left and right distributivity are straight-forward to show.
Hence $W$ is a ring.
Make $W$ a ring over $R$ via the inclusion $r\mapsto rb_0$.
Since the GWA relations hold for $b_1,b_{-1}$,
there is a unique $R$-ring homomorphism
$R[x,y;\sigma,z]\rightarrow W$
sending $x$ to $b_1$
and $y$ to $b_{-1}$.
Its inverse is the left $R$-module homomorphism that sends basis elements $b_i$
to $v_i$ for $i\in\Z$.
Hence $R[x,y;\sigma,z]$ is isomorphic to $W$ as a ring over $R$.

It is now clear that \ref{thm:gwa_model:1}-\ref{thm:gwa_model:3} hold,
with the right-handed version of \ref{thm:gwa_model:1} coming for free from Proposition \ref{prop:gwa_op}.
\end{excludeThis}
\begin{thesisOnly}
For assertion \ref{thm:gwa_model:5} let $S$ be the subring generated by $R$ and $x$.
Observe that $R[x,y;\sigma,z]$ is an over-ring of $S$ generated by $S$ and $y$
such that
\begin{equation}\label{eqn:gwa_noeth}
Sy+S=yS+S.
\end{equation}
By the skew Hilbert basis theorem,
$S=R[x;\sigma]$ is left (right) noetherian if $R$ is.
Using (\ref{eqn:gwa_noeth}), one can write a version of the Hilbert basis theorem that applies to $R[x,y;\sigma,z]$ over $S$;
see for example \cite[Theorem 2.10]{mcrob}.
That is, $R[x,y;\sigma,z]$ is left (right) noetherian if $S$ is.
\end{thesisOnly}
\end{pf}
\done

The following results are now routine.

\prop{ \label{prop:gwa_laurent2}
Let $W=R[x,y;z,\sigma]$ be a GWA.
The homomorphisms of Proposition \ref{prop:gwa_laurent}
are injective if and only if $z\in R$ is regular,
and they are isomorphisms if and only if $z\in R$ is a unit.
}
\begin{thesisOnly}
 \begin{pf}
 Assume that $z$ is regular.
 Let $\phi,\phi'$ be as in \ref{prop:gwa_laurent}.
 Consider any $w=\sum_{i\in\Z}a_iv_i\in W$ and assume that $\phi(w)=0$.
 Then
 \begin{align*}
 0
 =\phi(w)=\sum_{i\geq 0}{a_ix^i}
 +\sum_{i<0} a_iz^{-i}x^i,
 \end{align*}
 so all $a_i$ must vanish, since $z$ is regular.
 Thus $\phi$ is injective.
 Assume for the converse that $\phi$ is injective.
 Then for any $a\in R$,
 $$
 az=0 \ \Leftrightarrow \ 0 = azx^{-1} = \phi(ay) \ \Leftrightarrow \ ay = 0 \ \Leftrightarrow \ a = 0,
 $$
 so $z$ is regular.
 
 If $z$ is a unit, then $\phi$ is injective by the above, and it is surjective because it's image contains $x$ and $x^{-1}$.
 Assume for the converse that $\phi$ is an isomorphism.
 Then $x^{-1}=\phi(w)$ for some $w\in W$.
 Clearly $w=ay$ for some $a\in R$.
 Now $x^{-1} = \phi(ay) = (az)x^{-1}$, so $az=1$.
 
 The proof regarding $\phi'$ is similar.
 (Actually, $\phi$ and $\phi'$ are related via
 Proposition \ref{prop:gwa_op} and
 the isomorphism $R[x^\pm;\sigma]\cong R[x^\pm;\sigma^{-1}]$ sending $x\mapsto x^{-1}$.
 So the result for $\phi'$ follows from symmetry.)
 \end{pf}
 \done
\end{thesisOnly}

\cor{
\label{cor:gwa_domain}
A GWA $W=R[x,y;\sigma,z]$ is a domain if and only if $R$ is a domain and $z\neq 0$.
}
\begin{thesisOnly}
 \begin{pf}
 Proposition \ref{prop:gwa_laurent2} shows that if $R$ is a domain and $z\neq 0$ then
 $W$ embeds into the domain $R[x^\pm;\sigma]$.
 The converse follows from Proposition \ref{thm:gwa_model} and the fact that $yx=z$.
 \end{pf}
 \done
\end{thesisOnly}

\prop{\label{prop:gwa_zreg}
Let $W=R[x,y;\sigma,z]$ be a GWA.
Then $x,y\in W$ are regular if and only if $z\in R$ is regular.}
\begin{thesisOnly}
 \begin{pf}
 If $x$ and $y$ are regular in $W$, then so is $yx=z$.
 Suppose for the converse that $z$ is regular in $R$.
 By Proposition \ref{prop:gwa_laurent2},
 $W$ can be considered to be a subring of $R[x^\pm;\sigma]$ with $y=zx^{-1}$.
 Since $x\in R[x^\pm;\sigma]$ is regular, we conclude that $x$ and $y$ are regular.
 \end{pf}
 \done
\end{thesisOnly}

The center of a GWA is often easily described when its coefficient ring is a domain:
\prop{\label{prop:gwa_center}
Let $R$ be a domain, and let $\sigma$ be an automorphism of $R$
such that $\sigma|_{Z(R)}:Z(R)\rightarrow Z(R)$ has infinite order.
Then
$
Z(R[x,y;\sigma,z])$ is $ Z(R)^\sigma,
$
the subring of $Z(R)$ fixed by $\sigma$.
}
\begin{pf}
If $a\in Z(R)^\sigma$, then $a$ commutes with $R$, $x$, and $y$ and is therefore central.
Suppose for the converse that $a=\sum_{m\in\Z} a_mv_m$ is central.
Then $xa=ax$ and $ya=ay$ require that $\sigma(a_m)=a_m$ for all $m\in\Z$.
Given any nonzero $m\in\Z$, 
our hypothesis ensures that there is some $r\in Z(R)$ such that $\sigma^m(r)\neq r$.
Now $ra=ar$ requires $ra_m=a_m\sigma^m(r)$, so $a_m=0$.
Thus $a=a_0\in R^\sigma$.
Finally, $a$ commutes with $R$, so $a\in Z(R)^\sigma$.
\end{pf}
\done

There are similar and easily verified facts about skew Laurent polynomials and skew Laurent series:
\prop{\label{prop:Z_skew}
Let $R$ be a domain, and let $\sigma$ be an automorphism of $R$
such that $\sigma|_{Z(R)}:Z(R)\rightarrow Z(R)$ has infinite order.
Then
$
Z(R[x^\pm;\sigma]) = Z(R)^\sigma,
$
the subring of $Z(R)$ fixed by $\sigma$.
Similarly,
$Z(R((x^\pm;\sigma))) = Z(R)^\sigma$.
}
Under some stronger conditions, one can also characterize the normal elements of a GWA:
\prop{\label{prop:gwa_normalelements}
Let $R$ be a domain, and let $\sigma$ be an automorphism of $R$
such that there is an $r\in Z(R)$ which is not fixed by any nonzero power of $\sigma$.
Then the normal elements of $W=R[x,y;\sigma,z]$ are homogeneous.
}
\begin{pf}
Suppose that $a=\sum a_mv_m\in W$ is a nonzero normal element.
Then $ra=ab$ for some $b\in W$.
Looking at the highest degree and lowest degree terms of $ra$, and considering that $R$ is a domain,
$b$ must have degree $0$ in order for $ab$ to have the same highest and lowest degree terms.
Thus $b\in R$. Now $ra=ab$ becomes
$$ ra_m = a_m\sigma^m(b) $$
for all $m\in\Z$.
Since $r$ is central, we may cancel the $a_m$ whenever it is nonzero.
If $a_m$ is nonzero for multiple $m\in\Z$, then $r=\sigma^m(b)=\sigma^{m+n}(b)$
for some $m,n\in\Z$ with $n\neq 0$.
But $r=\sigma^n(r)$ would contradict our assumption on $r$, so $a$ must be homogeneous.
\end{pf}
\done

\prop{\label{prop:gwa_normalelements2}
Let $R$ be a commutative domain, $\sigma$ an automorphism, and $z\in R$ such that
$\sigma^m(z)$ is never a unit multiple of $z$ for nonzero $m\in\Z$.
Then the normal elements of $W=R[x,y;\sigma,z]$
are the $r\in R$ such that $\sigma(r)$ is a unit multiple of $r$.
}
\begin{pf}
Suppose that $r\in R$ and $\sigma(r)=ur$, where $u\in R^\times$.
Then $rR=Rr$ because $R$ is commutative,
$xr=r(ux)$, $rx=(u^{-1}x)r$,
$yr=r(yu^{-1})$, and $ry=(yu)r$.
Thus $r$ is normal in $W$.
Now assume for the converse that $a\in W$ is normal and nonzero.
By Proposition \ref{prop:gwa_normalelements},
using the fact that $z$ is not fixed by any nonzero powers of $\sigma$,
$a$ is homogeneous.
Write it as $a=a_mv_m$.

Suppose that $m\geq 0$, so that $a=a_mx^m$.
For some $b\in W$, $ax=ba$.
Clearly $b$ must have the form $b_1 x$ for some $b_1\in R$,
so we have $a_m=b_1 \sigma(a_m)$.
Thus $a_mR\subseteq \sigma(a_m)R$.
For some $c\in W$, $xa=ac$.
Then $c$ must have the form $c=c_1x$ for some $c_1\in R$,
so we have $\sigma(a_m)=a_m\sigma^m(c_1)$.
Thus $\sigma(a_m)R\subseteq a_mR$.
We conclude that $\sigma(a_m)R=a_mR$.

If $m\leq 0$, then we may use the $x\leftrightarrow y$ symmetry of Proposition 
\ref{prop:gwa_op} to apply the above argument and conclude that $\sigma^{-1}(a_m)R=a_mR$.
So in either case, $\sigma(a_m)=ua_m$ for some $u\in R^\times$.

Suppose that $m>0$, so that $a=a_mx^m$.
For some $d\in W$, $ay=da$.
Clearly $d$ must have the form $d_{-1}y$ for some $d_{-1}\in R$,
so we have
$$ a_m\sigma^m(z)=d_{-1}\sigma^{-1}(a_m)z=d_{-1}\sigma^{-1}(u^{-1})a_mz. $$
Thus, cancelling the $a_m$, $\sigma^m(z)R\subseteq zR$.
For some $e\in W$, $ae=ya$. Then $e$ must have the form $e=e_{-1}y$ for some $e_{-1}\in R$, so we have
$$ a_m\sigma^m(e_{-1})\sigma^m(z) = \sigma^{-1}(a_m)z = \sigma^{-1}(u^{-1})a_mz. $$
Thus, cancelling the $a_m$, $zR\subseteq \sigma^m(z)R$.
We conclude that $zR = \sigma^m(z)R$, contradicting the hypothesis on $z$.
Therefore one cannot have $m>0$.

If $m<0$, then we may use $x\leftrightarrow y$ symmetry to apply the above argument and conclude that
$\sigma(z)R=(\sigma^{-1})^m(\sigma(z))R$.
But this is equivalent to the contradiction $zR=\sigma^{-m}(z)R$,
so one cannot have $m<0$ either. Therefore $m=0$, and $a=a_0\in R$ with $\sigma(a)=ua$.
\end{pf}
\done

\subsection{Ideals} \label{sec:gwa:ideals}

We will establish in this section a notation for discussing the homogeneous ideals of a GWA.
We will also explore a portion of the prime spectrum of a GWA.
First, note that quotients by ideals in the coefficient ring work as they ought to:
\prop{\label{prop:gwa_quot}
Let $W=R[x,y;\sigma, z]$ be a GWA, with $J\triangleleft R$ an ideal such that $\sigma(J)=J$.
Let $I\triangleleft W$ be generated by $J$.
Then there is a canonical isomorphism 
\begin{equation} \label{eqn:gwa_quot} W/I \cong (R/J)[x,y;\hat{\sigma},z+J],\end{equation}
where $\hat{\sigma}$ is the automorphism of $R/J$ induced by $\sigma$.
}
\begin{thesisOnly}
 \begin{pf}
 Extend
 $$ R\twoheadrightarrow R/J \hookrightarrow (R/J)[x,y;\hat{\sigma},z+J] $$
 to $W$ by sending $x$ to $x$ and $y$ to $y$ and  checking that the needed GWA relations hold.
 Since the kernel of the resulting map contains $J$, it contains the ideal $I$ generated by $J$.
 This defines one direction of (\ref{eqn:gwa_quot}).
 For the other, observe that the kernel of
 $$ R \hookrightarrow W \twoheadrightarrow W/I $$
 contains $J$, and pass to the induced map $R/J \rightarrow W/I$.
 Extend this to $(R/J)[x,y;\hat{\sigma},z+J]$ by sending $x$ to $x$ and $y$ to $y$
 and checking that the needed GWA relations hold.
 The two homomorphisms just defined are inverse isomorphisms.
 \end{pf}
 \done
\end{thesisOnly}

We will generally abuse notation and reuse the labels ``$\sigma$'' and ``$z$''
instead of putting hats on things or writing out cosets.

\newcommand{\sz}[1]{\sigma^{#1}(z)}
\newcommand{\tzunderrightarrow}[3]{
\tikz[remember picture, overlay]
  \draw[-latex] ([shift={(0.8em,-0.5em)}]pic cs:#1) to [bend right] node[below]{{\scriptsize #3}}  ([shift={(0.3em,-0.5em)}]pic cs:#2);}
\newcommand{\tzunderleftarrow}[3]{
\tikz[remember picture, overlay]
  \draw[latex-] ([shift={(0.8em,-0.5em)}]pic cs:#1) to [bend right] node[below]{{\scriptsize #3}}  ([shift={(0.3em,-0.5em)}]pic cs:#2);}
\newcommand{\tzoverrightarrow}[3]{
\tikz[remember picture, overlay]
  \draw[-latex] ([shift={(0.8em,1em)}]pic cs:#1) to [bend left] node[above]{{\scriptsize #3}}  ([shift={(0.3em,1em)}]pic cs:#2);}
\newcommand{\tzoverleftarrow}[3]{
\tikz[remember picture, overlay]
  \draw[latex-] ([shift={(0.8em,1em)}]pic cs:#1) to [bend left] node[above]{{\scriptsize #3}}  ([shift={(0.3em,1em)}]pic cs:#2);}

\newcommand{\Inop}[2]{{{#1}^\text{op}_{#2}}}
\defn{ \label{defn:gwa_idealseq} Whenever $I$ is a subset of a GWA $R[x,y;\sigma,z]$ and $m\in\Z$,
$I_m$ shall denote the subset
$$I_m:=\{r\in R\bbel rv_m\in I\}$$
of $R$ and $\Inop{I}{m}$ shall denote
$$ \Inop{I}{m} := \{r\in R\bbel v_m r\in I\}. $$}

\remk{\label{remk:gwa_idealseq}
$\Inop{I}{m}$ is a notational device for working with
the symmetry $R[x,y;\sigma,z]^\text{op}=R^\text{op}[x,y;\sigma^{-1},\sigma(z)]$.
It transfers the definition of $I_m$ to the GWA structure on the opposite ring.
Note that the relation is that $\Inop{I}{m}=\sigma^{-m}(I_m)$ for all $m\in\Z$.
}

Propositions \ref{prop:gwa_homideals0} to \ref{prop:gwa_homideals2}
were essentially observed in \cite{bav_bilateral}.

\prop{\label{prop:gwa_homideals0}
Let $I$ be a right $R[x;\sigma]$-submodule
of $R[x,y;\sigma,z]$.
The $I_n$ are right ideals of $R$, and they satisfy
\begin{equation}  \label{eqn:gwa_homideals0}
\begin{array}{l@{\hspace{1cm}}l}
I_{-(n+1)}\sigma^{-n}(z)\subseteq I_{-n}  & I_n \subseteq I_{n+1} 
\end{array}
\end{equation}
for all $n\in\N$.
Thus, a homogeneous right $R[x;\sigma]$-submodule $I$
of $R[x,y;\sigma,z]$ has the form $\bigoplus_{n\in\Z}I_nv_n$ for a family
$(I_n)_{n\in\Z}$ of right ideals of $R$ satisfying
(\ref{eqn:gwa_homideals0}).
Further, any such family $(I_n)_{n\in\Z}$ defines a
right $R[x;\sigma]$-submodule of $R[x,y;\sigma,z]$ in this way.
}
\begin{thesisOnly}
 \begin{pf}
 Let $I$ be a right $R[x;\sigma]$-submodule
 of $R[x,y;\sigma,z]$,
 and let $n\in\N$.
 The $I_n$ are right ideals of $R$
 because $rs v_n=rv_n\sigma^{-n}(s)\in I$ whenever $rv_n\in I$ and $s\in R$.
 If $a\in I_n$, then
 $$
 I\ni (ax^n)x = ax^{n+1},
 $$
 so $a\in I_{n+1}$.
 And if $a\in I_{-(n+1)}$,
 then $a\sigma^{-n}(z)y^{n} = a y^{n+1}x\in I$,
 so $a\sigma^{-n}(z)\in I_{-n}$.
 This establishes (\ref{eqn:gwa_homideals0}).
 For the final assertion,
 assume that $(I_n)_{n\in\Z}$ is a family of right ideals of $R$ satisfying
 (\ref{eqn:gwa_homideals0}), and let $I=\bigoplus_{n\in\Z}I_nv_n$.
 The ring $R[x;\sigma]$ is generated by $R\cup\{x\}$,
 so one only needs to verify that $I$ is stable under right multiplication by $R$ and $x$.
 The former follows from the fact that the $I_n$ are right ideals of $R$,
 and the latter is ensured by (\ref{eqn:gwa_homideals0}).
 \end{pf}
 \done
\end{thesisOnly}

\prop{\label{prop:gwa_homideals1}
Let $I$ be a right ideal of $R[x,y;\sigma,z]$.
The $I_n$ are right ideals of $R$, and they satisfy
\begin{equation}  \label{eqn:gwa_homideals1}
\begin{array}{l@{\hspace{1cm}}l}
I_{-(n+1)} \supseteq I_{-n} &
I_n \subseteq I_{n+1} \\
I_{-(n+1)}\sigma^{-n}(z)\subseteq I_{-n} &
I_n\supseteq  I_{n+1}\sigma^{n+1}(z)
\end{array}
\end{equation}
for all $n\in\N$.
Thus, a homogeneous right ideal $I$
of $R[x,y;\sigma,z]$ has the form $\bigoplus_{n\in\Z}I_nv_n$ for a family
$(I_n)_{n\in\Z}$ of right ideals of $R$ satisfying
(\ref{eqn:gwa_homideals1}).
Further, any such family $(I_n)_{n\in\Z}$ defines a
right ideal of $R[x,y;\sigma,z]$ in this way.
}
\begin{thesisOnly}
 \begin{pf}
 To be a right ideal of $R[x,y;\sigma,z]$ is to be stable under the right multiplication by both of
 the subrings $R[x;\sigma]$
 and $R[y;\sigma^{-1}]$.
 We shall extend the assertions of Proposition \ref{prop:gwa_homideals0}
 using the symmetries indicated in Proposition \ref{prop:gwa_op}.
 We take $R[y,x;\sigma^{-1},\sigma(z)]$ to be \emph{equal} to the ring $R[x,y;\sigma,z]$,
 with our focus merely shifted to a different GWA structure.
 Applying Proposition \ref{prop:gwa_homideals0}
 to $R[y,x;\sigma^{-1},\sigma(z)]$ is a matter of swapping $x$ and $y$,
 replacing $\sigma$ by $\sigma^{-1}$, and replacing $z$ by $\sigma(z)$.
 The conditions in (\ref{eqn:gwa_homideals0}) become:
 \begin{equation}
 \begin{array}{l@{\hspace{1cm}}l}
 I_{n+1}\sigma^{n+1}(z)\subseteq I_{n}  & I_{-n} \subseteq I_{-{n+1}} .
 \end{array}
 \end{equation}
 Thus this Proposition is a consequence of Proposition \ref{prop:gwa_homideals0}
 applied to both $R[x,y;\sigma,z]$ and $R[y,x;\sigma^{-1},\sigma(z)]$.
 \end{pf}
 \done
\end{thesisOnly}

\prop{\label{prop:gwa_homideals1.5}
Let $I$ be a left ideal of $R[x,y;\sigma,z]$.
The $I_n$ are left ideals of $R$, and they satisfy
\begin{equation}  \label{eqn:gwa_homideals2}
\begin{array}{l@{\hspace{1cm}}l}
\sigma^{n+1}(I_{-(n+1)}) \supseteq \sigma^n(I_{-n}) &
\sigma^{-n}(I_n) \subseteq \sigma^{-(n+1)}(I_{n+1}) \\
\sigma^{n+1}(I_{-(n+1)})\sigma^{n+1}(z) \subseteq \sigma^n(I_{-n}) &
\sigma^{-n}(I_n) \supseteq \sigma^{-(n+1)}(I_{n+1})\sigma^{-n}(z)
\end{array}
\end{equation}
for all $n\in\N$.
Thus, a homogeneous left ideal $I$
of $R[x,y;\sigma,z]$ has the form $\bigoplus_{n\in\Z}I_nv_n$ for a family
$(I_n)_{n\in\Z}$ of left ideals of $R$ satisfying
(\ref{eqn:gwa_homideals2}).
Further, any such family $(I_n)_{n\in\Z}$ defines a
left ideal of $R[x,y;\sigma,z]$ in this way.
}
\begin{thesisOnly}
 \begin{pf}
 We again take advantage of the symmetries indicated in Proposition \ref{prop:gwa_op}.
 Applying Proposition \ref{prop:gwa_homideals1}
 to $ R[x,y;\sigma,z]^\text{op} = R^\text{op}[x,y;\sigma^{-1},\sigma(z)]$ is a matter of
 replacing $\sigma$ by $\sigma^{-1}$, and replacing $z$ by $\sigma(z)$.
 The conditions in (\ref{eqn:gwa_homideals1}) become:
 \begin{equation}  \label{eqn:gwa_homideals1_op}
 \begin{array}{l@{\hspace{1cm}}l}
 \Inop{I}{-(n+1)} \supseteq \Inop{I}{-n} &
 \Inop{I}{n} \subseteq \Inop{I}{n+1} \\
 \Inop{I}{-(n+1)}\sigma^{n+1}(z)\subseteq \Inop{I}{-n} &
 \Inop{I}{n}\supseteq  \Inop{I}{n+1}\sigma^{-n}(z),
 \end{array}
 \end{equation}
 Making the adjustment in Remark \ref{remk:gwa_idealseq} to (\ref{eqn:gwa_homideals1_op})
 yields (\ref{eqn:gwa_homideals2}).
 Thus this proposition is a consequence of Proposition \ref{prop:gwa_homideals1}
 applied to $R^\text{op}[x,y;\sigma^{-1},\sigma(z)]$.
 \end{pf}
 \done
\end{thesisOnly}

\prop{\label{prop:gwa_homideals2}
Let $I$ be an ideal of $R[x,y;\sigma,z]$.
The $I_n$ are ideals of $R$, and they satisfy
(\ref{eqn:gwa_homideals1}) and
(\ref{eqn:gwa_homideals2})
for all $n\in\N$.
Thus, a homogeneous ideal $I$
of $R[x,y;\sigma,z]$ has the form $\bigoplus_{n\in\Z}I_nv_n$ for a family
$(I_n)_{n\in\Z}$ of ideals of $R$ satisfying
(\ref{eqn:gwa_homideals1}) and (\ref{eqn:gwa_homideals2}).
Further, any such family $(I_n)_{n\in\Z}$ defines an
ideal of $R[x,y;\sigma,z]$ in this way.
}
\begin{thesisOnly}
 \begin{pf}
 Propositions \ref{prop:gwa_homideals1} and \ref{prop:gwa_homideals1.5}.
 \end{pf}
 \done
\end{thesisOnly}

We may depict
(\ref{eqn:gwa_homideals1}) and (\ref{eqn:gwa_homideals2})
by the following diagrams:

\begin{minipage}{\textwidth}
\begin{equation*}
\begin{array}{lllllllllllllllll}
\cdots & \supseteq & \tikzmark{im2}{I_{-2}} & \supseteq & \tikzmark{im1}{I_{-1}} & \supseteq & 
\tikzmark{i0}{I_0}
& \subseteq & \tikzmark{i1}{I_1} & \subseteq & \tikzmark{i2}{I_2} & \subseteq & \cdots\\[3em]
\cdots & \supseteq & \tikzmark{sim2}{\sigma^2(I_{-2})} & \supseteq & \tikzmark{sim1}{\sigma(I_{-1})} & \supseteq & 
\tikzmark{si0}{I_0}
& \subseteq & \tikzmark{si1}{\sigma^{-1}(I_1)} & \subseteq & \tikzmark{si2}{\sigma^{-2}(I_2)} & \subseteq & \cdots \quad .\\
\vspace{1em}\\
\end{array}
\end{equation*}
\tzunderrightarrow{im2}{im1}{$\sigma^{-1}(z)$}
\tzunderrightarrow{im1}{i0}{$z$}
\tzunderleftarrow{i0}{i1}{$\sigma(z)$}
\tzunderleftarrow{i1}{i2}{$\sigma^2(z)$}
\tzunderrightarrow{sim2}{sim1}{$\sigma^{2}(z)$}
\tzunderrightarrow{sim1}{si0}{$\sigma(z)$}
\tzunderleftarrow{si0}{si1}{$z$}
\tzunderleftarrow{si1}{si2}{$\sigma^{-1}(z)$}
\end{minipage}

We may also depict an alternative way of stating (\ref{eqn:gwa_homideals2}),
\begin{equation} \label{eqn:gwa_homideals2_alt}
\begin{array}{l@{\hspace{1cm}}l}
I_{-(n+1)} \supseteq \sigma^{-1}(I_{-n}) &
\sigma(I_n) \subseteq I_{n+1} \\
\sigma(I_{-(n+1)})\sigma(z) \subseteq I_{-n} &
I_n \supseteq \sigma^{-1}(I_{n+1})z,
\end{array}
\end{equation}
by the following diagram:\\[1em]
\begin{minipage}{\textwidth}
\begin{equation*}
\begin{array}{lllllllllllllllll}
\cdots &  & \tikzmark{lim2}{I_{-2}} & & \tikzmark{lim1}{I_{-1}} &  &
\tikzmark{li0}{I_0}
&  & \tikzmark{li1}{I_1} &  & \tikzmark{li2}{I_2} &  & \cdots \quad .
\end{array}
\end{equation*}
\tzoverleftarrow{lim2}{lim1}{$\sigma^{-1}$}
\tzoverleftarrow{lim1}{li0}{$\sigma^{-1}$}
\tzoverrightarrow{li0}{li1}{$\sigma$}
\tzoverrightarrow{li1}{li2}{$\sigma$}
\tzunderrightarrow{lim2}{lim1}{$\sigma,\sigma(z)$}
\tzunderrightarrow{lim1}{li0}{$\sigma,\sigma(z)$}
\tzunderleftarrow{li0}{li1}{$\sigma^{-1},z$}
\tzunderleftarrow{li1}{li2}{$\sigma^{-1},z$}
\end{minipage}

\def\p{\mathfrak{p}}

The following lemma will be useful for working out the prime spectrum of certain GWAs.
\lem{\label{lem:gwa_centralspec}
Let $W=R[x,y;\sigma,z]$ be a GWA such that $R^\sigma\subseteq R$ has the following property:
\begin{equation}\label{eqn:gwa_centralspec1}
\forall I\triangleleft R^\sigma,\ RIR\,\cap R^\sigma=I.
\end{equation}
Then there are mutually inverse inclusion-preserving bijections
\begin{equation}\label{eqn:gwa_centralspec1.5}\begin{array}{ccc}
\{I\bbel I\triangleleft R^\sigma\} & \leftrightarrow & \{ W{I}W \bbel I\triangleleft R^\sigma \} \\
I & \mapsto & W{I}W \\
\mathcal{I} \cap R^\sigma & \mapsfrom & \mathcal{I}.
\end{array}\end{equation}
Now let $S=\{W{\p}W\bbel \p\in\spec{R^\sigma}\}$ and assume that $S\subseteq\spec{W}$.
Assume also that extension of ideals to $R$ preserves intersections in the following sense:
for any family $(I_\alpha)_{\alpha\in \mathfrak{A}}$ of ideals of $R^\sigma$,
\begin{equation}\label{eqn:gwa_centralspec1.6}
\bigcap_{\alpha\in\mathfrak{A}} RI_\alpha R = R(\bigcap_{\alpha\in\mathfrak{A}}I_\alpha)R.
\end{equation}
Then (\ref{eqn:gwa_centralspec1.5}) restricts to a homeomorphism
$$ \spec{R^\sigma}  \approx S. $$
}
\begin{pf}
Given any $I\triangleleft R^\sigma$,
\begin{equation}\label{eqn:gwa_centralspec2}
W{I}W = \bigoplus_{m\in\Z} RIR\,v_m,
\end{equation}
because the right hand side satisfies the conditions of Proposition \ref{prop:gwa_homideals2}
needed to make it an ideal of $W$.
Given $I,J\triangleleft R^\sigma$ with $W{I}W\subseteq W{J}W$,
we have $RIR\subseteq RJR$ from looking at the degree zero component.
From (\ref{eqn:gwa_centralspec1}) we can then deduce that $I\subseteq J$.
The converse of this is clear: $I\subseteq J\Rightarrow W{I}W\subseteq W{J}W$.
Putting this information together, we have the inclusion-preserving correspondence (\ref{eqn:gwa_centralspec1.5}).

Now assume that $S\subseteq \spec{W}$ and that (\ref{eqn:gwa_centralspec1.6}) holds.
Let $\phi:\spec{R^\sigma}\rightarrow S$ be the restriction of (\ref{eqn:gwa_centralspec1.5}).
We show that the bijection $\phi$ is a homeomorphism.

\emph{$\phi$ is a closed map:}
Given any $I\triangleleft R^\sigma$, one has that
$\p\supseteq I$ if and only if $W{\p}W\supseteq W{I}W$,
for all $\p\in\spec{R^\sigma}$.
That is,
the collection of $\p\in\spec{R^\sigma}$
that contain $I$
is mapped by $\phi$ onto
the collection of $P\in S$ that contain $W{I}W$.

\emph{$\phi$ is continuous:}
Let $K\triangleleft W$.
Define $\mathcal{J}:=\{J\triangleleft R^\sigma\bbel K\subseteq W{J}W\}$
and $I:=\bigcap\mathcal{J}$ (with an intersection of the empty set being $R^\sigma$).
For $\p\in\spec{R^\sigma}$,
if $W{\p}W \supseteq K$, then $\p\in\mathcal{J}$, so $\p\supseteq I$.
And if $I\subseteq \p$, then
\begin{align*}
K
&\subseteq \bigcap_{J\in\mathcal{J}}W{J}W 
= W \left(\bigcap_{J\in\mathcal{J}} RJR\right) W 
= WIW 
\subseteq W{\p}W, 
\end{align*}
where
the first equality
is an application of Proposition \ref{prop:gwa_centralspec} to $R\subseteq W$,
and
the second equality 
is due to the assumption (\ref{eqn:gwa_centralspec1.6}).
We have therefore shown that the collection of $P\in S$ that contain $K$
pulls back via $\phi$ to the collection of $\p\in\spec{R^\sigma}$ that contain $I$.
\end{pf}
\done

We identify in the following proposition one situation in which the condition (\ref{eqn:gwa_centralspec1.6}) holds for a given family of ideals.
\prop{\label{prop:gwa_centralspec}
Let $A\subseteq B$ be rings such that $B$ is a free left $A$-module
with a basis $(b_j)_{j\in\mathcal{J}}$ for which
$Ab_j = b_jA$ for all $j\in\mathcal{J}$.
Let $(I_\alpha)_{\alpha\in \mathfrak{A}}$
be a family of ideals of $A$ satisfying $b_jI_\alpha\subseteq I_\alpha b_j$ for all $j$ and $\alpha$. Then
\begin{align} \label{eqn:intofideals}
\bigcap_{\alpha\in\mathfrak{A}} BI_\alpha B
&= B \left( \bigcap_{\alpha\in\mathfrak{A}}  I_\alpha \right) B.
\end{align}
}
\begin{pf}
We begin by showing that $b_j\left(\bigcap_\alpha I_\alpha\right)\subseteq
\left(\bigcap_\alpha I_\alpha\right)b_j$ for all $j$.
Consider any $j\in\mathcal{J}$ and any
$r\in \bigcap_\alpha I_\alpha$.
There is, for each $\alpha\in\mathfrak{A}$, an $r'_\alpha\in I_\alpha$
such that $b_jr=r'_\alpha b_j$.
Since $b_j$ came from a basis for $_AB$, all the $r'_\alpha$ are equal,
and so we've shown that $b_j\left(\bigcap_\alpha I_\alpha\right)\subseteq
\left(\bigcap_\alpha I_\alpha\right)b_j$ for all $j$.

Let $I$ be any ideal of $A$ satisfying $b_j I \subseteq I b_j$ for all $j$.
Observe that $\bigoplus_{j\in\mathcal{J}} Ib_j$ is then an ideal of $B$,
and hence it is the extension of $I$ to an ideal of $B$.
Applying this principle to $I=I_\alpha$ for $\alpha\in\mathfrak{A}$,
and also applying it to $I=\bigcap_\alpha I_\alpha$,
(\ref{eqn:intofideals}) follows from the fact that
\begin{align*}
& \bigcap_{\alpha\in\mathfrak{A}} \left( \bigoplus_{j\in\mathcal{J}} I_\alpha b_j \right)  
= \bigoplus_{j\in\mathcal{J}} \left( \bigcap_{\alpha\in\mathfrak{A}}  I_\alpha \right) b_j .
\tag*{\done}
\end{align*}
\end{pf}

\subsection{Localizations} \label{sec:gwa:loc}

\renewcommand{\S}{\mathcal{S}}

\begin{excludeThis}
\prop{
Let $W=R[x,y;z,\sigma]$ be a GWA with $z$ regular.
Then $\{1,x,x^{2},\ldots\}$
and $\{1,y,y^{2},\ldots\}$ are Ore sets.
\begin{pf}
Let $\S=\{1,x,x^{2},\ldots\}$, a multiplicatively closed set of regular elements.
Consider any $n\in\N$ and $w=\sum_{i\in\Z} b_i v_i\in W$.
To show that $\S$ is right Ore, we must find $m\in\N$ and $w'\in W$ so that
$wx^m=x^nw'$.
Let $I\in\Z$ be such that $b_i=0$ when $i\leq I$.
Choose $m\in\N$ such that $m\geq n-I$.
Define
$$
w' = \sum_{i\in\N} \sigma^{-n}(b_{i+n-m}\zz{i+n-m}{m})x^i.
$$
Observe that this is a well-defined element of $W$, since all but finitely many of the $b$ vanish.
We have
\begin{align*}
x^nw'
&= \sum_{i\geq 0}
x^n\sigma^{-n}(b_{i+n-m}\zz{i+n-m}{m})x^i \\
&= \sum_{i\geq 0}
b_{i+n-m}\zz{i+n-m}{m}x^{n+i} \\
&= \sum_{i\geq n-m}
b_{i}\zz{i}{m}x^{i+m} \\
&= \sum_{i\geq n-m}
b_{i}v_iv_m \\
&= \sum_{i\in\Z}
b_{i}v_iv_m
= wx^m,
\end{align*}
so $\S$ is indeed a right Ore set.
Now use Proposition \ref{prop:gwa_op}:
applying the above result to $W^\text{op} = R^\text{op}[x,y,\sigma^{-1},\sigma(z)]$
shows that $\S$ is also a \emph{left} Ore set in $W$.
And expressing $W$ as $R[y,x;\sigma^{-1},\sigma(z)]$ shows that
$\{1,y,y^{2},\ldots\}$ is an Ore set in $W$.
\end{pf}
\done

\end{excludeThis}

\prop{ \label{prop:gwa_xloc}
Let $W=R[x,y;z,\sigma]$ be a GWA with $z$ regular.
Then $\S:=\{1,x,x^{2},\ldots\}$
is an Ore set of regular elements,
and the corresponding ring of fractions is given by
the homomorphism $\phi:W\rightarrow R[x^\pm;\sigma]$ of Proposition \ref{prop:gwa_laurent}.
}
\begin{pf}
That the elements of $\S$ are regular comes from Proposition \ref{prop:gwa_zreg}.
If we can show that $\phi$ is the localization homomorphism
for a right ring of fractions of $W$ with respect to $\S$,
then by \cite[Theorem 6.2]{NNR} we will have that $\S$ is a right Ore set.
Then it will also be a left Ore set due to Proposition \ref{prop:gwa_op},
of course with the same ring of fractions, by \cite[Proposition 6.5]{NNR}.

So we have only two things to verify:
that $\phi(\S)$ is a collection of units and
that elements of $R[x^\pm;\sigma]$ have the form $\phi(w)\phi(s)^{-1}$ with $w\in W$ and $s\in S$.
The former statement is obvious.
For the latter, consider an arbitrary $p\in R[x^\pm;\sigma]$.
There is some $n\in\N$ such that $px^n\in R[x;\sigma]$.
Observe that, by Proposition \ref{prop:gwa_laurent2},
$\phi$ maps the $R$-subring $R[x;\sigma]$ of $W$ generated by $x$
isomorphically to the $R$-subring $R[x;\sigma]$ of $R[x^\pm;\sigma]$.
So
$$
p=\phi(\phi^{-1}(px^n)) \phi(x^n)^{-1},
$$
proving that $\phi$ gives a right ring of fractions.
That $\phi$ also works as a left ring of fractions is obtained for free using Proposition \ref{prop:gwa_op}.
\end{pf}
\done

\prop{\label{prop:gwa_loc}
Let $W=R[x,y;\sigma,z]$ be a GWA.
Let $\S\subseteq R$ be a right denominator set,
and assume that $\sigma(\S)= \S$.
Then $\S$ is a right denominator set of $W$,
and the associated localization map has the following description:
Let $\phi_0:R\rightarrow R\S^{-1}$ be the localization map for the right ring of fractions of $R$.
Let $\hat{\sigma}$ be the automorphism of $R\S^{-1}$ induced by $\sigma$,
and let $\hat{z}=\phi_0(z)$.
Let $\phi:W \rightarrow R\S^{-1}[x,y;\hat{\sigma},\hat{z}]$ be the extension of
$R\xrightarrow{\phi_0}R\S^{-1}\hookrightarrow R\S^{-1}[x,y;\hat{\sigma},\hat{z}]$
to $W$ that sends $x$ to $x$ and $y$ to $y$.
This is the desired localization map. In short,
$$
W\S^{-1}\ =\ (R\S^{-1})[x,y;\hat{\sigma},\hat{z}].
$$
An analogous statement holds for left denominator sets.
}
\begin{pf}
Note that $\hat{\sigma}$ exists due to our hypothesis $\sigma(\S)=\S$.
And the extension $\phi$ of $\phi_0$ exists because GWA relations hold where needed.
If we can show that $\phi$ really does define a right ring of fractions
of $R[x,y;\sigma,z]$ with respect to $\S$, then
it will follow that $\S$ is a right denominator set in $R[x,y;\sigma,z]$
(by \cite[Theorem 10.3]{NNR} for example).
Thus, three things need to be verified:
that $\phi(\S)$ is a collection of units,
that elements of $R\S^{-1}[x,y;\hat{\sigma},\hat{z}]$ have the form $\phi(w)\phi(s)^{-1}$ with $w\in W$ and $s\in S$,
and that the kernel of $\phi$ is $\{w\in W\bbel ws=0 \text{\  for some $s\in\S$}\}$.
That $\phi(\S)$ is a collection of units is obvious.

Let $\sum_{i\in\Z}a_iv_i$ 
be an arbitrary element
of $R\S^{-1}[x,y;\hat{\sigma},\hat{z}]$.
Get a ``common right denominator'' $s\in\S$ and elements $r_i$ of $R$
so that $a_i=\phi_0(r_i)\phi_0(s)^{-1}$ for all $i\in\Z$
(see \cite[Lemma 10.2a]{NNR}, noting that all but finitely many of the $a_i$ vanish).
Then
\begin{align*}
\sum a_iv_i 
&= \sum \phi_0(r_i) \phi_0(s)^{-1} v_i 
= \sum \phi_0(r_i) v_i \hat{\sigma}^{-i}(\phi_0(s)^{-1})  
= \sum \phi(r_i v_i) \phi_0(\sigma^{-i}(s))^{-1}  \\
&= \sum \phi(r_i v_i) \phi(\sigma^{-i}(s))^{-1}.
\end{align*}
After a further choice of common denominator, we see that $\sum_{i\in\Z}a_iv_i$ has the needed form.
It remains to examine the kernel of $\phi$.
Let $w=\sum r_iv_i$ be an arbitrary element of $W$.
If $ws=0$ with $s\in\S$, then $\phi(w)$ must vanish because $\phi(s)$ is a unit.
Assume for the converse that $0=\phi(w)=\sum \phi_0(r_i) v_i$.
Then $r_i\in\ker(\phi_0)$ for all $i$,
so there are $s_i\in\S$
such that $r_is_i=0$ for all $i$.
By \cite[Lemma 4.21]{NNR},
there are $b_i\in R$ such that the products $\sigma^{-i}(s_i)b_i$ are all equal to a single $s\in\S$.
Then
\begin{align*}
ws
&= \sum r_i v_i \sigma^{-i}(s_i)b_i 
= \sum r_i s_i v_i b_i =0.
\end{align*}
Thus $\ker(\phi) =\{w\in W\bbel ws=0 \text{\quad for some $s\in\S$}\}$,
and this completes the proof of the right-handed version of the theorem.
The left-handed version can be obtained for free from Proposition \ref{prop:gwa_op}.
\done
\end{pf}

We will generally abuse notation and reuse the labels ``$\sigma$'' and ``$z$''
instead of putting hats on things.

\cor{\label{cor:gwa_loc}
The localization of $W=R[x,y;\sigma,z]$
at the multiplicative set $\S$ generated by $\{\sigma^{i}(z)\bbel i\in\Z\}$
is a skew Laurent ring $(R\S^{-1})[x^\pm;\sigma]$,
where the localization map extends the one $R\rightarrow R\S^{-1}$ by
sending $x$ to $x$ and $y$ to $zx^{-1}$.
}
\begin{pf}
Use Proposition \ref{prop:gwa_loc} to describe the localization.
Then observe that it is isomorphic to a skew Laurent ring by
Proposition \ref{prop:gwa_laurent2}, since $z$ has become a unit.
\end{pf}
\done

\subsection{GK Dimension} \label{sec:gwa:gkdim}

\newcommand{\GK}[1]{{\operatorname{GK}(#1)}}
\newcommand{\subalg}[1]{{k\gen{#1}}}

Throughout this section, $R$ denotes an algebra over a field $k$,
$z$ a central element, $\sigma:R\rightarrow R$ an algebra automorphism,
and $W$ the GWA $R[x,y;\sigma,z]$.

\prop{\label{prop:gwa_gkdim}
$\GK{W}\geq \GK{R}+1$
}
\begin{pf}
Since $W$ contains a copy of the skew polynomial ring $R[x;\sigma]$,
the problem reduces to showing that $\GK{R[x;\sigma]}\geq \GK{R}+1$.
The proof is standard
(c.f. \cite[Lemma 3.4]{KL}).
\begin{thesisOnly}
We provide it for completeness.

Let $A$ be any affine subalgebra of $R$.
Let $V$ be a finite dimensional generating subspace for $A$ with $1\in V$.
Let $X=V+kx$.
Then for $n\geq 1$,
\begin{align*}
X^{2n} &= (V+kx)^{2n}
\supseteq V^n +V^nx+\cdots+V^nx^n,
\end{align*}
so $\dim(X^{2n})\geq (n+1)\dim(X^n)$. Thus,
\begin{align*}
\GK{W}
&\geq  \limsup \log_n\dim(X^{2n}) \\
&\geq  \limsup \log_n((n+1)\dim(V^n)) \\
&=  \lim_{n\rightarrow \infty}\log_n(n+1) + \limsup \log_n(\dim(V^n)) \\
&=  1  + \GK{A} .
\end{align*}
Since $A$ was an arbitrary affine subalgebra of $R$, this gives $\GK{W}\geq \GK{R}+1$.
\end{thesisOnly}
\end{pf}
\done

Under what conditions can
Proposition \ref{prop:gwa_gkdim}
be upgraded to an equality?
We look to the skew Laurent case,
i.e. the case in which $z$ is a unit, for some guidance.

\defn{An algebra automorphism $\sigma:R\rightarrow R$ is \emph{locally algebraic}
iff for each $r\in R$, $\{\sigma^n(r)\bbel n\geq 0 \}$
spans a finite dimensional subspace of $R$. Equivalently, 
$\sigma$ is locally algebraic if and only if every finite dimensional subspace of $R$
is contained in some $\sigma$-stable finite dimensional subspace of $R$.}

It was shown in \cite[Prop. 1]{LMO} that if $\sigma$ is locally algebraic,
then $\GK{R[x^\pm;\sigma]}=\GK{R}+1$.
The locally algebraic assumption was also shown to be partly necessary in \cite{Zhang},
for example when $R$ is a commutative domain with finitely generated fraction field.
So we should at least adopt the locally algebraic assumption.
Unfortunately, it is difficult to apply the result of \cite{LMO}
to a general GWA;
the process of inverting $z$, as in Corollary \ref{cor:gwa_loc},
does not make it simple to carry along GK dimension information.
For one thing, $z$ is typically not central or even normal in $W$.
Also, a locally algebraic $\sigma$ can fail to induce a locally algebraic automorphism of the localized algebra.
So we instead proceed with a direct calculation:

\def\bX{\bar{X}}

\thm{\label{thm:gwa_gkdim}
Assume that the automorphism $\sigma:R\rightarrow R$ is locally algebraic. 
Then
$$
\GK{W} = \GK{R} + 1.
$$
}
\begin{pf}
Given Proposition \ref{prop:gwa_gkdim}, it remains to show that $\GK{W}\leq\GK{R}+1$.
Let $Z$ denote the linear span of $\{\sigma^i(z)\bbel i\in\Z\}\cup \{1\}$.
Consider any affine subalgebra of $W$; let $X$ be a finite dimensional generating subspace for it.
We first enlarge $X$ to a subspace $\bX$ of the form
\begin{equation}\label{eqn:gwa_gkdim1}
\bX := \bigoplus_{|m|\leq m_0} U\,v_m,
\end{equation}
where $U$ is a finite dimensional $\sigma$-stable subspace of $R$ with $Z\subseteq U$.
Here is a procedure for doing this:
for $m\in \Z$, let $\pi_m:W\rightarrow R$ denote the $m^\text{th}$ projection map
coming from the left $R$-basis $(v_m)_{m\in\Z}$ of $W$.
Let $m_0=\max\{ |m| \bbel \pi_m(X)\neq 0 \}$.
Now $\sum_{|m|\leq m_0}\pi_m(X)$ is a finite dimensional subspace of $R$,
so it is contained in a $\sigma$-stable subspace $U$ of $R$.
It is harmless to include $Z$ in $U$
(note that $Z$ is finite dimensional because $\sigma$, and hence also $\sigma^{-1}$, is locally algebraic).
This gives us $\bX$ defined by (\ref{eqn:gwa_gkdim1}), with $X\subseteq \bX$.

\def\spacing{\hspace{2mm}}

Next, we show that 
\begin{equation}\label{eqn:gwa_gkdim2}
\bX^n\subseteq \bigoplus_{|m|\leq nm_0} U^{n+(n-1)m_0}v_m
\end{equation}
for $n\geq 1$.
It holds by definition when $n=1$, so assume that $n>1$ and that (\ref{eqn:gwa_gkdim2}) holds for $\bX^{n-1}$.
Then the induction goes through:
\begin{align}
\bX^n 
\spacing&=\spacing\bX^{n-1}\bX  
\spacing\subseteq\spacing
\left( \bigoplus_{|m|\leq (n-1)m_0} U^{n-1+(n-2)m_0}v_m \right)
\left( \bigoplus_{|m|\leq m_0} U\,v_m \right) \nonumber \\
&=\spacing
\bigoplus_{|m|\leq nm_0}
\sum_{\substack{m_1+m_2=m\\ |m_1|\leq (n-1)m_0\\ |m_2|\leq m_0}}
 U^{n-1+(n-2)m_0}v_{m_1}\, U\,v_{m_2} \nonumber\\
\spacing&\subseteq\spacing
\bigoplus_{|m|\leq nm_0}
\sum_{\substack{m_1+m_2=m\\ |m_1|\leq (n-1)m_0\\ |m_2|\leq m_0}}
U^{n+(n-1)m_0}v_{m_1+m_2} 
\spacing=\spacing
\bigoplus_{|m|\leq nm_0}
U^{n+(n-1)m_0}v_{m}. \nonumber
\end{align}
For the inclusion in the final line we used the fact,
evident from (\ref{eqn:zz}),
that $\zz{m_1}{m_2}\in Z^{\min(|m_1|,|m_2|)}\subseteq U^{\min(|m_1|,|m_2|)}$.
With (\ref{eqn:gwa_gkdim2}) established, we have
$$
\dim(\bX^n)\leq (2nm_0+1)\dim(U^{n+(n-1)m_0})
$$
for all $n\geq 1$. The theorem follows:
\begin{align*}
\GK{\subalg{X}}
\leq \GK{\subalg{\bX}}
\leq 1 + \GK{R}.\tag*{\done}
\end{align*}
\end{pf}

\subsection{Representation Theory} \label{sec:gwa:rep}

Modules over GWAs have been explored and classified under various hypotheses by several authors.
A classification of simple $R[x,y;\sigma,z]$-modules is obtained in \cite{bav_simplemod}
for $R$ a Dedekind domain with restricted minimum condition and with a condition placed on $\sigma$:
that maximal ideals of $R$ are never fixed by any nonzero power of $\sigma$.
These results are expanded in \cite{bav_bek}
and further in \cite{drozd},
where indecomposable weight modules with finite length as $R$-modules are classified
for $R$ commutative.
In the latter work, the authors introduce \emph{chain}
and \emph{circle} categories to handle maximal ideals of $R$ that have infinite and finite
$\sigma$-orbit respectively.
Another expansion of the work of \cite{bav_simplemod} was carried out in \cite{nordstrom},
where the simple $R$-torsion modules were classified relaxing all assumptions on $R$ (even commutativity),
but with the assumption that $\sigma$ acts freely on the set of maximal left ideals of $R$.
In order to establish notation and put the spotlight on a particular setting that will be of use to us,
we proceed with our own development.

\subsubsection{Simple Modules of Finite Dimension}

\def\m{\mathfrak{m}}

Let $R$ be a commutative $k$-algebra and let $W=R[x,y;\sigma,z]$ be a GWA.
Let $_WV$ be a finite dimensional simple left $W$-module.
It contains some simple left $R$-module $V_0$, which has an annihilator $\m:=\ann_RV_0\in\maxspec R$.
The automorphism $\sigma$ acts on $\maxspec R$,
and the behavior of $V$ depends largely on whether $\m$ sits in a finite or an infinite orbit.
We'd like to deal with the infinite orbit case, so assume that
$\sigma^i(\m)=\sigma^j(\m)\Rightarrow i=j$ for $i,j\in\Z$.

Let $e_0$ be a nonzero element of $V_0$,
so we have $\m=\ann_Re_0$.
For $i\in\Z$, let $e_i=v_i.e_0$.
Notice that for $i\in\Z$ and $r\in\m$, we have
$$\sigma^i(r).e_i =\sigma^i(r)v_i.e_0=v_ir.e_0=0, $$
so $\sigma^i(\m)\subseteq\ann_Re_i$.
So whenever $e_i\neq 0$, $\sigma^i(\m)=\ann_Re_i$.
We use this to argue that the subspaces $Re_i$
are independent:
Consider a vanishing combination
\begin{equation}\label{eqn:gwa:reps:1}
\sum_{i\in \mathscr{I}} r_ie_i=0
\end{equation}
where $\mathscr{I}\subseteq \Z$  is finite and $e_i\neq 0$ for $i\in\mathscr{I}$.
For any $j\in\mathscr{I}$, choose a $c_j\in \left(\prod_{i\in\mathscr{I}\setminus\{j\}}\sigma^i(\m)\right)\setminus\sigma^j(\m)$,
and apply it to (\ref{eqn:gwa:reps:1}).
The result is $c_jr_je_j=0$,
which implies that $c_jr_j\in\sigma^j(\m)$, so $r_j\in\sigma^j(\m)$ and $r_je_j=0$.

Since we assumed $V$ to be finite dimensional,
only finitely many of the $e_i$ may be nonzero.
In particular, there is some $e_{i_0}\neq 0$ such that $e_{i_0-1}=0$
(a ``lowest weight vector'').
We may as well shift our original indexing so that this $e_{i_0}$ is $e_0$.
(After all, $e_0$ was only assumed to be a nonzero element of some simple $R$-submodule of $V$
with annihilator having infinite $\sigma$-orbit,
and $e_{i_0}$ would have fit the bill just as well.)
Similarly, on the other end, there is some $n\geq 0$ so that
$e_{n-1}\neq 0$ and $e_n=0$.
Note that these definitions imply that $e_i=x^i.e_0$ is nonzero for $0\leq i\leq n-1$.

It is now clear that $\bigoplus_{i=0}^{n-1}Re_i$ is a $W$-submodule of $V$:
\begin{equation}\label{eqn:gwa:reps:2}
\begin{array}{l@{\,}l@{\,}l}
x(re_i) &= \sigma(r)xe_i &= \sigma(r)e_{i+1}\\
y(re_i) &= \sigma^{-1}(r)ye_i &= \sigma^{-1}(r)ze_{i-1} .
\end{array}
\end{equation}
So, since $_WV$ is simple, $\bigoplus_{i=0}^{n-1}Re_i=V$.
Each $Re_i$ for $0\leq i\leq n-1$
is isomorphic as an $R$-module to $R/\sigma^i(\m)$.
Knowing this and knowing that the $W$-action is described by
(\ref{eqn:gwa:reps:2}),
we have pinned down $_WV$ up to isomorphism.
Let us also pin down $e_0$ and $\m$.

Applying $xy$ and $yx$ to the extreme ``edges'' of the module shows that
$\sigma(z),\sigma^{-n+1}(z)\in\m$:
\begin{align*}
\sigma(z). e_0 = x.(y.e_0) = 0 \ &\Rightarrow\  \sigma(z)\in\m \\
z.e_{n-1} = y.(x.e_{n-1}) = 0 \ &\Rightarrow\  z\in \sigma^{n-1}(\m).
\end{align*}
Further, $n>0$ is \emph{minimal} with respect to this property:
if we had $0<i<n$ with $\sigma^{-i+1}(z)\in\m$,
then $y.e_i=0$,
so $Re_i+\cdots+Re_n$ would be a proper nontrivial submodule of $V$.

This allows us to characterize $Re_0$ as $\ann_V(y)$, as follows.
The inclusion $Re_0\subseteq \ann_V(y)$ is obvious since $y$ normalizes $R$.
Suppose that $y.\left( \sum_{i=0}^{n-1} r_ie_i \right)=0$, where $r_i\in R$.
Then $0=\sum_{i=1}^{n-1}\sigma^{-1}(r_i)ze_{i-1}$,
so for each $1\leq i\leq n-1$ we have $\sigma^{-1}(r_i)z\in\sigma^{i-1}(\m)$.
The minimality of $n$ discussed above implies that $z\notin\sigma^{i-1}(\m)$,
so we have $\sigma^{-1}(r_i)\in\sigma^{i-1}(\m)$, and hence $r_i\in\sigma^i(\m)=\ann_R(e_i)$,
for $1\leq i\leq n-1$.
So $\sum_{i=0}^{n-1}r_ie_i=r_0e_0\in Re_0$,
proving that $\ann_V(y)=Re_0$.
We have also gained a nice internal description for $\m$:
it is $\ann_R(\ann_V(y))$.
Let us record what has been established so far:

\lem{\label{lem:gwa:reps:1}
Let $_WV$ be a finite dimensional simple left $W$-module,
where $W=R[x,y;\sigma,z]$ and $R$ is a commutative $k$-algebra.
Assume that $V$ contains some simple $R$-submodule with annihilator having infinite $\sigma$-orbit.
Then $\ann_V(y)$ is just such an $R$-submodule.
Let $\m=\ann_R(\ann_V(y))$.
Then $\sigma(z)\in\m$,
there is a minimal $n>0$ such that $\sigma^{-n+1}(z)\in\m$,
and $V$ is isomorphic to
\begin{equation}\label{eqn:gwa:reps:3}
\bigoplus_{i=0}^{n-1}R/\sigma^i(\m)
\end{equation}
as an $R$-module.
Let $e_i$ denote $1\in R/\sigma^i(\m)$ as an element of (\ref{eqn:gwa:reps:3})
for $0\leq  i\leq n-1$, and let $e_{-1}=e_n=0$.
Then $_WV$ is isomorphic to (\ref{eqn:gwa:reps:3}) if (\ref{eqn:gwa:reps:3})
is given the following $W$-action:
\begin{equation*}
\begin{array}{l@{\,}l}
x(re_i) &= \sigma(r)e_{i+1}\\
y(re_i) &= \sigma^{-1}(r)ze_{i-1} .
\end{array}
\end{equation*}
}

One could check explicitly that forming the $R$-module (\ref{eqn:gwa:reps:3})
and defining actions of $x$ and $y$ according to (\ref{eqn:gwa:reps:2})
yields a well-defined, simple, and finite-dimensional module over $W$.
But we can learn a bit more about $W$ by instead realizing these modules as quotients by
certain left ideals.
We will run into a family of infinite dimensional simple modules along the way;
the construction mimics the Verma modules typical to the treatment of representations
of $\mathfrak{sl}_{2}$
\cite[II.7]{Humph}
and $U_q(\mathfrak{sl}_2)$
\cite[I.4]{KGnB}.

\defn{
Let $R$ be a commutative ring, $W=R[x,y;\sigma,z]$,
and $\m$ a maximal ideal of $R$ with infinite $\sigma$-orbit.
Define $I_\m:=W\m$ to be the left ideal of $W$ that $\m$ generates,
and define $M_\m$ to be the $\Z$-graded left $W$-module $M_\m:=W/I_\m$.
Define $e_i$ to be the image of $v_i$ in $M_\m$ for $i\in\Z$.
}

Note that
$$ I_\m = \bigoplus_{i\in \Z} \sigma^i(\m)v_i; $$
the inclusion $\supseteq$ is due to the fact that $v_i\m=\sigma^i(\m)v_i$,
and $\subseteq$ holds because the right hand side is a left ideal of $W$
(condition (\ref{eqn:gwa_homideals2}) is satisfied).

\def\I{\mathscr{I}}
\lem{
\label{lem:gwa:reps:submod}
Let $R$ be a commutative $k$-algebra, $W=R[x,y;\sigma,z]$,
and $\m$ a maximal ideal of $R$ with infinite $\sigma$-orbit.
The submodules of $M_\m$ are of the following types:
\vspace{-1em}\begin{enumerate}
\item \label{itm:gwa:submod:1} 0 or $M_\m$
\item \label{itm:gwa:submod:2} $\bigoplus_{i\geq j}Re_i$ for some $j>0$ with $\sigma^{-j+1}(z)\in\m$
\item \label{itm:gwa:submod:3} $\bigoplus_{i\leq -j'}Re_i$ for some $j'>0$ with $\sigma^{j'}(z)\in\m$
\item \label{itm:gwa:submod:4} a sum of a submodule of type \ref{itm:gwa:submod:2} and one of type \ref{itm:gwa:submod:3}.
\end{enumerate}
}
\begin{pf}
Let $S$ be a proper nontrivial submodule of $M_\m$.
We first show that $S$ is homogeneous, so that if $\sum a_ie_i\in S$ with a certain $a_je_j\neq 0$, then $e_j\in S$.

\claim{$S$ is homogeneous.}{
Suppose that $a\in S$, say $a=\sum_{i\in\I}a_ie_i$ with $\I\subseteq\Z$ finite and $a_i\in R\setminus \sigma^i(\m)$ for $i\in\I$.
Let $j\in\I$, and choose an element $c$ of 
$\left( \prod_{i\in\I\setminus \{j\}} \sigma^i(\m) \right) \setminus \sigma^j(\m)$.
Then $ca=ca_je_j\in S$. Since $c,a_j\in R\setminus \sigma^j(\m)$, $ca_j$ is a unit mod $\sigma^j(\m)$.
Hence $e_j\in S$.}

Define vector subspaces
$M^+:=\bigoplus_{i>0}Re_i$ and
$M^-:=\bigoplus_{i< 0}Re_i$
of $M_\m$.
Since $S$ is proper and homogeneous,
$$ S=(S\cap M^+)\oplus (S\cap M^-). $$
To show that $S$ is of type \ref{itm:gwa:submod:2}, \ref{itm:gwa:submod:3}, or \ref{itm:gwa:submod:4}, then,
it suffices to show that $S\cap M^+$ is a type $\ref{itm:gwa:submod:2}$ submodule when it is nonzero,
and that $S\cap M^-$ is a type $\ref{itm:gwa:submod:3}$ submodule when it is nonzero.

Assume that $S\cap M^+\neq  0$.
Then $e_j\in S$ for some $j>0$; let $j>0$ be minimal such that this happens.
By applying powers of $x$, we see that $S\cap M^+ = \bigoplus_{i\geq j}Re_i$.
Since $e_{j-1}\notin S$, $ye_j=ze_{j-1}$ must vanish.
This happens if and only if $z\in \sigma^{j-1}(\m)$, i.e. if and only if
\begin{equation}
\label{eqn:gwa:reps:4}
\sigma^{-j+1}(z)\in\m.
\end{equation}
Now assume that $S\cap M^-\neq 0$.
Let $j'>0$ be minimal such that $e_{-j'}\in S$.
By applying powers of $y$, we see that $S\cap M^- = \bigoplus_{i\leq -j'}Re_i$.
Since $e_{-j'+1}\notin S$, $xe_{-j'}=\sigma(z)e_{-j'+1}$ must vanish.
This happens if and only if $\sigma(z)\in\sigma^{-j'+1}(\m)$, i.e. if and only if
\begin{equation}
\label{eqn:gwa:reps:5}
\sigma^{j'}(z)\in\m.
\end{equation}
Finally, it is routine to check that \ref{itm:gwa:submod:1}-\ref{itm:gwa:submod:4}
are actually submodules of $M_\m$,
considering the equivalences mentioned in
(\ref{eqn:gwa:reps:4}) and 
(\ref{eqn:gwa:reps:5}).
\end{pf}
\done

This shows that $M_\m$ has a unique largest proper submodule, $N_\m$, given by
\begin{equation}\label{eqn:gwa:reps:Nm}
\hspace{75pt} N_\m:= \bigoplus_{i\leq -n'}Re_i \ \oplus\ \bigoplus_{i\geq n}Re_i \hspace{25pt} \text{(with $n,n'$ possibly $\infty$)}
\end{equation}
where $n>0$ is chosen to be minimal such that $\sigma^{-n+1}(z)\in\m$ (or $\infty$ if this never occurs)
and $n'>0$ is chosen to be minimal such that $\sigma^{n'}(z)\in\m$ (or $\infty$ if this never occurs).
For example, if $\m$ is disjoint from $\{\sigma^i(z)\bbel i\in\Z\}$, then $N_\m=0$ and $M_\m$ is simple.

\def\M{\mathscr{M}}
\thm{ \label{thm:gwa:reps}
Let $R$ be a commutative $k$-algebra and $W=R[x,y;\sigma,z]$.
\vspace{-1em}
\begin{enumerate}
\item\label{itm:gwa:reps:ass1}
Let $\m$ be a maximal ideal of $R$ with infinite $\sigma$-orbit.
Assume that $R$ is affine.
The simple module $V_\m:=M_\m/N_\m$ is finite dimensional if and only if there are $n,n'>0$ such that
$$ \sigma^{-n+1}(z),\sigma^{n'}(z)\in\m. $$
\item\label{itm:gwa:reps:ass2}
Let
\begin{equation}\label{eqn:gwa:reps:M}
\M=
\{\m\in\maxspec R\bbel \text{$\m$ has infinite $\sigma$-orbit, $\sigma(z)\in\m$, and $\sigma^{-n+1}(z)\in\m$ for some $n>0$}\}.
\end{equation}
Any finite dimensional simple left $W$-module $V$
that contains a simple $R$-submodule with annihilator having infinite $\sigma$-orbit
is isomorphic to $V_\m$ for exactly one $\m\in\mathscr{M}$, namely $\m=\ann_R(\ann_{V}(y))$.
\item\label{itm:gwa:reps:ass3}
If $\m\in\M$ and $n>0$ is minimal such that $\sigma^{-n+1}(z)\in\m$,
then $V_{\m}\cong W/(W\m+Wy+Wx^n)$.
\end{enumerate}
}
\begin{pf}
Assertion \ref{itm:gwa:reps:ass1} follows from Lemma \ref{lem:gwa:reps:submod},
the definition of $N_\m$, and the fact (due to the Nullstellensatz) that each $R/\sigma^i(\m)$ is finite dimensional when $R$ is affine.
For assertion \ref{itm:gwa:reps:ass2}, suppose that $_WV$ is simple, finite dimensional,
and contains a simple $R$-submodule with annihilator having infinite $\sigma$-orbit.
Lemma \ref{lem:gwa:reps:1} pins $V$ down as isomorphic to the left $W$-module in (\ref{eqn:gwa:reps:3}).
This construction is in turn isomorphic to $V_\m$, where $\m=\ann_R(\ann_V(y))$,
and the lemma guarantees that $\sigma(z),\sigma^{-n+1}(z)\in\m$ for some $n>0$.
Hence $\m\in\M$ and $V\cong V_\m$.
Since
$ \m = \ann_R(\ann_{V_\m}(y)) $,
no two $V_\m$ for $\m\in\M$ can be isomorphic.
Assertion \ref{itm:gwa:reps:ass3} amounts to the fact that, under the given hypotheses, $N_\m$
is the submodule of $M_\m$ generated by the cosets $y+I_\m$ and $x^n+I_\m$.
\end{pf}
\done

\subsubsection{Weight Modules of Finite Dimension}

\def\m{\mathfrak{m}}
In further pursuit of finite dimensional modules,
we now explore a class of modules that includes the semisimple ones.
We continue with the notation $W=R[x,y;\sigma,z]$
and the assumption that $R$ is a commutative $k$-algebra.
Let $_WX$ be finite dimensional and semisimple.
Consider the $R$-submodule spanned by annihilators of maximal ideals,
$$
S:=\sum_{\m\in\maxspec R}\ann_X\m.
$$
It is in fact a $W$-submodule of $X$, since $x$ and $y$ map
$\ann_X\m$ into $\ann_X\sigma(\m)$ and $\ann_X\sigma^{-1}(\m)$ respectively.
Since we assumed $X$ to be semisimple, $S$ has a direct sum complement $S'$ in $_WX$.
If $S'$ were nonzero, then it contains some simple $R$-submodule which is then annihilated by some maximal ideal of $R$, contradicting
$S'\cap S=0$.
Thus our assumption that $_WX$ is semisimple requires $X=S$.
We now wonder when this condition is sufficient for semisimplicity.

\newcommand{\supp}{\operatorname{supp}}
\defn{\label{def:gwa:weightmod}
Let $R$ be a commutative $k$-algebra.
A $W$-module where $W=R[x,y;\sigma,z]$ is a \emph{weight} module
if and only if it is semisimple as an $R$-module.
Note that this is equivalent to saying that $X$ is spanned by annihilators of maximal ideals of $R$.
The \emph{support} $\supp X$ of an $R$-module $X$
is the collection of maximal ideals $\m$ of $R$ such that
$\ann_X\m$ is nonzero.
}

Let us collect some elementary facts about weight modules for use in the coming semisimplicity theorem.

\def\n{\mathfrak{n}}
\prop{\label{prop:gwa:weightmod3}
Let $R$ be a commutative $k$-algebra,
$X$ a semisimple $R$-module, and $_RY\leq {_R}X$.
Then
$$ X=\bigoplus_{\m\in\maxspec R}\ann_X\m $$
and the canonical map $X\rightarrow X/Y$ induces an isomorphism of $R$-modules
$$
(\ann_X\m)/(\ann_Y\m) \cong \ann_{X/Y}\m
$$
for each $\m\in\maxspec R$.
}
\begin{pf}
By assumption, $X$ is a direct sum of simple $R$-submodules.
Each simple $R$-submodule is isomorphic to $R/\m$ for some $\m\in\maxspec R$.
Thus $X$ is a direct sum of the $\ann_X\m$;
each $\ann_X\m$ is actually just the $(R/\m)$-homogeneous component of $X$.
Since $Y$ is a submodule of $X$, it is semisimple and has its own decomposition
\begin{equation}\label{eqn:gwa:reps:weight3}
Y=\bigoplus_{\m\in\maxspec R}\ann_Y\m = \bigoplus_{\m\in\maxspec R} Y\cap\ann_X\m.
\end{equation}
Fix an $\m\in\maxspec R$.
It is clear that the canonical map $X\rightarrow X/Y$
restricts to an $R$-homomorphism $\ann_X\m\rightarrow \ann_{X/Y}\m$
with kernel $Y\cap \ann_X\m=\ann_Y\m$.
To see that it is surjective, consider any $x+Y\in \ann_{X/Y}\m$.
Write $x$ as $\sum_{\n\in \maxspec R}x_\n$,
where $x_\n\in \ann_X\n$.
Since $\m x\subseteq Y$, the decomposition (\ref{eqn:gwa:reps:weight3}) gives $\m x_\n\subseteq Y$ for all $\n$.
When $\n\neq \m$, this implies that $x_\n\in Y$ since $\m$ contains a unit mod $\n$.
Thus $x+Y$ is the image of $x_\m$ under $\ann_X\m\rightarrow \ann_{X/Y}\m$.
\end{pf}
\done

\newcommand{\chaintype}{chain-type} 

Since we only focused on simple finite dimensional $W$-modules of a certain type,
we will only attempt to get at the weight modules whose composition factors are of that type.
Adapting the ``chain'' and ``circle'' terminology from \cite{drozd}:
\defn{
Let $R$ be a commutative $k$-algebra and let
$W=R[x,y;\sigma,z]$.
A finite dimensional module
$_WX$ is of \emph{\chaintype{}} if and only if every $\m\in\supp X$
has infinite $\sigma$-orbit.
}

\prop{\label{prop:gwa:weightmod2}
Let $R$ be a commutative $k$-algebra,
$W=R[x,y;\sigma,z]$, and $_WX$ a \chaintype{} finite dimensional weight module.
Let $\M$ be as in (\ref{eqn:gwa:reps:M}).
Then each composition factor of $X$ has the form $V_\m$ for some  $\m\in\M$, and
\begin{equation}\label{eqn:gwa:reps:weightmod2}
\supp X\cap \M = \{ \m\in\M\bbel \text{$V_\m$ is a composition factor of $_WX$} \}.
\end{equation}
}
\begin{pf}
Choose a $W$-module composition series $0=X_0 \subsetneq X_1 \subsetneq \cdots \subsetneq X_r = X$.
It can be refined into a composition series for $_RX$,
so since $_RX$ is semisimple we have:
\begin{equation}\label{eqn:gwa:reps:weightmod5}
_RX\cong \bigoplus_{i=1}^r\bigoplus\{(R/\m)^{(k)} \bbel \text{$R/\m$ is a composition factor of $X_i/X_{i-1}$ with multiplicity $k$}\}.
\end{equation}
In particular,
each $X_i/X_{i-1}$ contains some simple $R$-submodule
whose annihilator comes from $\supp X$ and therefore has infinite $\sigma$-orbit.
Theorem \ref{thm:gwa:reps} applies:
for $1\leq i\leq r$, $X_i/X_{i-1}\cong V_{\m_i}$ for a unique $\m_i\in \M$.
The right hand side of (\ref{eqn:gwa:reps:weightmod2}) is then $\{\m_1,\ldots,\m_r\}$.
For each $1\leq i\leq r$, let $n_i>0$ be minimal such that $\sigma^{-n_i+1}(z)\in\m_i$.

Knowing that $X_{i}/X_{i-1}\cong V_{\m_i}$, we can read off the support of $X$ from (\ref{eqn:gwa:reps:weightmod5}):
$$\supp X = \{ \sigma^\ell(\m_i) \bbel \text{$1\leq i\leq r$ and $0\leq \ell\leq n_i-1$}  \}.$$
Suppose that $\sigma^{\ell}(\m_i)\in \M$, with
$1\leq i\leq r$ and $0\leq \ell\leq n_i-1$.
Then $\sigma(z)\in \sigma^{\ell}(\m_i)$, so $\sigma^{-\ell+1}(z)\in\m_i$.
The minimality of $n_i$ forces $\ell=0$.
This proves that $\supp X\cap \M = \{\m_1,\ldots,\m_r\}$,
and the latter is the right hand side of (\ref{eqn:gwa:reps:weightmod2}).
\end{pf}
\done

Next, we identify a condition on $\supp X\cap \M$
that we will show guarantees semisimplicity for $_WX$.

\newcommand{\hasSepStrings}{has separated chains} 
\newcommand{\haveSepStrings}{have separated chains} 
\newcommand{\havingSepStrings}{having separated chains} 
\def\M{\mathscr{M}}
\def\m{\mathfrak{m}}
\defn{ \label{defn:gwa:reps:sepStrings}
Let $R$ be a commutative $k$-algebra, $\sigma$ an automorphism, and $z\in R$.
Let $\M$ be as in (\ref{eqn:gwa:reps:M}).
A subset $S\subset\M$
\emph{\hasSepStrings{}} if and only if the following holds:
whenever $\m\in S$ and $n>0$ is minimal such that $\sigma^{-n+1}(z)\in\m$,
it follows that $\sigma^n(\m)\notin S$.
}

\prop{\label{prop:gwa:reps:sepStrings}
Let $R$ be a commutative $k$-algebra, $\sigma$ an automorphism, and $z\in R$.
Let $\M$ be as in (\ref{eqn:gwa:reps:M}),
and suppose that $S\subset\M$ \hasSepStrings{}.
Then given $\m,\m'\in S$ and $n,n'>0$ minimal such that $\sigma^{-n+1}(z)\in\m$ and $\sigma^{-n'+1}(z)\in\m'$,
\vspace{-1em}\begin{enumerate}
\item \label{itm:gwa:reps:sepStrings:1}
$\m'\in \{\m,\sigma(\m),\ldots,\sigma^{n-1}(\m)\}$ \hspace{1pt} only if \hspace{1pt} $\m'=\m$
\item \label{itm:gwa:reps:sepStrings:2}
$\sigma^{-1}(\m'),\sigma^{n'}(\m')\notin \{\m,\sigma(\m),\ldots,\sigma^{n-1}(\m)\}$.
\end{enumerate}
}
\begin{pf}
Suppose that $\m'=\sigma^\ell(\m)$, where $0\leq \ell \leq n-1$.
Then $\sigma(z)\in \sigma^\ell(\m)$,
so $\sigma^{-\ell+1}(z)\in\m$.
The minimality of $n$ then forces $\ell=0$, whence $\m'=\m$.

Suppose that $\sigma^{-1}(\m')=\sigma^\ell(\m)$, where $0\leq \ell \leq n-1$.
Then $\sigma(z)\in\m'=\sigma^{\ell+1}(\m)$, so $\sigma^{-(\ell+1)+1}(z)\in\m$.
The minimality of $n$ then forces $\ell+1=n$, which gives $\sigma^n(\m)=\m'\in S$.
This contradicts the assumption that $S$ \hasSepStrings{}.

Suppose that $\sigma^{n'}(\m')=\sigma^\ell(\m)$, where $0\leq \ell \leq n-1$.
Then we have $\sigma^{-\ell+1}(z)\in\m$, since $\sigma^{-n'+1}(z)\in\m'$.
The minimality of $n$ then forces $\ell=0$, which gives $\sigma^{n'}(\m')=\m\in S$,
contradicting the assumption that $S$ \hasSepStrings{}.
\end{pf}
\done

\def\n{\mathfrak{n}}
\thm{\label{thm:gwa:reps:sepStrings}
Let $R$ be a commutative $k$-algebra, let $W=R[x,y;\sigma,z]$, and let $\M$ be as in (\ref{eqn:gwa:reps:M}).
Let $X$ be a \chaintype{} finite dimensional weight left $W$-module.
If $\supp X \cap \M$ \hasSepStrings{},
then $X$ is semisimple.
}
\begin{pf}
Assume the hypotheses. Choose a composition series for $_WX$:
$$
0=X_0 \subsetneq X_1 \subsetneq \cdots \subsetneq X_r = X.
$$
For $1\leq i\leq r$, $X_i/X_{i-1}\cong V_{\m_i}$ for a unique $\m_i\in \M$,
and $\{\m_1,\ldots,\m_r\}$ \hasSepStrings{}
(Proposition \ref{prop:gwa:weightmod2}).
Let $\n_1,\ldots, \n_s$ be the distinct items among $\m_1,\ldots,\m_r$, 
with respective multiplicities $t_1,\ldots,t_s$.
For each $1\leq j\leq s$, let $n_j>0$ be minimal such that $\sigma^{-n_j+1}(z)\in\n_j$.

\def\a{\mathfrak{a}}
For any $\a\in\maxspec R$, we iteratively apply Proposition \ref{prop:gwa:weightmod3} to obtain:
\begin{align}
\dim_{R/\a}\ann_X\a
&= \dim_{R/\a}\ann_{X_r/X_{r-1}} \a + \dim_{R/\a}\ann_{X_{r-1}}\a
= \cdots 
= \sum_{i=1}^r \dim_{R/\a}\ann_{X_i/X_{i-1}} \a \nonumber\\
&= \sum_{i=1}^r \dim_{R/\a}\ann_{V_{\m_i}} \a \label{eqn:gwa:reps:theabove}
\end{align}
Fix a $j$ with $1\leq j\leq s$.
Apply (\ref{eqn:gwa:reps:theabove}) to the case $\a=\n_j$ and use
Proposition \ref{prop:gwa:reps:sepStrings}.\ref{itm:gwa:reps:sepStrings:1}
to obtain
$$ \dim_{R/\n_j}\ann_X\n_j = t_j. $$
Apply (\ref{eqn:gwa:reps:theabove}) to the cases $\a=\sigma^{-1}(\n_j)$ and $\a=\sigma^{n_j}(\n_j)$ and use
Proposition \ref{prop:gwa:reps:sepStrings}.\ref{itm:gwa:reps:sepStrings:2}
to obtain
$$ \ann_X(\sigma^{-1}(\n_j)) = \ann_X(\sigma^{n_j}(\n_j)) = 0. $$
Let $b^j_1,\ldots,b^j_{t_j}$ be an $(R/\n_j)$-basis for $\ann_X\n_j$.
Each $Wb^j_u$ is a nonzero homomorphic image of $_WW$ in which $\n_j$, $y$, and $x^{n_j}$ are killed:
$\n_j$ is killed because $b^j_u$ came from $\ann_X\n_j$, and $y$ and $x^{n_j}$ are killed because they
map $\ann_X\n_j$ into $\ann_X(\sigma^{-1}(\n_j))$ and $\ann_X(\sigma^{n_j}(\n_j))$ respectively.
By Theorem \ref{thm:gwa:reps}.\ref{itm:gwa:reps:ass3}, it follows that
each $Wb^j_u$ is isomorphic to $M_{\n_j}/N_{\n_j}=:V_{\n_j}$.

Do this for all $j$. Let
$$ S=\sum_{j=1}^s\sum_{u=1}^{t_j} Wb_u^j, $$
a semisimple $W$-submodule of $X$.
Since $\ann_X\n_j$ is $R$-spanned by $b^j_1,\ldots,b^j_{t_j}$, we have
(by Proposition \ref{prop:gwa:weightmod3})
$$
\ann_{X/S}\n_j \cong (\ann_X\n_j)/(\ann_S\n_j) = 0
$$
for all $j$.
By the Jordan-H\"older theorem, any simple $W$-submodule of $X/S$ is isomorphic to $V_{\n_j}$
for some $j$. Therefore $X/S$ must be $0$. That is, $X=S$ is semisimple.
\done
\end{pf}

If $\M$ as a whole \hasSepStrings{}, then we conclude from this theorem that
\emph{all} \chaintype{} finite dimensional weight modules are semisimple.
There is a converse:
\prop{\label{prop:gwa:reps:sepStrings4}
Let $R$ be an affine commutative $k$-algebra, let $W=R[x,y;\sigma,z]$, and let $\M$ be as in (\ref{eqn:gwa:reps:M}).
If $\M$ does not \haveSepStrings{}, then there is a \chaintype{} finite dimensional weight left  $W$-module that is not semisimple.
}
\begin{pf}
If $\M$ does not \haveSepStrings{}, there is some $\m\in\M$ such that $\sigma^n(\m)\in\M$, where
$n>0$ is minimal such that $\sigma^{-n+1}(z)\in\m$.
Let $n'>0$ be minimal such that $\sigma^{-n'+1}(z)\in\sigma^n(\m)$.
Then $\m$ contains $\sigma(z)$, $\sigma^{-n+1}(z)$, and $\sigma^{-(n+n')+1}(z)$.
Hence 
$$S:=\left( \bigoplus_{i\leq -1}Re_i \right) \oplus \left( \bigoplus_{i\geq n+n'} Re_i \right)$$
is a submodule of $M_\m$, by Lemma \ref{lem:gwa:reps:submod}.
Let $_WX=M_\m/S$. 
This is isomorphic to $\bigoplus_{0\leq i < n+n'}R/\sigma^i(\m)$ as an $R$-module, so $_WX$ is
a \chaintype{} finite dimensional weight left  $W$-module.
Since $M_\m$ contains a unique largest proper submodule
$$N_\m = \left( \bigoplus_{i\leq -1}Re_i \right) \oplus \left( \bigoplus_{i\geq n} Re_i \right) $$
and $N_\m$ properly contains $S$,
$X$ contains a unique largest proper nontrivial submodule $N_\m/S$. Therefore $X$ cannot be semisimple.
\end{pf}
\done

\thm{\label{thm:gwa:reps:sepStrings2}
Let $R$ be an affine commutative $k$-algebra, let $W=R[x,y;\sigma,z]$, and let $\M$ be as in (\ref{eqn:gwa:reps:M}).
The following are equivalent:
\vspace{-1em}\begin{enumerate}
\item\label{itm:gwa:reps:sepStrings3}
All \chaintype{} finite dimensional weight left $W$-modules are semisimple.
\item\label{itm:gwa:reps:sepStrings4}
$\M$ \hasSepStrings{}.
\item\label{itm:gwa:reps:sepStrings6}
For any maximal ideal $\m$ of $R$ with infinite $\sigma$-orbit,
there are no more than two integers $i$ such that $\sigma^i(z)\in\m$.
\item\label{itm:gwa:reps:sepStringsNew}
For any $\m\in\M$, there is exactly one $n>0$ such that $\sigma^{-n+1}(z)\in\m$.
\end{enumerate}
}
\begin{pf}
The equivalence \ref{itm:gwa:reps:sepStrings3}$\Leftrightarrow$\ref{itm:gwa:reps:sepStrings4}
is due to Theorem \ref{thm:gwa:reps:sepStrings} and Proposition \ref{prop:gwa:reps:sepStrings4}.

\ref{itm:gwa:reps:sepStrings4} $\Rightarrow$ \ref{itm:gwa:reps:sepStringsNew}:
Assume that \ref{itm:gwa:reps:sepStringsNew} fails.
Let $\m$ be in $\M$ and let $i<j$ be positive integers such that $\sigma^{-i+1}(z),\sigma^{-j+1}(z)\in\m$.
We may assume that $i>0$ is minimal such that $\sigma^{-i+1}(z)\in\m$.
Observe that $\sigma(z),\sigma^{-(j-i)+1}(z)\in\sigma^i(\m)$.
This implies that $\sigma^i(\m)\in\M$, so $\M$ does not \haveSepStrings{}.

\ref{itm:gwa:reps:sepStringsNew} $\Rightarrow$ \ref{itm:gwa:reps:sepStrings6}:
Assume that \ref{itm:gwa:reps:sepStrings6} fails;
let $\m$ be a maximal ideal of $R$ with infinite $\sigma$-orbit and with
$\sigma^i(z),\sigma^j(z),\sigma^k(z)\in\m$, where $i<j<k$ are integers.
We may assume 
that $j>i$ is minimal such that $\sigma^j(z)\in\m$ and
that $k>j$ is minimal such that $\sigma^k(z)\in\m$.
Let $\n:=\sigma^{-k+1}(\m)$.
Observe that $\n\in\M$ since $\sigma(z)\in \n$ and $\sigma^{-(k-j)+1}(z)\in\n$.
Since $\sigma^{-(k-i)+1}(z)\in\n$ as well, with $k-i\neq k-j$, we see that \ref{itm:gwa:reps:sepStringsNew} fails.

\ref{itm:gwa:reps:sepStrings6} $\Rightarrow$ \ref{itm:gwa:reps:sepStrings4}:
Suppose that $\M$ does not \haveSepStrings{}.
Then there is some $\m\in\M$ such that $\sigma^n(\m)\in\M$, where
$n>0$ is minimal such that $\sigma^{-n+1}(z)\in\m$.
Let $n'>0$ be minimal such that $\sigma^{-n'+1}(z)\in\sigma^n(\m)$.
Since $\sigma(z),\sigma^{-n+1}(z),\sigma^{-(n+n')+1}(z) \in \m$,
\ref{itm:gwa:reps:sepStrings6} fails to hold.
\end{pf}
\done

\section{The $2\times 2$ Reflection Equation Algebra} \label{sec:rea}

We now shift our focus to a specific GWA, the algebra $\A$ defined in (\ref{rels2}).
Define an automorphism $\sigma$ of the polynomial ring $k[u_{22},u_{11},z]$ by
$$\begin{array}{lll}
\sigma(u_{22})&=&q^2u_{22}\\
\sigma(u_{11})&=&u_{11}+(q^{-2}-1)u_{22}\\
\sigma(z)&=&z+(q^{-2}-1)u_{22}(u_{22}-u_{11}).
\end{array}$$
The algebra $\A$ is a GWA over the above polynomial algebra,
with $x$ being $u_{21}$ and $y$ being $u_{12}$:
$$\A\cong k[u_{22},u_{11},z][x,y;\sigma,z].$$
This can be verified by defining mutually inverse homomorphisms in both directions using universal properties.
One checks that the reflection equation relations (\ref{rels2}) hold in the GWA,
and that the GWA relations (\ref{gwa_relns}) hold in $\A$.
\prop{$\A$ is a noetherian domain of GK dimension $4$.}
\begin{pf}
In \cite[Proposition 3.1]{DL}, 
polynormal sequences and Gr\"obner
basis techniques are used to show that $\A_q(\Mat_n)$ is a noetherian domain for all $n$.
Theorem \ref{thm:gwa_model} and Corollary \ref{cor:gwa_domain}
give an alternative way to see this for $\A=\A_q(\Mat_2)$.

It is also observed in \cite{DL} that the Hilbert series of $\A_q(\Mat_n)$
can be determined using \cite[(7.37)]{MaBook}.
One may deduce from the Hilbert series that
the GK dimension of $\A_q(\Mat_n)$ is $n^2$.
Theorem \ref{thm:gwa_gkdim} gives an alternative way to see this for $\A=\A_q(\Mat_2)$,
since $\sigma$ is locally algebraic.
\done
\end{pf}

By a change of variables in $k[u_{22},u_{11},z]$ we can greatly simplify
the expression of $\A$ as a GWA.
Consider the change of variables:
\begin{equation} \label{eqn:change_var}
\begin{array}{l}
u = u_{22} \\
t = u_{11}+q^{-2}u_{22}\\
d = z-q^{-2}u_{11}u_{22}
\end{array}
\end{equation}
Now we have
\begin{equation} \label{eqn:simple_gwa}
\A \cong k[u,t,d][x,y;\sigma,z],
\end{equation}
where
$z=d+q^{-2}tu-q^{-4}u^2$ and
$$\begin{array}{lll}
\sigma(u)&=&q^2u\\
\sigma(t)&=&t\\
\sigma(d)&=&d.
\end{array}$$
The special elements $t$ and $d$ of $\A$ are, up to a scalar multiple,
the quantum trace and quantum determinant explored in \cite{Ma}.

Since $q$ is not a root of unity, $\sigma$ has infinite order.
We may therefore apply Proposition \ref{prop:gwa_center} to determine the center of $\A$:
\prop{\label{prop:aq2_center}$Z(\A)=k[t,d]$.}
This was also computed in \cite{KSk}, and a complete description of the center of $\A_q(\Mat_n)$
for arbitrary $n$ is given in \cite{JW}.

Using the fact that $q$ is not a root of unity,
the elements
$$\sigma^m(z)=d+q^{2m-2}tu-q^{4m-4}u^2$$
of $k[u,t,d]$,
for $m\in\Z$, are pairwise coprime.
This allows us to get at the normal elements of $\A$,
which gives us a handle on its automorphism group:
\thm{\label{thm:aq2_aut}
The automorphism group of $\A$ is isomorphic to $(k^\times)^2$,
with $(\alpha,\gamma)\in(k^\times)^2$
corresponding to the automorphism given by
$$\begin{array}{llcll}
u_{11} & u_{12} &  &  \alpha \, u_{11} & \frac{\alpha}{\gamma}\,u_{12} \\
       &        & \mapsto &         &   \\
u_{21} & u_{22} &  &  \alpha\gamma \, u_{21} & \alpha\,u_{22} \, .
\end{array}$$
}
\begin{pf}
Let $\psi:\A\rightarrow\A$ be an automorphism.
By Proposition \ref{prop:gwa_normalelements2},
the nonzero normal elements of $\A$ are the $\sigma$-eigenvectors in $k[u,t,d]$.
That is, they all have the form $u^if(t,d)$ for some polynomial $f(t,d)$ and some $i\in\N$.
Since $k[u,t,d]$ is the linear span of such elements, it is preserved by $\psi$.
Since $u$ is normal, $\psi(u)=u^if(t,d)$ for some $i$ and $f$,
and similarly $\psi^{-1}(u)=u^jg(t,d)$ for some $j$ and $g$.
Note that $k[t,d]$, being the center of $\A$, is also preserved by $\psi$.
Therefore
$u=\psi(\psi^{-1}(u))=u^{ij}f^j\psi(g)$
implies that $i=1$ and $f$ is a unit.
So $\psi(u)=\alpha u$ for some $\alpha\in k^\times$.

Observe that $\psi(x)u=\alpha^{-1}\psi(xu)=\alpha^{-1}q^2\psi(ux)=q^2u\psi(x)$.
Any $a\in\A$ with the property that $au=q^2ua$ is a sum of homogeneous such $a$'s,
and a homogeneous such $a$ is $bv_m$ for some $b\in k[u,t,d]$ and some $m\in \Z$ such that
$$ q^2 u b v_m = b v_m u = q^{2m} b u v_m. $$
This equation requires that either $b=0$ or $m=1$.
Therefore $\psi(x)=bx$ for some nonzero $b\in k[u,t,d]$.
The same argument applies to $\psi^{-1}$, and it is easy to deduce from this that $b$ must be a unit,
i.e. $b\in k^\times$.
Similarly, using the fact that $\psi(y)u=q^{-2}u\psi(y)$, we get that
$\psi(y)=cy$ for some $c\in k^\times$.

For any $m>0$, we have
\begin{align*}
b^{m-1}\psi(\sigma^m(z))x^{m-1}
&= \psi(\sigma^m(z)x^{m-1}) \\
&= \psi(x^my) \\
&= b^mcx^my \\
&= b^mc\sigma^m(z)x^{m-1}.
\end{align*}
It follows that $\psi(\sigma^m(z))=bc\sigma^m(z)$ for all $m>0$.
Considering that $\sigma^m(z)=d+q^{2m-2}tu-q^{4m-4}u^2$,
the linear span of $\{\sigma(z),\sigma^2(z),\sigma^3(z)\}$, for instance,
contains $\{d,tu,u^2\}$.
So $bc\,u^2=\psi(u^2)=\alpha^2u^2$, i.e. $bc=\alpha^2$.
And $bc\,tu=\psi(tu)=\alpha\psi(t)u$, so $\psi(t)=\alpha t$.
And $bc\,d=\psi(d)$, so $\psi(d)=\alpha^2  d$.
Letting $\gamma=b\alpha^{-1}$,
so that $\psi(x)=(\alpha\gamma)x$ and $\psi(y)=(\alpha\gamma^{-1})y$,
we see that $\psi$ is the automorphism corresponding to $(\alpha,\gamma)$ in the theorem statement.
One easily checks that there is such an automorphism for every $(\alpha,\gamma)\in (k^\times)^2$,
and that composition of automorphisms corresponds to multiplication in $(k^\times)^2$.
\end{pf}
\done

\subsection{Finite Dimensional Simple Modules}
\label{sec:aq2:simples}

The finite dimensional simple modules over $\A$ come in two types:
the ones annihilated by $u_{22}$ and the ones on which $u_{22}$ acts invertibly.
This observation follows from the fact that since $u_{22}$ is normal, its annihilator in any $\A$-module is a submodule.
The former are modules over $\A/\langle u_{22}\rangle$, a three-variable polynomial ring.
The latter are addressed by Theorem \ref{thm:gwa:reps} given the GWA structure (\ref{eqn:simple_gwa});
we proceed to apply the theorem and state a classification.

\def\m{\mathfrak{m}}
Assume that $k$ is algebraically closed. 
Let $R$ denote the coefficient ring $k[u,t,d]$ of $\A$ as a GWA.
Maximal ideals of $R$ take the form $\m(u_0,t_0,d_0):=\gen{u-u_0,t-t_0,d-d_0}$
for some scalars $u_0,t_0,d_0\in k$.
They get moved by $\sigma^n$ to $\m(q^{-2n}u_0,t_0,d_0)$ for $n\in\Z$,
so $\m(u_0,t_0,d_0)$ has infinite $\sigma$-orbit if and only if $u_0\neq 0$.
Therefore a finite dimensional simple left $\A$-module contains
a simple $R$-submodule with annihilator having infinite $\sigma$-orbit
if and only if $u=u_{22}$ acts nontrivially.
Theorem \ref{thm:gwa:reps} requires us to consider the condition $\sigma^{-n+1}(z),\sigma^{n'}(z)\in \m(u_0,t_0,d_0)$
where $n,n'>0$.
Since
\begin{align*}
&\sigma^{-n+1}(z) = d+q^{-2n}tu-q^{-4n}u^2\\
&\sigma^{n'}(z) = d+q^{2n'-2}tu - q^{4n'-4}u^2,
\end{align*}
a straightforward calculation shows that, as long as $u_0\neq 0$,
one has $\sigma^{-n+1}(z),\sigma^{n'}(z)\in \m(u_0,t_0,d_0)$ if and only if
\begin{align} \label{eqn:aq2:reps1}
d_0 = -q^{2(n'-n-1)} u_0^2
\hspace{1cm}
t_0=(q^{-2n}+q^{2(n'-1)})u_0.
\end{align}

Define for $u_0\in k^\times$ and $t_0,d_0\in k$ the left $\A$-module $M(u_0,d_0,t_0) := \A/(\A\m(u_0,t_0,d_0))$,
and let $e_i$ denote the image of $v_i$ in it for all $i\in\Z$.
Let $N(u_0,t_0,d_0)$ be the submodule
$ \bigoplus_{i\leq -n'}Re_i \oplus \bigoplus_{i\geq n} Re_i, $
where $n>0$ is chosen to be minimal such that $d_0+q^{-2n}t_0u_0-q^{-4n}u_0^2=0$ (or $\infty$ if this does not occur),
and $n'>0$ is chosen to be minimal such that $d_0+q^{2n'-2}t_0u_0 - q^{4n'-4}u_0^2=0$ (or $\infty$ if this does not occur).
We observed in the general setting (\ref{eqn:gwa:reps:Nm}) that this is the unique largest proper submodule of $M(u_0,d_0,t_0)$.
Define $V(u_0,t_0,d_0)$ to be the simple left $\A$-module $M(u_0,t_0,t_0)/N(u_0,t_0,d_0)$.
As an $R$-module, this is isomorphic to
$$ \bigoplus_{-n'<i<n} R/\sigma^i(\m(u_0,t_0,d_0)),$$
so it has dimension $n+n'-1$ when $n$ and $n'$ are finite.
Putting together our observations and applying Theorem \ref{thm:gwa:reps}, we have:

\def\M{\mathscr{M}}
\thm{ \label{thm:aq2:reps}
Assume that $k$ is algebraically closed. 
\vspace{-1em}
\begin{enumerate}
\item\label{itm:aq2:reps:ass1}
Let $u_0\in k^\times$ and let $t_0,d_0\in k$.
The simple left $\A$-module $V(u_0,t_0,d_0)$ is finite dimensional if and only if there are $n,n'>0$ such that
(\ref{eqn:aq2:reps1}) holds.
\item\label{itm:aq2:reps:ass2}
Let $n>0$.
Any $n$-dimensional simple left $\A$-module $V$
that is not annihilated by $u=u_{22}$
is isomorphic to 
$$V_n(u_0)\ :=\ V(u_0,\ t_0=(q^{-2n}+1)u_0 ,\ d_0=-q^{-2n} u_0^2)$$
for a unique
$u_0\in k^\times$,
namely the eigenvalue of $u_{22}$ on $\ann_V(u_{12})$.
\end{enumerate}
}

\def\U{U_q(\mathfrak{sl}_2)}
\def\Vpb{\overline{V}(n-1,+)}
\def\spacing{\hspace{0.5em}}

These simple modules are all pullbacks of simple 
$\U$-modules along homomorphisms.
Define, for each $\alpha\in k^\times$, an algebra homomorphism
$\psi_\alpha:\A\rightarrow\U$:
$$\begin{array}{lclclcl}
u_{11} &\spacing & u_{12} &  & q^{-1}(q-q^{-1})^2 \alpha EF +\alpha K^{-1} & \spacing & \alpha E \\
       & &       & \mapsto &         &  & \\
u_{21} &\spacing & u_{22} &  & q^{-1}(q-q^{-1})^2 \alpha KF & \spacing & \alpha K  \, .
\end{array}$$
Such homomorphisms can be shown to exist by checking that the relations (\ref{rels2}) hold inside $\U$
for the desired images of the $u_{ij}$.
The definition we use for $\U$
is given in \cite[I.3]{KGnB}.
For $n>0$, 
consider the $n$-dimensional simple left $\U$-module $V(n-1,+)$
defined in \cite[I.4]{KGnB}.
By using $x$ and $y$ as ``raising'' and ``lowering'' operators in the usual way,
one can easily verify that the pullback $\Vpb$ of $V(n-1,+)$
along $\psi_\alpha$ is a simple $\A$-module.
Identifying $\ann_{\Vpb}(u_{12})$ as ``$\overline{m}_0$''
from \cite[I.4]{KGnB}, which has a $u$-eigenvalue of $\alpha q^{n-1}$,
we conclude that $\Vpb \cong V_n(\alpha q^{n-1})$.
This gives:
\thm{
Assume that $k$ is algebraically closed. 
Every finite dimensional simple left $\A$-module
that is not annihilated by $u=u_{22}$
is the pullback of some simple left $\U$-module
along $\psi_\alpha$ for some $\alpha\in k^\times$.
}

\subsection{Finite Dimensional Weight Modules}
\label{sec:aq2:weight}

Keep the notation and assumptions of the previous section.
The \emph{weight} $\A$-modules are the ones that decompose into simultaneous eigenspaces for the actions of $u$, $t$, and $d$;
this is what it means to be semisimple over $R=k[u,t,d]$ when $k$ is algebraically closed.
In this section, we simply apply Theorem \ref{thm:gwa:reps:sepStrings2} to $\A$.

We observed in the previous section that the only maximal ideals $\m(u_0,t_0,d_0)$ of $R$ with finite $\sigma$-orbit
are ones with $u_0=0$.
Hence the \emph{\chaintype{}} finite-dimensional weight $\A$-modules
are exactly the ones on which $u$ acts as a unit.

\def\M{\mathscr{M}}

In the previous section we identified the set $\M$ defined in (\ref{eqn:gwa:reps:M}) as
\begin{align*}
\M&=\{\m\in\maxspec R\bbel \text{$\m$ has infinite $\sigma$-orbit, $\sigma(z)\in\m$, and $\sigma^{-n+1}(z)\in\m$ for some $n>0$}\}\\
&=\{\m(u_0,\ t_0=(q^{-2n}+1)u_0 ,\ d_0=-q^{-2n} u_0^2)\bbel \text{$u_0\in k^\times$ and $n>0$}\}.
\end{align*}
We will show that statement \ref{itm:gwa:reps:sepStringsNew} of Theorem \ref{thm:gwa:reps:sepStrings2} holds for $\A$.
Let $\m=\m(u_0,(q^{-2n}+1)u_0,-q^{-2n}u_0^2)$
be an element of $\M$.
Suppose that $\sigma^{-n'+1}(z)\in\m$, where $n'>0$.
Then, using (\ref{eqn:aq2:reps1}), we have:
\begin{align}
 (q^{-2n }+1)u_0 &=
 (q^{-2n'}+1)u_0'  \label{eqn:aq2:reps2}\\
 q^{-2n}u_0^2 &= q^{-2n'}u_0'^2. \label{eqn:aq2:reps3}
\end{align}
Using (\ref{eqn:aq2:reps2}) to eliminate  $u_0'$ from (\ref{eqn:aq2:reps3}), we obtain
$$
q^{-2n} = q^{-2n'}\left(\frac{q^{-2n}+1}{q^{-2n'}+1}\right)^2,
$$
which simplifies to
$$ (q^{2n}-q^{2n'}) = q^{-2n-2n'}(q^{2n}-q^{2n'}). $$
This requires that $n=n'$. 
Therefore
Theorem \ref{thm:gwa:reps:sepStrings2} applies to $\A$ and gives:
\thm{\label{thm:aq2:reps2}
Finite-dimensional weight left $\A$-modules on which $u=u_{22}$ acts as a unit are semisimple.
}

\subsection{Prime Spectrum}\label{sec:rea:spec}

We rely on the expression of $\A$ as a GWA in (\ref{eqn:simple_gwa}):
\begin{align*}
&k[u,t,d][x,y;\sigma,z]\\
&\sigma: u\mapsto q^2u, t\mapsto t, d\mapsto d\\
&z = d+q^{-2}tu-q^{-4}u^2.
\end{align*}
We can get at all the prime ideals of $\A$ by considering various quotients and localizations.
Let us begin by laying out notation for the algebras to be considered:
\begin{itemize}
\item
$\A/\gen{u}$ is simply a polynomial ring,
$$ \A/\gen{u} \cong k[u_{11},u_{12},u_{21}]. $$
A glance at the reflection equation relations (\ref{rels2}) is enough to see this.
\item
Let $\A_u$ denote the localization of $\A$ at the set of powers of $u$,
a denominator set because $u$ is normal.
By Proposition \ref{prop:gwa_loc}, this is $k[u^\pm,t,d][x,y;\sigma,z]$.
By Proposition \ref{prop:gwa_quot}, $\A_{u}/\gen{t,d} = k[u^\pm][x,y;\sigma,z]$.
In this quotient, $z$ is a unit: $z=-q^{-4}u^2$.
Hence, by Proposition \ref{prop:gwa_laurent2},
$$
\A_{u}/\gen{t,d} = k[u^\pm][x^\pm;\sigma].
$$
\item
Let $\A_{ud}$ denote the localization of $\A_u$ at the set of powers of $d$,
a denominator set because $d$ is central.
By Proposition \ref{prop:gwa_loc}, this is $k[u^\pm,t,d^\pm][x,y;\sigma,z]$.
By Proposition \ref{prop:gwa_quot},
$$
\A_{ud}/\gen{t} = k[u^\pm,d^\pm][x,y;\sigma,z].
$$
\item
Let $\A_{ut}$ denote the localization of $\A_u$ at the set of powers of $t$,
a denominator set because $t$ is central.
By Proposition \ref{prop:gwa_loc}, this is $k[u^\pm,t^\pm,d][x,y;\sigma,z]$.
Let $\A_{utx}$ denote the localization of this at the set of powers of $x$;
by Proposition \ref{prop:gwa_xloc}, this is indeed a denominator set, and we
obtain $\A_{utx} = k[u^\pm,t^\pm,d][x^\pm;\sigma]$. Then
$$
\A_{utx}/\gen{d} = k[u^\pm,t^\pm][x^\pm;\sigma].
$$
\item
Let $\A_{utxd}$ denote the localization of $\A_{utx}$ at the set of powers of $d$:
$$
\A_{utxd} = k[u^\pm,t^\pm,d^\pm][x^\pm;\sigma].
$$
\end{itemize}
What will turn out to be missing from this list is an algebra that gives us access to those prime ideals
of $\A_{ut}$ that contain some power of $x$.
We cover this in the next section.

\begin{excludeThis}
By Proposition \ref{prop:gwa_loc}, $\A_{ut}$ is a GWA,
$$ k[u^\pm,t^\pm, d][x,y;\sigma,z]. $$
Let $R$ denote the base ring $k[u^\pm,t^\pm, d]$.
A homogeneous prime ideal of ${\A_{ut}}$
that does not contain any power of $x$
corresponds (via Theorem \ref{thm:spec_loc})
to a prime ideal of the localization at the set of powers of $x$.
We've seen in Proposition \ref{prop:gwa_xloc}
that this gives a skew Laurent algebra
$R[x^\pm;\sigma]$.
The homogeneous prime spectrum here is readily analyzed:
$$
\grspec{R[x^\pm;\sigma]}
\cong \sspec{\sigma}{R}
= \operatorname{proj}(k[t^\pm, d][u^\pm])
\cong \spec{k[t^\pm, d]}.
$$
The first homeomorphism is from Proposition \ref{prop:spec_skewlaurent}.
The middle equality is due to the fact that $q$ is not a root of unity;
here ``$\operatorname{proj}(k[t^\pm, d][u^\pm])$'' refers to prime ideals of $k[t^\pm, d][u^\pm]$ that are homogeneous
in terms of $u$-degree.
Such primes contract to primes of $k[t^\pm, d]$, and this gives the final homeomorphism
(it can be thought of as an instance of $\mathbb{P}^0_{C[x]}\cong \spec{C}$ for commutative rings $C$).
Hence we have
\begin{equation}\label{eqn:spec_aut_no_x}
\{ P\in\grspec{\A_{ut}} \bbel x^n\notin P\quad\forall n\in\N \}\ \cong\ \spec{k[t^\pm, d]}.
\end{equation}
\end{excludeThis}

\subsubsection{
Primes of $\A_{ut}$ That Contain a Power of $x$
}
\label{sec:primespwrx}

We write $\A_{ut}$ as 
$$ \A_{ut} = R[x,y;\sigma,z],$$
where $R=k[u^\pm,t^\pm,d]$.
A reminder about our notation: a subscript on a subset of a GWA indicates a certain subset of
its base ring, seen in Definition \ref{defn:gwa_idealseq}.
Define
\begin{equation}\label{eqn:specaq2_rdef} r_n = (q^{2n}+1)^2d+q^{2n}t^2 \end{equation}
for $n\in\Z$;
these elements of $R$ will help us to understand the ideal of $\A_{ut}$ generated by a power of $x$:
\prop{\label{prop:specaq2_rs}
Let $n\in\Z_{>0}$. Then
\begin{equation} \label{eqn:specaq2_rs}
\prod_{j=n-i+1}^n r_n \in\gen{x^n}_{n-i}
\end{equation}
for all $0\leq i\leq n$.
}
\begin{pf}
The induction will rely on the following observations:
\vspace{-1em}
\begin{enumerate}
\item For $n\geq 1$, $r_n\in\gen{\sigma^n(z),z}_R$.
\label{rs_obs_1}
\item Let $I$ be an ideal of a GWA $R[x,y;\sigma,z]$. For $n\geq 1$,
$ ( R^\sigma\cap I_n ) \gen{\sigma^n(z),z}_R\subseteq I_{n-1}$.
\label{rs_obs_2}
\end{enumerate}
Direct calculation verifies observation \ref{rs_obs_1},
$$
r_n = \frac{q^{2n+2}}{q^{2n}-1}tu^{-1}(\sigma^n(z)-z) +
\frac{q^{2n}+1}{q^{2n}-1}(q^{4n}z-\sigma^n(z)),
$$
and observation \ref{rs_obs_2} follows from Proposition \ref{prop:gwa_homideals2}.
The $i=0$ case, $1\in\gen{x^n}_n$, is trivial.
Assume that $0\leq i < n$ and that (\ref{eqn:specaq2_rs}) holds for $i$.
Then $a:=r_nr_{n-1}\cdots r_{n-(i-1)}\in\gen{x^n}_{n-i}$.
By observation \ref{rs_obs_2},
$a\gen{\sigma^{n-i}(z),z}_R\subseteq \gen{x^n}_{n-(i+1)}$.
Hence, by observation \ref{rs_obs_1}, $ar_{n-i}\in \gen{x^n}_{n-(i+1)}$,
proving (\ref{eqn:specaq2_rs}) for $i+1$.
\end{pf}
\done

\prop{\label{prop:specaq2_rs2}
Assume that $n\geq 1$ and $P\in \spec{\A_{ut}}$.
If $x^n\in P$ and $x^{n-1}\notin P$,
then $r_n\in P$.
}
\begin{pf}
From the $i=n$ case of Proposition \ref{prop:specaq2_rs},
$$r_1r_2\cdots r_n\in P.$$
Since this is a product of central elements in $\A_{ut}$,
we conclude that that $r_{n'}\in P$ for some $n'$.
In particular, $r_{n'}\in P_{n-1}$.
Applying Proposition \ref{prop:specaq2_rs}
with $i=1$, we also have $r_n\in P_{n-1}$.
Since $t$ is a unit, and since $q$ is not a root of unity,
it is clear from
(\ref{eqn:specaq2_rdef})
that $1\in\gen{r_n, r_{n'}}_R$ if $n\neq n'$.
We assumed that $x^{n-1}\notin P$,
so $n'=n$.
\end{pf}
\done

So when considering homogeneous prime ideals $P$ of $\A_{ut}$
that contain a power of $x$, we can 
eliminate a variable by
factoring out the ideal generated by
one of the $r_i$.
Namely, we may factor out $\gen{r_n}$
if $n\geq 1$ is taken to be minimal such that $x^n\in P$,
and we may then consider $P$ as a prime ideal of $A_{(n)}:=\A_{ut}/\gen{r_n}$.
Using Proposition \ref{prop:gwa_quot},
this algebra is isomorphic to
$$
k[u^\pm, t^\pm][x,y;\sigma,z_n],
$$
where
\begin{equation}\label{eqn:specaq2_zndef}
z_n = \frac{-q^{2n}}{(q^{2n}+1)^2}t^2+q^{-2}ut-q^{-4}u^2.
\end{equation}
Let $R_{(n)}$ denote $k[u^\pm,t^\pm]$, thought of as $R/\gen{r_n}_R$.
The ideal generated by $x^n$ can be pinned down completely
in $A_{(n)}$.
We again start by defining some special elements of the base ring that will help us break things down.
Make the following definitions for $n\in\Z$:
\newcommand{\snj}[2]{s^{#1}_{#2}}
\newcommand{\enj}[2]{e^{#1}_{#2}}
\newcommand{\pinm}[2]{\pi^{#1}_{#2}}
\newcommand{\Jnm}[2]{\mathcal{J}^{#1}_{#2}}
\begin{equation}
\begin{array}{ll}
\snj{n}{j} = u - \dfrac{q^{2j}}{q^{2n}+1}t
&\text{\ for\ } j\in\Z. \\[1em]
\begin{cases}
\Jnm{n}{ m} &= \{ j\in\Z \bbel 1\leq j\leq n-m \}\\
\Jnm{n}{-m} &= \{ j\in\Z \bbel m+1\leq j\leq n \}
\end{cases} & \text{\ for\ } m\geq 0 \\[2em]
\pinm{n}{m} = \prod_{j\in \Jnm{n}{m}} \snj{n}{j} & \text{\ for\ } m\in\Z.
\end{array}
\label{eqn:specaq2_Jm}
\end{equation}
Here is a way to visually organize these definitions for the example $n=3$:
\begin{center}
\newcommand{\thing}[1]{$\snj{3}{#1}$}
\newcommand{\midthing}{$\cdot$}
\begin{tikzpicture}[scale=0.5]
  \draw
    (0,0) node {\thing{2}}
    (1,1) node {\thing{2}}
    (-1,-1) node {\thing{2}}
    (-1,1) node {\thing{1}}
    (0,2) node {\thing{1}}
    (-2,0) node {\thing{1}}
    (1,-1) node {\thing{3}}
    (2,0) node {\thing{3}}
    (0,-2) node {\thing{3}}
    (0,1) node {\midthing{}}
    (1,0) node {\midthing{}}
    (-1,0) node {\midthing{}}
    (0,-1) node {\midthing{}}
    (3,2) node {$=$}    (4,2) node {$\pinm{3}{2}$}
    (3,1) node {$=$}    (4,1) node {$\pinm{3}{1}$}
    (3,0) node {$=$}    (4,0) node {$\pinm{3}{0}$}
    (3,-1) node {$=$}   (4,-1) node {$\pinm{3}{-1}$}
    (3,-2) node {$=$}   (4,-2) node {$\pinm{3}{-2}$}
    (0,3) node {$1$} (3,3) node {$=$} (4,3) node {$\pinm{3}{3}$}
    (0,-3) node {$1$} (3,-3) node {$=$} (4,-3) node {$\pinm{3}{-3}$};
\end{tikzpicture}
\end{center}

Observe that $\sigma(z_n)=-\snj{n}{n}\snj{n}{0}$ and that
\begin{equation}\sigma^{-1}(\snj{n}{j})=q^{-2}\snj{n}{j+1}, \label{eqn:s_obs_4} \end{equation}
so that
\begin{equation}\sigma^i(z_n) = -q^{4i-4}\snj{n}{n-i+1}\snj{n}{1-i} \label{eqn:specaq2_sigma_z} \end{equation}
for $n,i\in\Z$. Finally, observe that the $\snj{n}{j}$ are pairwise coprime over various $j$, since $q$ is not a root of unity.

For the next results, we abstract this situation.

\prop{\label{prop:specaq2_xnm}
Fix $n\geq 1$.
Consider an arbitrary
GWA $A = R[x,y;\sigma,z]$
and sequence $(s_j)_{j\in\Z}$ of elements of $R$
such that
\vspace{-1em}
\begin{enumerate}
\item \label{item:specaq2_xn_1}
$R$ is commutative,
\item \label{item:specaq2_xn_2}
$z$ is a unit multiple of $s_1s_{n+1}$,
\item \label{item:specaq2_xn_3}
$\sigma^{-1}(s_{j})$ is a unit multiple of $s_{j+1}$,
\item \label{item:specaq2_xn_4}
and $\gen{s_i,s_j}_R=R$ for all $i\neq j$.
\end{enumerate}
Define $\Jnm{}{m}$ as $\Jnm{n}{m}$ is defined in (\ref{eqn:specaq2_Jm}) and
let $\pi_m=\prod_{j\in \Jnm{}{m}} s_j$.
Then we have 
$$\gen{x^n}_m = \gen{\pi_m}_{R}$$
for $m\in\Z$.}
\begin{pf}
The sequence of ideals $\gen{\pi_m}_{R}$ satisfies
the conditions needed in Proposition \ref{prop:gwa_homideals2} in order for
$\bigoplus_{m\in\Z}\gen{\pi_m}_{R}v_m$ to define an ideal of $A$,
as can be checked using our assumptions \ref{item:specaq2_xn_2} and \ref{item:specaq2_xn_3}.
Since $\pi_n=1$, the latter ideal contains $x^n$.
This gives the inclusion $\gen{x^n}_m \subseteq \gen{\pi_m}_{R}$
for $m\in\Z$.
To get equality
we must show that 
\begin{equation} \label{eqn:specaq2_ss} \pi_m\in\gen{x^n}_m\end{equation}
for all $m\in\Z$.

For $m\geq n$, (\ref{eqn:specaq2_ss}) holds trivially.
Assume that (\ref{eqn:specaq2_ss}) holds for a given $m$, with $1\leq m\leq n$.
Then:
\begin{align}
\label{eqn:specaq2_ind1}
\gen{x^n}_{m-1}
&\supseteq \gen{\pi_m\sigma^m(z),\sigma^{-1}(\pi_m)z}_{R} \\
\label{eqn:specaq2_ind2}
&= \gen{  \left(\prod_{j=1}^{n-m}s_{j}\right) (s_{n-(m-1)}s_{1-m}) ,
\left( \prod_{j=2}^{n-(m-1)}s_{j} \right)(s_{n+1}s_{1})  }_{R} \\
&= \pi_{m-1} \gen{ s_{1-m} , s_{n+1} }_{R} 
\label{eqn:specaq2_ind4}
=\pi_{m-1}R.
\end{align}
Line (\ref{eqn:specaq2_ind1}) is due to the induction hypothesis and Proposition \ref{prop:gwa_homideals2}.
Line (\ref{eqn:specaq2_ind2}) uses assumptions \ref{item:specaq2_xn_2} and \ref{item:specaq2_xn_3}.
And line (\ref{eqn:specaq2_ind4}) uses assumption \ref{item:specaq2_xn_4}.
Hence, by induction, (\ref{eqn:specaq2_ss}) holds for $m\geq 0$.

Now assume that $1-n\leq m\leq 0$ and that (\ref{eqn:specaq2_ss}) holds for $m$.
We can apply a similar strategy to what was done for (\ref{eqn:specaq2_ind1})-(\ref{eqn:specaq2_ind4}):
\begin{align*}
\gen{x^n}_{m-1}&\supseteq \gen{\pi_m , \sigma^{-1}(\pi_m)}_{R} 
= \gen{  \prod_{j=-m+1}^{n}s_{j}\ ,\ \prod_{j=-m+2}^{n+1}s_{j}  }_{R} \\
&=\pi_{m-1} \gen{   s_{-m+1} , s_{n+1}  }_{R} 
= \pi_{m-1}R.
\end{align*}
Hence, by induction, (\ref{eqn:specaq2_ss}) holds for all $m\geq -n$.
In particular (the case $m=-n$), $y^n\in\gen{x^n}$.
Thus (\ref{eqn:specaq2_ss}) holds trivially for $m<-n$.
\end{pf}
\done

\cor{ \label{cor:specaq2_xnm}
In the setup of Proposition \ref{prop:specaq2_xnm},
$\gen{x^n}=\gen{y^n}$.
}
\begin{pf}
We shall make use of Proposition \ref{prop:gwa_op}
to exploit symmetries in the hypotheses of Proposition \ref{prop:specaq2_xnm}.
Let us use hats to denote our new batch of input data
to Proposition \ref{prop:specaq2_xnm}.
Consider $A$ as a GWA $R[\hat{x},\hat{y};\hat{\sigma},\hat z]$,
with elements $(\hat s_j)_{j\in\Z}$ of $R$, where
\begin{equation}\label{eqn:specaq2_xyflipdef1}\begin{array}{l@{\hspace{7mm}}l@{\hspace{7mm}}l}
\hat{x}=y&
\hat{y}=x&
\hat{z}=\sigma(z)\\
\hat{\sigma} = \sigma^{-1} &
\hat s_j = s_{n+1-j} &
\end{array}\end{equation}
This satisfies the hypotheses of Proposition \ref{prop:specaq2_xnm}.
Following along the notations needed to state the conclusion, define
\begin{equation}\label{eqn:specaq2_xyflipdef2}\begin{array}{l}
\hat\pi_m=\prod_{j\in \Jnm{}{m}}\hat s_j
\end{array}\end{equation}
and also define
\begin{equation}\label{eqn:specaq2_xyflipdef3}\begin{array}{l}
\hat{I}_m = I_{-m} \hspace{1cm} ( \text{which is\ } \{r\in R\bbel rv_{-m}\in  I\})
\end{array}\end{equation}
whenever $I\subseteq A$,
to match Definition \ref{defn:gwa_idealseq} with the new GWA structure.
Observe that 
\begin{equation}
\{n+1-j\bbel j\in \Jnm{}{m}\} = \Jnm{}{-m}\label{eqn:specaq2_Jsym1}
\end{equation}
for all $m\in \Z$,
so that $\hat{\pi}_m = \pi_{-m}$.
The conclusion of Proposition \ref{prop:specaq2_xnm} for the items with hats is then that
$ \widehat{\gen{\hat{x}^n}}_m = \gen{\hat{\pi}_m}_R $. That is,
$$
\gen{y^n}_{-m} = \gen{\pi_{-m}}_R
$$
for all $m\in\Z$.
\end{pf}
\done

In order to get at the homogeneous primes of $A_{(n)}$ that contain $x^n$,
we now seek to describe \emph{all} the homogeneous ideals of $A_{(n)}$ that contain $x^n$.
Statements of the next few results remain in a general GWA setting,
in order to continue taking advantage of the symmetry of GWA expressions.

\prop{\label{prop:specaq2_CRTbasis}
Assume the setup of Proposition \ref{prop:specaq2_xnm}.
Fix arbitrary integers $\ell_1 \leq \ell_2$.
There is an element $e_0$ of $R$ such that,
setting $e_j=\sigma^{-j}(e_0)$ for $j\in\Z$, the family $(e_j)_{j\in\Z}$ satisfies:
\vspace{-1em}
\begin{enumerate}
\item \label{itm:specaq2_CRTbasis1}
$e_j \equiv 1$ mod $s_j$ for $j\in\Z$.
\item \label{itm:specaq2_CRTbasis2}
$e_j \equiv 0$ mod $s_i$ for distinct $i,j\in\{\ell_1,\ldots,\ell_2\}$.
\item \label{itm:specaq2_CRTbasis3}
$\bar{e}_{\ell_1},\ldots,\bar{e}_{\ell_2}$ is a collection of
orthogonal idempotents that sum to $1$,
where bars denote cosets with respect to $\gen{\prod_{i=\ell_1}^{\ell_2}s_i}$.
\end{enumerate}
}
\begin{pf}
The $s_{j}$, for $j\in\Z$, are pairwise coprime as elements of $R$.
The Chinese Remainder Theorem (CRT) provides an $e_{0}\in R$
which is congruent to $1$ mod $s_{0}$ and congruent to $0$ mod
$s_{i}$ for all nonzero $i\in\{\ell_1-\ell_2,\ldots,\ell_2-\ell_1\}$.
Then for $j\in\Z$ we have that $\sigma^{-j}(e_{0})$
is congruent to $1$ mod $s_{j}$ and congruent to $0$
mod $s_{i}$ for all $i\in\{\ell_1-\ell_2+j,\ldots,\ell_2-\ell_1+j\}$ with $i\neq j$.
Setting $e_{j} = \sigma^{-j}(e_{0})$ gives us \ref{itm:specaq2_CRTbasis1} and \ref{itm:specaq2_CRTbasis2}.
Part of the CRT says that
$$\gen{\prod_{i=\ell_1}^{\ell_2}s_i} = \bigcap_{i=\ell_1}^{\ell_2}\gen{s_i},$$
and \ref{itm:specaq2_CRTbasis3} easily follows from this using \ref{itm:specaq2_CRTbasis1} and \ref{itm:specaq2_CRTbasis2}.
\end{pf}
\done

\prop{
\label{prop:specaq2_homideals0}
Assume the setup of Proposition \ref{prop:specaq2_xnm}.
Let $(e_j)_{j\in\Z}$ be as in
Proposition \ref{prop:specaq2_CRTbasis} with $\ell_1=1$ and $\ell_2=n$.
There are mutually inverse inclusion-preserving bijections
\begin{equation} \label{eqn:specaq2_corr0} \begin{array}{ccc}
\settext{ \parbox{17em}{homogeneous right $R[x;\sigma]$-submodules $I$ of $A$ containing $x^n$} }
& \leftrightarrow
& \settext{ \parbox{17em}{families $(I_{mj} \bbel m\in \Z,\ j\in \Jnm{}m)$ 
of ideals of $R$
satisfying (\ref{eqn:specaq2_cond0})
with $s_j\in I_{mj}$ for all $m,j$
} }
\\[2em]
I & \mapsto & (I_m+\gen{s_j} \bbel m\in \Z,\ j\in \Jnm{}m) \\[0.5em]
\bigoplus_{m\in\Z} \left( \gen{\pi_m} + \sum_{j\in \Jnm{}m}I_{mj}e_{j} \right)v_m & \mapsfrom & (I_{mj} \bbel m\in \Z,\ j\in \Jnm{}m),\\[1.5em]
\end{array} \end{equation}
where the condition (\ref{eqn:specaq2_cond0}) is that
\begin{equation}  \label{eqn:specaq2_cond0}
\begin{array}{l@{\hspace{0.6cm}}l@{\hspace{0.6cm}}l}
I_{-(m+1),j}\subseteq I_{-m,j}\ \forall j\in \Jnm{}{-(m+1)} &
\text{and} & 
I_{mj}\subseteq I_{m+1,j}\ \forall j\in \Jnm{}{m+1}
\end{array}
\end{equation}
for all $m\in \N$.
}
\begin{pf}
Combining Propositions \ref{prop:gwa_homideals0} and \ref{prop:specaq2_xnm},
we obtain the following correspondence:
\begin{equation} \label{eqn:specaq2_corr1} \begin{array}{ccc}
\settext{ \parbox{17em}{homogeneous right $R[x;\sigma]$-submodules $I$ of $A$ containing $x^n$} }
& \leftrightarrow
& \settext{ \parbox{17em}{sequences $(I_m)_{m\in\Z}$ of ideals of $R$
satisfying the conditions (\ref{eqn:gwa_homideals0}) of Proposition \ref{prop:gwa_homideals0}
with $\pi_m\in I_m$ for all $m$} }
\\[2em]
I & \mapsto & (I_m)_{m\in\Z} \\[0.5em]
\bigoplus_{m\in\Z} I_mv_m & \mapsfrom & (I_m)_{m\in\Z}.\\[1.5em]
\end{array} \end{equation}
The $s_{j}$, for $j\in\Z$, are pairwise coprime as elements of $R$.
So an ideal $I_m$ of $R$ containing $\pi_m$ corresponds,
via the CRT,
to a collection of ideals $(I_{mj})_{j\in \Jnm{}m}$ such
that $s_{j}\in I_{mj}$ for $j\in \Jnm{}m$.
Explicitly, the correspondence is:
\begin{equation} \label{eqn:specaq2_corr2} \begin{array}{ccc}
\settext{ \parbox{17em}{sequences $(I_m)_{m\in\Z}$ of ideals of $R$
with $\pi_m\in I_m$ for all $m$} }
& \leftrightarrow
& \settext{ \parbox{17em}{families $(I_{mj} \bbel m\in \Z,\ j\in \Jnm{}m)$ 
of ideals of $R$ with $s_{j}\in I_{mj}$ for all $m,j$} }
\\[2em]
(I_m)_{m\in\Z} & \mapsto & (I_{mj} = I_m + \gen{s_{j}}\bbel m\in\Z,\ j\in \Jnm{}m) \\[0.5em]
\left(\gen{\pi_m} + \sum_{j\in \Jnm{}m} I_{mj}e_{j}\right)_{m\in\Z} & \mapsfrom & (I_{mj} \bbel m\in \Z,\ j\in \Jnm{}m).\\[1.5em]
\end{array} \end{equation}
In order to make use of this with (\ref{eqn:specaq2_corr1}),
we need to express the condition (\ref{eqn:gwa_homideals0}) of Proposition \ref{prop:gwa_homideals0}
in terms of the $I_{mj}$.
Let $(I_m)_{m\in\Z}$  be a sequence of ideals of $R$
with $\pi_m\in I_m$ for all $m$, and let $(I_{mj} \bbel m\in \Z,\ j\in \Jnm{}m)$
be the family of ideals it corresponds to in (\ref{eqn:specaq2_corr2}).
For $m\in\N$,
\begin{align}
\label{eqn:specaq2_condline1}
I_m\subseteq I_{m+1}\ 
&\Leftrightarrow\ \gen{\pi_m}     + \sum_{j=1}^{n-m}  I_{mj}   e_j
\subseteq         \gen{\pi_{m+1}} + \sum_{j=1}^{n-m-1}I_{m+1,j}e_j \\
\label{eqn:specaq2_condline2}
&\Rightarrow\     \gen{s_i}       + I_{mi}
\subseteq         \gen{s_i}       + I_{m+1,i}
\hspace{3em} \forall i\in \Jnm{}{m+1}\\
\label{eqn:specaq2_condline3}
&\Rightarrow\                       I_{mi}
\subseteq                           I_{m+1,i}
\hspace{3em} \forall i\in \Jnm{}{m+1}\\
\label{eqn:specaq2_condline4}
&\Rightarrow\ I_m\subseteq I_{m+1}.
\end{align}
Line (\ref{eqn:specaq2_condline2})
is obtained by adding $\gen{s_i}$
to both sides of the inclusion in line (\ref{eqn:specaq2_condline1}),
and using the properties of the $e_j$ from Proposition \ref{prop:specaq2_CRTbasis}.
Line (\ref{eqn:specaq2_condline3}) is due to the fact that
$s_i\in I_{m+1,i}$.
Line (\ref{eqn:specaq2_condline4})
can be seen by looking at (\ref{eqn:specaq2_condline1})
and noting that $e_{n-m}\in\gen{\pi_{m+1}}$ because $e_{n-m}$
vanishes mod $s_{j}$ for $j\in \Jnm{}{m+1}$.
For similar reasons we also have, for $m\in\N$,
\begin{align}
\nonumber I_{-(m+1)}&\sigma^{-m}(z)\subseteq I_{-m}\\
&\Leftrightarrow\   \left( \gen{\pi_{-(m+1)}} + \sum_{j=m+2}^{n}  I_{-(m+1),j} e_j \right)s_{n+m+1}s_{m+1}
\subseteq                  \gen{\pi_{-m}}     + \sum_{j=m+1}^{n}  I_{-m,j}     e_j \\
\label{eqn:specaq2_condline5}
&\Leftrightarrow\ s_{n+m+1} \gen{\pi_{-m}}    + \sum_{j=m+2}^{n}  I_{-(m+1),j} s_{n+m+1}s_{m+1} e_j
\subseteq                   \gen{\pi_{-m}}    + \sum_{j=m+1}^{n}  I_{-m,j}     e_j \\
\label{eqn:specaq2_condline6}
&\Rightarrow\ \gen{s_i} + I_{-(m+1),i} s_{n+m+1}s_{m+1}
\subseteq     \gen{s_i} + I_{-m,i}
\hspace{3em } \forall i\in \Jnm{}{-(m+1)}\\
\label{eqn:specaq2_condline7}
&\Rightarrow\ I_{-(m+1),i}
\subseteq     I_{-m,i}
\hspace{3em } \forall i\in \Jnm{}{-(m+1)}\\
\label{eqn:specaq2_condline8}
&\Rightarrow\ I_{-(m+1)}\sigma^{-m}(z)\subseteq I_{-m}.
\end{align}
The only subtlety this time is that line
(\ref{eqn:specaq2_condline7})
relies on the fact that $s_{n+m+1}$ and $s_{m+1}$ are units modulo $s_i$ for all $i\in \Jnm{}{-(m+1)}$.

We conclude that
the condition (\ref{eqn:gwa_homideals0}) of Proposition \ref{prop:gwa_homideals0} holds for $(I_m)_{m\in\Z}$
if and only if the condition (\ref{eqn:specaq2_cond0}) holds for 
$(I_{mj} \bbel m\in \Z,\ j\in \Jnm{}m)$.
Combining this fact with the correspondences 
(\ref{eqn:specaq2_corr1})
and
(\ref{eqn:specaq2_corr2})
yields the desired correspondence
(\ref{eqn:specaq2_corr0}).
\end{pf}
\done

\cor{\label{cor:specaq2_homideals1}
Assume the setup of Propositions \ref{prop:specaq2_xnm}
and \ref{prop:specaq2_homideals0}.
There are mutually inverse inclusion-preserving bijections
\begin{equation} \label{eqn:specaq2_corr6} \begin{array}{ccc}
\settext{ \parbox{18em}{homogeneous right $R[y;\sigma^{-1}]$-submodules $I$ of $A$ containing $y^n$} }
& \leftrightarrow
& \settext{ \parbox{17em}{families $(I_{mj} \bbel m\in \Z,\ j\in \Jnm{}m)$ 
of ideals of $R$
satisfying (\ref{eqn:specaq2_cond1})
with $s_j\in I_{mj}$ for all $m,j$
} }
\\[2em]
I & \mapsto & (I_m+\gen{s_j} \bbel m\in \Z,\ j\in \Jnm{}m) \\[0.5em]
\bigoplus_{m\in\Z} \left( \gen{\pi_m} + \sum_{j\in \Jnm{}m}I_{mj}e_{j} \right)v_m & \mapsfrom & (I_{mj} \bbel m\in \Z,\ j\in \Jnm{}m),\\[1.5em]
\end{array} \end{equation}
where the condition (\ref{eqn:specaq2_cond1}) is that
\begin{equation}  \label{eqn:specaq2_cond1}
\begin{array}{l@{\hspace{0.6cm}}l@{\hspace{0.6cm}}l}
 I_{-(m+1),j}\supseteq  I_{-m,j}\ \forall j\in \Jnm{}{-(m+1)} &
\text{and} & 
I_{mj}\supseteq I_{m+1,j}\ \forall j\in \Jnm{}{m+1}
\end{array}
\end{equation}
for all $m\in \N$.
}
\begin{pf}
We shall apply Proposition \ref{prop:specaq2_homideals0}
while viewing $A$ as a GWA with the alternative
GWA structure  $R[y,x;\sigma,^{-1},\sigma(z)]$.
Make the definitions
(\ref{eqn:specaq2_xyflipdef1})-(\ref{eqn:specaq2_xyflipdef3}),
and also define
\begin{equation}\label{eqn:specaq2_xyflipdef4}
\hat e_j = e_{n+1-j}.
\end{equation}
This data satisfies the hypotheses of Proposition \ref{prop:specaq2_homideals0},
and allows us to conclude that there is a correspondence
\begin{equation} \label{eqn:specaq2_corr5} \begin{array}{ccc}
\settext{ \parbox{18em}{homogeneous right $R[y;\sigma^{-1}]$-submodules $I$ of $A$ containing $y^n$} }
& \leftrightarrow
& \settext{ \parbox{17em}{families $(\hat I_{mj} \bbel m\in \Z,\ j\in \Jnm{}m)$ 
of ideals of $R$
satisfying (\ref{eqn:specaq2_cond0hat})
with $\hat s_j\in \hat I_{mj}$ for all $m,j$
} }
\\[2em]
I & \mapsto & (\hat I_m+\gen{\hat s_j} \bbel m\in \Z,\ j\in \Jnm{}m) \\[0.5em]
\bigoplus_{m\in\Z} \left( \gen{\hat \pi_m} + \sum_{j\in \Jnm{}m}\hat I_{mj}\hat e_{j} \right)v_m & \mapsfrom & (\hat I_{mj} \bbel m\in \Z,\ j\in \Jnm{}m),\\[1.5em]
\end{array} \end{equation}
where the condition (\ref{eqn:specaq2_cond0hat}) is that
\begin{equation}  \label{eqn:specaq2_cond0hat}
\begin{array}{l@{\hspace{0.6cm}}l@{\hspace{0.6cm}}l}
\hat I_{-(m+1),j}\subseteq \hat I_{-m,j}\ \forall j\in \Jnm{}{-(m+1)} &
\text{and} & 
\hat I_{mj}\subseteq \hat I_{m+1,j}\ \forall j\in \Jnm{}{m+1}
\end{array}
\end{equation}
for all $m\in \N$.
Using the observation (\ref{eqn:specaq2_Jsym1}) and
reindexing by $(m,j)\mapsto (-m,n+1-j)$, this becomes
the correspondence (\ref{eqn:specaq2_corr6}).
\end{pf}
\done

\prop{\label{prop:specaq2_homideals_allare}
Assume the setup of Proposition \ref{prop:specaq2_xnm}.
All ideals of $A$ containing $x^n$ are homogeneous.
}
\begin{pf}
Let $(e_j)_{j\in\Z}$ be as in
Proposition \ref{prop:specaq2_CRTbasis} with $\ell_1=-n+2$ and $\ell_2=2n-1$.
Let $I$ be any ideal of $A$ containing $x^n$,
and let $\sum_{m\in\Z}a_mv_m\in I$ be an arbitrary element.
Then since $\gen{x^n}=\gen{y^n}$, from Corollary \ref{cor:specaq2_xnm},
we have $v_m\in I$ for $m\geq n$ and for $m\leq -n$,
and the problem is reduced to considering 
$$\sum_{m=-n+1}^{n-1}a_mv_m\in I$$
and needing to show that $a_mv_m\in I$ for $m\in\{-n+1,\ldots,n-1\}$.
Consider any $j,j'\in \{1,\ldots, n\}$.
Multiplying on the left by $e_j$ and on the right by $e_{j'}$ yields
\begin{align*}
I &\ni \sum_{m=-n+1}^{n-1} e_j a_m v_m e_{j'} 
=     \sum_{m=-n+1}^{n-1} a_m e_j \sigma^{m}(e_{j'}) v_m 
=     \sum_{m=-n+1}^{n-1} a_m e_j e_{j'-m} v_m .
\end{align*}
When $-n+1\leq m \leq n-1$ and $1\leq j'\leq n$, we have $-n+2\leq j'-m \leq 2n-1$.
So, since $\pi_0\in I$ (due to Proposition \ref{prop:specaq2_xnm}),
the product $e_j e_{j'-m}$  that appears above vanishes mod $I$ unless $j'-m=j$,
in which case it is congruent to $e_j$ mod $I$.
Thus we have
$$ a_{j'-j}e_jv_{j'-j}\in I $$
for all $j,j'\in \{1,\ldots, n\}$.
When $j\in \Jnm{}m$, we have $j,m+j\in \{1,\ldots,n\}$,
so this shows that $a_me_jv_m \in I$, for all $m\in\{-n+1,\ldots,n-1\}$ and $j\in \Jnm{}m$,
and in  particular that
$$
\sum_{j\in \Jnm{}m} a_m e_j v_m \in I.
$$
Fix an $m\in\{-n+1,\ldots,n-1\}$.
Since $\pi_{m} \in I_m$ (due to Proposition \ref{prop:specaq2_xnm}),
$\sum_{j\in \Jnm{}m}e_j \equiv 1$ mod $I_m$.
Hence $a_mv_m\in I$.
\end{pf}
\done

\cor{\label{cor:specaq2_homideals2}
Assume the setup of Propositions \ref{prop:specaq2_xnm}
and \ref{prop:specaq2_homideals0}.
There are mutually inverse inclusion-preserving bijections
\begin{equation} \label{eqn:specaq2_corr7} \begin{array}{ccc}
\settext{ \parbox{12em}{ideals $I$ of $A$ containing $x^n$} }
& \leftrightarrow
& \settext{ \parbox{17em}{families $(I_{mj} \bbel m\in \Z,\ j\in \Jnm{}m)$ 
of ideals of $R$
satisfying (\ref{eqn:specaq2_cond2})
with $s_j\in I_{mj}$ for all $m,j$
} }
\\[2em]
I & \mapsto & (I_m+\gen{s_j} \bbel m\in \Z,\ j\in \Jnm{}m) \\[0.5em]
\bigoplus_{m\in\Z} \left( \gen{\pi_m} + \sum_{j\in \Jnm{}m}I_{mj}e_{j} \right)v_m & \mapsfrom & (I_{mj} \bbel m\in \Z,\ j\in \Jnm{}m),\\[1.5em]
\end{array} \end{equation}
where the condition (\ref{eqn:specaq2_cond2}) is that
\begin{equation}  \label{eqn:specaq2_cond2}
\begin{array}{l@{\hspace{0.6cm}}l@{\hspace{0.6cm}}l}
 I_{-(m+1),j}=  I_{-m,j}\ \forall j\in \Jnm{}{-(m+1)}, &
 & I_{mj}= I_{m+1,j}\ \forall j\in \Jnm{}{m+1}, \\
 \sigma(I_{-(m+1),j})=  I_{-m,j-1}\ \forall j\in \Jnm{}{-(m+1)}, &
\text{and} & \sigma(I_{m,j+1})= I_{m+1,j}\ \forall j\in \Jnm{}{m+1}
\end{array}
\end{equation}
for all $m\in \N$.
}
\begin{pf}
We shall deduce left-handed versions of
(\ref{eqn:specaq2_corr0}) and (\ref{eqn:specaq2_corr6})
by viewing $A^\text{op}$ as a GWA $R[x,y;\sigma^{-1},\sigma(z)]$.
Recall the notation $\Inop{I}{m}$ from Definition \ref{defn:gwa_idealseq} and Remark \ref{remk:gwa_idealseq}.
Let us use hats to keep track of things in terms of this GWA structure: $A^\text{op}=R[\hat{x},\hat{y};\hat{\sigma},\hat z]$;
define
\begin{equation}\label{eqn:specaq2_rlflipdef1}\begin{array}{l@{\hspace{7mm}}l@{\hspace{7mm}}l@{\hspace{7mm}}l}
\hat{x}=x&
\hat{y}=y&
\hat{z}=\sigma(z)\\
\hat{\sigma} = \sigma^{-1} &
\hat s_j = s_{n+1-j} &
\hat e_j = e_{n+1-j} &
\hat\pi_m=\prod_{j\in \Jnm{}m}\hat s_j.
\end{array}\end{equation}
This data satisfies the hypotheses of Proposition
\ref{prop:specaq2_homideals0}
and Corollary
\ref{cor:specaq2_homideals1},
so we obtain correspondences
\begin{equation} \label{eqn:specaq2_corr8} \begin{array}{ccc}
\settext{ \parbox{17em}{homogeneous left $R[x;\sigma]$-submodules of $A$ containing $x^n$} }
& \leftrightarrow
& \settext{ \parbox{17em}{families $(\hat I_{mj} \bbel m\in \Z,\ j\in \Jnm{}m)$ 
of ideals of $R$
satisfying (\ref{eqn:specaq2_cond2a})
with $\hat s_j\in \hat I_{mj}$ for all $m,j$
} }
\\[2em]
\settext{ \parbox{18em}{homogeneous left $R[y;\sigma^{-1}]$-submodules of $A$ containing $y^n$} }
& \leftrightarrow
& \settext{ \parbox{17em}{families $(\hat I_{mj} \bbel m\in \Z,\ j\in \Jnm{}m)$ 
of ideals of $R$
satisfying (\ref{eqn:specaq2_cond2b})
with $\hat s_j\in \hat I_{mj}$ for all $m,j$
} }
\\[2em]
I & \mapsto & (\Inop{I}{m}+\gen{\hat s_j} \bbel m\in \Z,\ j\in \Jnm{}m) \\[0.5em]
\bigoplus_{m\in\Z} \sigma^{m}\left( \gen{\hat \pi_m} + \sum_{j\in \Jnm{}m}\hat I_{mj}\hat e_{j} \right)v_m & \mapsfrom & (\hat I_{mj} \bbel m\in \Z,\ j\in \Jnm{}m),\\[1.5em]
\end{array} \end{equation}
where the specified conditions are that
\begin{equation}  \label{eqn:specaq2_cond2a} \begin{array}{l@{\hspace{0.6cm}}l@{\hspace{0.6cm}}l}
\hat I_{-(m+1),j}\subseteq \hat I_{-m,j}\ \forall j\in \Jnm{}{-(m+1)}, &
 & 
\hat I_{mj}\subseteq \hat I_{m+1,j}\ \forall j\in \Jnm{}{m+1},
\end{array}\end{equation}
\begin{equation}  \label{eqn:specaq2_cond2b} \begin{array}{l@{\hspace{0.6cm}}l@{\hspace{0.6cm}}l}
\hat I_{-(m+1),j}\supseteq  \hat I_{-m,j}\ \forall j\in \Jnm{}{-(m+1)}, &
\text{and} & 
\hat I_{mj}\supseteq \hat I_{m+1,j}\ \forall j\in \Jnm{}{m+1}
\end{array} \end{equation}
for all $m\in \N$.
To make this useful, we transform the expression for the families of ideals
in the right hand side of (\ref{eqn:specaq2_corr8}) as follows:
$$
I_{mj} = \sigma^m(\hat I_{m,n+1-(j+m)}).
$$
The index sets $\Jnm{}m$ have symmetries that can be used to reindex sums and products after applying this transformation:
$
\Jnm{}m = \{n+1-j\bbel j\in \Jnm{}{-m}\}
= \{j-m\bbel j\in \Jnm{}{-m}\}.
$
A consequence is that $\sigma^m(\hat \pi_m)$ is a unit multiple of $\pi_m$.
One now makes the routine substitutions and reindexings in 
(\ref{eqn:specaq2_corr8})-(\ref{eqn:specaq2_cond2b})
to obtain correspondences
\begin{equation} \label{eqn:specaq2_corr9} \begin{array}{ccc}
\settext{ \parbox{17em}{homogeneous left $R[x;\sigma]$-submodules of $A$ containing $x^n$} }
& \leftrightarrow
& \settext{ \parbox{17em}{families $( I_{mj} \bbel m\in \Z,\ j\in \Jnm{}m)$ 
of ideals of $R$
satisfying (\ref{eqn:specaq2_cond2c})
with $ s_j\in  I_{mj}$ for all $m,j$
} }
\\[2em]
\settext{ \parbox{18em}{homogeneous left $R[y;\sigma^{-1}]$-submodules of $A$ containing $y^n$} }
& \leftrightarrow
& \settext{ \parbox{17em}{families $(I_{mj} \bbel m\in \Z,\ j\in \Jnm{}m)$ 
of ideals of $R$
satisfying (\ref{eqn:specaq2_cond2d})
with $ s_j\in I_{mj}$ for all $m,j$
} }
\\[2em]
I & \mapsto & ( I_m+\gen{s_j} \bbel m\in \Z,\ j\in \Jnm{}m) \\[0.5em]
\bigoplus_{m\in\Z} \left( \gen{ \pi_m} + \sum_{j\in \Jnm{}m} I_{mj} e_{j} \right)v_m & \mapsfrom & ( I_{mj} \bbel m\in \Z,\ j\in \Jnm{}m),\\[1.5em]
\end{array} \end{equation}
where the specified conditions are now that
\begin{equation}  \label{eqn:specaq2_cond2c} \begin{array}{l@{\hspace{0.6cm}}l@{\hspace{0.6cm}}l}
 \sigma(I_{-(m+1),j})\subseteq I_{-m,j-1}\ \forall j\in \Jnm{}{-(m+1)}, &
 & 
\sigma(I_{m,j+1})\subseteq I_{m+1,j}\ \forall j\in \Jnm{}{m+1},
\end{array}\end{equation}
\begin{equation}  \label{eqn:specaq2_cond2d} \begin{array}{l@{\hspace{0.6cm}}l@{\hspace{0.6cm}}l}
\sigma (I_{-(m+1),j})\supseteq  \hat I_{-m,j-1}\ \forall j\in \Jnm{}{-(m+1)}, &
\text{and} & 
\sigma( I_{m,j+1} )\supseteq I_{m+1,j}\ \forall j\in \Jnm{}{m+1}
\end{array} \end{equation}
for all $m\in \N$.

Note that, by Corollary \ref{cor:specaq2_xnm}, an ideal of $A$ contains $x^n$ if and only if it contains $y^n$.
And note that, by Proposition \ref{prop:specaq2_homideals_allare}, all ideals of $A$ are homogeneous.
Hence we may combine
(\ref{eqn:specaq2_corr0}),
(\ref{eqn:specaq2_corr6}), and
(\ref{eqn:specaq2_corr9})
to obtain the correspondence
(\ref{eqn:specaq2_corr7}),
and the condition in 
(\ref{eqn:specaq2_cond2})
is just the conjunction of conditions
(\ref{eqn:specaq2_cond0}),
(\ref{eqn:specaq2_cond1}),
(\ref{eqn:specaq2_cond2c}), and
(\ref{eqn:specaq2_cond2d}).
\end{pf}
\done

We now specialize back to the algebra $A_{(n)}=\A_{ut}/\gen{r_n}$.
Corollary \ref{cor:specaq2_homideals2} applies to $A_{(n)}$ with the elements of $R_{(n)}$
defined in (\ref{eqn:specaq2_Jm}).

\newcommand{\Iddot}[1]{{\widetilde{#1}_{\Cdot[1.9]\Cdot[1.9]}}}

\prop{ \label{prop:specaq2_homideals3}
For $n\geq 1$,
there are mutually inverse inclusion-preserving bijections
\begin{equation} \label{eqn:specaq2_corr10} \begin{array}{ccc}
\settext{ \parbox{13em}{ideals $I/\gen{x^n}$ of $A_{(n)}/\gen{x^n}$} }
& \leftrightarrow
& \settext{ \parbox{7.5em}{ideals $\Iddot{I}$ of $k[t^\pm]$} } 
\\[2em]
I/\gen{x^n} & \mapsto & (I_0 + \gen{\snj{n}{1}}_{R_{(n)}})\cap k[t^\pm] \\[0.5em]
\left( \bigoplus_{m\in\Z} \left(
  \gen{\pinm{n}{m}}+\gen{\Iddot{I}}
\right)v_m \right) /\gen{x^n} & \mapsfrom & \Iddot{I}\ .\\[1.5em]
\end{array} \end{equation}
}
\begin{pf}
Let $\enj{n}{j}$ for $j\in\Z$ be as in 
Proposition \ref{prop:specaq2_CRTbasis} with $\ell_1=1$ and $\ell_2=n$.
In particular they are elements of $R_{(n)}$ such that
$\enj{n}{j}$ is congruent to $1$ mod $\snj{n}{j}$ for $j\in\Z$ and
congruent to 0 mod $\snj{n}{i}$ for all distinct $i,j\in \{1,\ldots,n\}$,
and $\sigma^{-1}(\enj{n}{j})=\enj{n}{j+1}$ for all $j\in\Z$.
For $j\in\Z$, the algebra $R_{(n)}/\gen{\snj{n}{j}}$ is isomorphic to $k[t^\pm]$, and the isomorphism
is the composite
$$k[t^\pm]\hookrightarrow R_{(n)} \twoheadrightarrow R_{(n)}/\gen{\snj{n}{j}}. $$
We obtain from this a correspondence of ideals for each $j$:
\begin{equation} \label{eqn:specaq2_corr3} \begin{array}{ccc}
\settext{ideals of $R_{(n)}$ containing $\snj{n}{j}$} & \leftrightarrow & \settext{ideals of $k[t^\pm]$} \\
J &\mapsto & \widetilde{J} = J\cap k[t^\pm]\\
\widetilde{J}+\gen{\snj{n}{j}} & \mapsfrom &\widetilde{J}.
\end{array} \end{equation}
This allows us to restate the correspondence that we obtain from
Corollary \ref{cor:specaq2_homideals2} as
\begin{equation} \label{eqn:specaq2_corr11} \begin{array}{ccc}
\settext{ \parbox{13em}{ideals $I$ of $A_{(n)}$ containing $x^n$} }
& \leftrightarrow
& \settext{ \parbox{17em}{families $(\widetilde{I}_{mj} \bbel m\in \Z,\ j\in \Jnm{n}{m})$ 
of ideals of $k[t^\pm]$
satisfying (\ref{eqn:specaq2_cond3})
} }
\\[2em]
I & \mapsto & ((I_m+\gen{\snj{n}{j}})\cap k[t^\pm] \bbel m\in \Z,\ j\in \Jnm{n}{m}) \\[0.5em]
\bigoplus_{m\in\Z} \left( \gen{\pinm{n}{m}} + \sum_{j\in \Jnm{n}{m}}(\widetilde{I}_{mj}+\gen{\snj{n}{j}})\enj{n}{j} \right)v_m & \mapsfrom 
& (\widetilde{I}_{mj} \bbel m\in \Z,\ j\in \Jnm{n}{m}),\\[1.5em]
\end{array} \end{equation}
where the condition (\ref{eqn:specaq2_cond3}) is that
\begin{equation}  \label{eqn:specaq2_cond3}
\begin{array}{l@{\hspace{0.6cm}}l@{\hspace{0.6cm}}l}
 \widetilde{I}_{-(m+1),j}=  \widetilde{I}_{-m,j}\ \forall j\in \Jnm{n}{-(m+1)}, &
 & \widetilde{I}_{mj}= \widetilde{I}_{m+1,j}\ \forall j\in \Jnm{n}{m+1}, \\
 \widetilde{I}_{-(m+1),j}=  \widetilde{I}_{-m,j-1}\ \forall j\in \Jnm{n}{-(m+1)}, &
\text{and} & \widetilde{I}_{m,j+1}= \widetilde{I}_{m+1,j}\ \forall j\in \Jnm{n}{m+1}
\end{array}
\end{equation}
for all $m\in \N$.
The $\sigma$ has disappeared from the condition (\ref{eqn:specaq2_cond2})
because $\sigma$ fixes $k[t^\pm]$.
Notice that (\ref{eqn:specaq2_cond3}) simply says
that all the ideals
in the family $(\widetilde{I}_{mj} \bbel m\in \Z,\ j\in \Jnm{n}{m})$ 
are equal. So we may as well give them all one name,
$\Iddot{I} := \widetilde{I}_{01}$.
We may also simplify the expression of the left hand side of (\ref{eqn:specaq2_corr11}):
for $m\in \Z$, we have
\begin{align}
\gen{\pinm{n}{m}} + \sum_{j\in \Jnm{n}{m}} (\Iddot{I}+\gen{\snj{n}{j}}) \enj{n}{j} 
& = \gen{\pinm{n}{m}} + \langle \snj{n}{j}\enj{n}{j} \bbel j\in \Jnm{n}{m} \rangle + \sum_{j\in \Jnm{n}{m}} \gen{\Iddot{I}}\enj{n}{j}
\nonumber \\
& = \gen{\pinm{n}{m}} + \langle \snj{n}{j}\enj{n}{j} \bbel j\in \Jnm{n}{m} \rangle + \gen{\Iddot{I}}
\label{eqn:specaq2_specAut1} \\
& = \gen{\pinm{n}{m}} +\gen{\Iddot{I}}.
\label{eqn:specaq2_specAut2}
\end{align}
Line (\ref{eqn:specaq2_specAut1})
is due to the fact that $\sum_{j\in \Jnm{n}{m}} \enj{n}{j}$
is congruent to $1$ mod $\pinm{n}{m}$.
Line (\ref{eqn:specaq2_specAut2})
is due to the fact that $\snj{n}{j}\enj{n}{j}$
is congruent to $0$ mod $\pinm{n}{m}$ for all $j\in \Jnm{n}{m}$.
Thus we obtain (\ref{eqn:specaq2_corr10}).
\end{pf}
\done

\prop{ \label{prop:specaq2_homideals_prod}
Products of ideals are preserved by the correspondence (\ref{eqn:specaq2_corr10}).
}
\begin{pf}
Let $\mathfrak{a},\mathfrak{b}$ be ideals of $k[t^\pm]$,
and let $I/\gen{x^n},J/\gen{x^n}$ be the respective corresponding ideals of $A_{(n)}/\gen{x^n}$
via (\ref{eqn:specaq2_corr10}).
We must show that the product $\mathfrak{a}\mathfrak{b}$ corresponds via
(\ref{eqn:specaq2_corr10})
to $(I/\gen{x^n})(J/\gen{x^n})
=(IJ+\gen{x^n})/\gen{x^n}$.
That is,
we must show that
$ 
( (IJ+\gen{x^n})_0 + \gen{\snj{n}{1}} )\cap k[t^\pm]
= \mathfrak{a}\mathfrak{b} $.
Using the fact that all of the ideals on the right hand side of
(\ref{eqn:specaq2_corr11})
are equal,
\begin{equation}\label{draft:081217:32:star}
\begin{array}{l}
(I_m+\gen{\snj{n}{j}})\cap k[t^\pm] = \mathfrak{a}\\[0.2em]
(J_m+\gen{\snj{n}{j}})\cap k[t^\pm] = \mathfrak{b}
\end{array}
\end{equation}
for all $m\in\Z$, $j\in\Jnm{n}{m}$.
The contraction $(IJ)_0$ of the product $IJ$ consists of
sums of products of homogeneous terms of opposite degree;
i.e. terms of the form
$$
(av_m) \cdot (bv_{-m}) = a\sigma^m(b) \zz{m}{-m} 
$$
for $m\in\Z$.
Hence $(IJ+\gen{x^n})_0+\gen{\snj{n}{1}}$ can be written as
$$
(IJ + \gen{x^n})_0 + \gen{\snj{n}{1}} = 
(IJ)_0 + \gen{\pinm{n}{0}} + \gen{\snj{n}{1}} = 
(IJ)_0 + \gen{\snj{n}{1}} = 
\sum_{m\in\Z} \zz{m}{-m} I_m \sigma^{m}(J_{-m})
+\gen{\snj{n}{1}}.
$$
Observe the following:
\claim{
If $m\in\{0,\ldots,n-1\}$, then
$\zz{m}{-m}$ is a unit mod $\gen{\snj{n}{1}}_{R_{(n)}}$.
Otherwise, it is in $\gen{\snj{n}{1}}_{R_{(n)}}$.
}{
If $m<0$, then $\zz{m}{-m}=\slinezn{m+1}{0}$ is divisible by $z_n$, which is divisible by $\snj{n}{1}$.
If $m>n-1$, then $\zz{m}{-m}=\slinezn{1}{m}$ is divisible by $\sigma^n(z_n)$, which is also divisible by $\snj{n}{1}$.
If $m=0$, then $\zz{m}{-m}=1$ is a unit mod $\snj{n}{1}$.
Finally, assume that $m\in\{1,\ldots,n-1\}$.
Then $\zz{m}{-m}=\slinezn{1}{m}$ is a unit multiple of the product
$$ \prod_{i=1}^m \snj{n}{n-i+1} \prod_{i=1}^m\snj{n}{1-i}. $$
Observe that the assumption $1\leq m\leq n-1$ precludes $\snj{n}{1}$ from being a factor in the product above.
Since the $\snj{n}{j}$ are pairwise coprime, it follows that $\zz{m}{-m}$ is a unit mod $\snj{n}{1}$.
}
This simplifies the expression above:
$$
\sum_{m\in\Z} \zz{m}{-m} I_m \sigma^{m}(J_{-m}) + \gen{\snj{n}{1}}
= \sum_{m=0}^{n-1}  I_m \sigma^{m}(J_{-m}) + \gen{\snj{n}{1}}.
$$
Now we calculate what is needed:
\begin{align}
( (IJ+\gen{x^n})_0 + \gen{\snj{n}{1}} )\cap k[t^\pm]
&= \left( \sum_{m=0}^{n-1}I_m \sigma^{m}(J_{-m}) + \gen{\snj{n}{1}} \right)\cap k[t^\pm] \nonumber\\
&= \left( \sum_{m=0}^{n-1}(I_m+\gen{\snj{n}{1}}) \sigma^{m}(J_{-m}+\gen{\snj{n}{m+1}}) + \gen{\snj{n}{1}} \right)\cap k[t^\pm] \nonumber\\
&= \left( \sum_{m=0}^{n-1}(\mathfrak{a}+\gen{\snj{n}{1}}) \sigma^{m}(\mathfrak{b}+\gen{\snj{n}{m+1}}) + \gen{\snj{n}{1}} \right)
   \cap k[t^\pm] \label{eqn:specaq2_explain0}\\
&= \left( \sum_{m=0}^{n-1}\mathfrak{a}\mathfrak{b}+\gen{\snj{n}{1}} \right) \cap k[t^\pm] \nonumber\\
&= \left( \mathfrak{a}\mathfrak{b}+\gen{\snj{n}{1}} \right) \cap k[t^\pm]  \nonumber \\
&=  \mathfrak{a}\mathfrak{b} \label{eqn:specaq2_explain1}.
\end{align}
Line (\ref{eqn:specaq2_explain0}) uses (\ref{draft:081217:32:star}),
and lines (\ref{eqn:specaq2_explain0}) and (\ref{eqn:specaq2_explain1})
both make use of the correspondence (\ref{eqn:specaq2_corr3}).
\end{pf}
\done

\cor{\label{cor:specaq2_specAut}
For $n\geq 1$,
there is a homeomorphism
$$ \spec{A_{(n)}/\gen{x^n}} \ \approx\ \spec{k[t^\pm]} $$
given by
\begin{equation} \begin{array}{ccc}
P/\gen{x^n} & \mapsto & (P_0 + \gen{\snj{n}{1}}_{R_{(n)}})\cap k[t^\pm] \\[0.5em]
\bigoplus_{m\in\Z} \left( 
\gen{\pinm{n}{m}}+\gen{\mathfrak{p}}
\right)v_m 
\ /\ \gen{x^n}
& \mapsfrom & \mathfrak{p}\ .\\[1.5em]
\end{array} \end{equation}
}
\begin{pf}
Propositions
\ref{prop:specaq2_homideals3} and
\ref{prop:specaq2_homideals_prod}.
\end{pf}
\done

\subsubsection{
The Prime Spectrum of $\A$
}

Express the algebra $\A$ as a GWA according to (\ref{eqn:simple_gwa}).
Let $X$ denote the set of positive powers of $x$.
Define $r_n \in \A$ for $n\geq 1$ as in (\ref{eqn:specaq2_rdef}).
Also define $\snj{n}{j}$, $\Jnm{n}{m}$, and $\pinm{n}{m}$ for $n\geq 1$ and $j,m\in\Z$ as in (\ref{eqn:specaq2_Jm}),
but with everything taking place in $\A$.
Define the following subsets of $\spec{\A}$:
\begin{equation*}
\begin{array}{llll}
T_1     & = & \{P\in\spec{\A} \bbel u\in P \},                                          &                       \\[0.4em]
T_2     & = & \{P\in\spec{\A} \bbel P=\gen{P\cap k[t,d]} \}, \text{\ and}               &                       \\[0.4em]
T_{3n}  & = & \{P\in\spec{\A} \bbel u,t\notin P$, $P\cap X = \{x^n, x^{n+1},\ldots\} \} & \text{for $n\geq 1$}.
\end{array}
\end{equation*}

\def\p{\mathfrak{p}}
\def\S{\mathcal{S}}
\def\G{\mathcal{G}}

\thm{\label{thm:specaq2}
The prime spectrum of $\A$ is, as a set, the disjoint union of $T_1$, $T_2$, and $T_{3n}$ for $n\geq 1$.
Each of these subsets is homeomorphic to the prime spectrum of a commutative algebra as follows:
\vspace{-1em}
\begin{itemize}
\item $\spec{k[u_{11},u_{12},u_{21}]}\approx T_1$ via $\p\mapsto \gen{u}+\gen{\p}$.
\item $\spec{ k[t,d] }\approx T_2$ via $\p\mapsto \gen{\p}$.
\item $\spec{ k[t^\pm] }\approx T_{3n}$ via $\p\mapsto 
\gen{\pinm{n}{m}v_m \bbel -n\leq m \leq n} + \gen{r_n} + \gen{\p\cap k[t]}
$, for all $n\geq 1$.
\end{itemize}
}
Our proof will make use of the localizations and quotients of $\A$
that were described in the introduction to section \ref{sec:rea:spec}.
Many of them are quantum tori, so it will help that the prime spectrum of a quantum torus is known.
\defn{A \emph{quantum torus} over a field $k$ is an iterated skew Laurent algebra
$$
k[x_1^\pm][x_2^\pm; \tau_2]\cdots[x_n^\pm;\tau_n]
$$
for some $n\in\N$ and some automorphisms $\tau_2,\ldots,\tau_n$
such that $\tau_i(x_j)$ is a nonzero scalar multiple of $x_j$
for all $i\in\{2,\ldots, n\}$ and $j\in\{1,\ldots, i-1\}$.
}
\begin{mylemma}[{\cite[Corollary 1.5b]{GL98}}]
\label{lem:spec_qmtorus}
Contraction and extension provide mutually inverse homeomorphisms
between the prime spectrum of a quantum torus and the prime spectrum of its center.
\end{mylemma}
\textbf{Proof of Theorem \ref{thm:specaq2}:}
Consider the partition of $\spec{\A}$ into subsets $S_1,\ldots, S_6$ given by
the following tree, in which branches represent mutually exclusive possibilities:
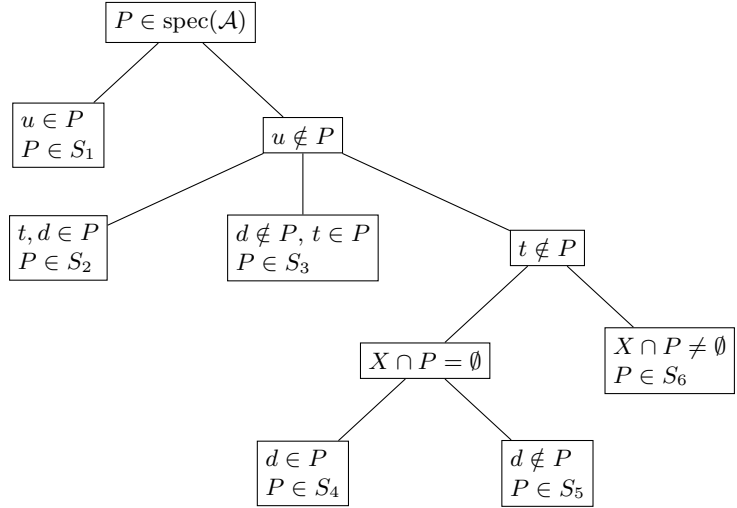
\begin{center}
\begin{tikzpicture}[sibling distance=10em,
  every node/.style = {shape=rectangle,
    draw, align=left}]]
  \node {$P\in\spec{\A}$}
    child { node {$u\in P$\\$P\in S_1$} }
    child { node {$u\notin P$}
      child { node {$t,d\in P$\\$P\in S_2$} }
      child { node {$d\notin P$, $t\in P$\\$P\in S_3$} }
      child { node {$t\notin P$}
        child { node {$X\cap P  = \emptyset$} 
          child { node {$d\in P$\\$P\in S_4$} }
          child { node {$d\notin P$\\$P\in S_5$} }
        }
        child { node {$X\cap P \neq \emptyset$\\$P\in S_{6}$} }
      }
    };
\end{tikzpicture}
\end{center}
It is easy to verify that
\begin{align*}
&S_1=T_1, \\
&S_6=\bigsqcup_{n\geq 1} T_{3n}.
\end{align*}
To establish that $\{T_1, T_2\}\cup\{T_{3n}\bbel n\geq 1\}$ is a partition of $\spec{\A}$,
we will show that
\begin{equation} \label{eqn:specaq2_partition}
S_2\cup S_3\cup S_4\cup S_5 = T_2.
\end{equation}
Let $P\in T_2$ and let $\p=P\cap k[t,d]$.
Then, using the same reasoning as in (\ref{eqn:gwa_centralspec2}),
$P_m=\p\,k[u,t,d]$ for all $m\in\Z$.
In particular, $u\notin P$, so $P\notin S_1$,
and $P_n=P_0$ for all $n\geq 1$, so $P\notin S_6$.
This establishes the inclusion $\supseteq$ of (\ref{eqn:specaq2_partition}).
We now address the reverse inclusion.

\newcommand{\stepheader}[1]{\paragraph{#1:}}

\stepheader{$S_2\subseteq T_2$}
Since $u$ is normal,
a prime ideal of $\A$ that excludes $u$ also excludes any power of $u$.
So $S_2\approx\spec{\A_u/\gen{t,d}}$.
Since $\A_u/\gen{t,d} = k[u^\pm][x^\pm;\sigma]$ is a quantum torus,
Lemma \ref{lem:spec_qmtorus} and Proposition \ref{prop:Z_skew}
give that $S_2\approx\spec{k}$.
Let $\p$ be the single point of $\spec{k}$.
It corresponds to the zero ideal of $k[u^\pm][x^\pm;\sigma]$,
which corresponds to $\gen{t,d}\triangleleft \A_u$.
Now $\gen{t,d}\triangleleft\A$ is prime
because $\A/\gen{t,d} = k[u][x,y;\sigma,z]$ and $z\notin\gen{t,d}$.
So by
Lemma \ref{lem:contract_gens} (appendix),
$\p$ corresponds to $\gen{t,d}\triangleleft\A$, and we have
$$ S_2 = \{ \gen{t,d} \} \subseteq T_2.$$

\stepheader{$S_3\subseteq T_2$}
Since $d$ is central, a prime ideal of $\A$ that excludes $d$ also excludes any power of $d$.
So $S_3\approx \spec{\A_{ud}/\gen{t}}$.
\claim{
In the algebra $\A_{ud}/\gen{t} = k[u^\pm,d^\pm][x,y;\sigma, z]$,
one has
$\gen{x^n} = \gen{1}$ for all $n\in\N$.
}{
Let $n\geq 1$.
Multiplying $x^n$ by $y$ on either side shows that
$\gen{x^n}_{n-1}$ contains $z$ and $\sigma^{n}(z)$.
Here these are $d-q^{-4}u^2$ and $d-q^{4n-4}u^2$.
Since $u$ is invertible and $q$ is not a root of unity,
this implies that $\gen{x^n}_{n-1}=\gen{1}$;
i.e. $x^{n-1}\in\gen{x^n}$.
This works for all $n\geq 1$, so we conclude by induction that $1\in \gen{x^n}$.
}
Thus all prime ideals of $k[u^\pm,d^\pm][x,y;\sigma,z]$
are disjoint from the set of powers of $x$.
Therefore by localization, using
Proposition \ref{prop:gwa_xloc} and Theorem \ref{thm:spec_loc},
$\spec{k[u^\pm, d^\pm][x,y;\sigma,z]}\approx \spec{k[u^\pm,d^\pm][x^\pm;\sigma]}$.
Lemma \ref{lem:spec_qmtorus} and Proposition \ref{prop:Z_skew}
give that $S_3\approx \spec{k[d^\pm]}$.
Let's start with a $\p\in\spec{k[d^\pm]}$
and follow it back to $S_3$:
$\p$ corresponds to its extension $\gen{\p}\triangleleft k[u^\pm,d^\pm][x^\pm;\sigma]$.
Now $\gen{\p}\triangleleft \A_{ud}/\gen{t}$ is prime
because the quotient by it is $(k[u^\pm,d^\pm]/\gen{\p})[x,y;\sigma,z]$,
a GWA over a domain with $z\neq 0$.
Hence by
Lemma \ref{lem:contract_gens},
$\p$ corresponds to $\gen{\p}\triangleleft \A_{ud}/\gen{t}$.
This in turn corresponds to $\gen{t}+\gen{\p} = \gen{t}+\gen{\p\cap k[d]}\triangleleft \A_{ud}$.
Now $\gen{t}+\gen{\p\cap k[d]}\triangleleft \A$ is prime
because the quotient by it is $(k[u,d]/\gen{\p\cap k[d]})[x,y;\sigma,z]$,
a GWA over a domain with $z\neq 0$.
Hence by
Lemma \ref{lem:contract_gens},
$\p$ corresponds to $\gen{t}+\gen{\p\cap k[d]}\triangleleft \A$.
So $$ S_3 = \{\gen{t}+\gen{\p\cap k[d]} \bbel \p\in\spec{k[d^\pm]} \} \subseteq T_2. $$

\stepheader{$S_4\subseteq T_2$}
Since $t$ is central,
$S_4\approx \spec{\A_{utx}/\gen{d}}$.
Since $\A_{utx}/\gen{d}=k[u^\pm,t^\pm][x^\pm;\sigma]$ is a quantum torus,
Lemma \ref{lem:spec_qmtorus} and Proposition \ref{prop:Z_skew}
give that $S_4\approx \spec{k[t^\pm]}$.
Let $\p\in\spec{k[t^\pm]}$,
and let us follow $\p$ back to $S_4$:
$\p$ corresponds to its extension $\gen{\p}\triangleleft \A_{utx}/\gen{d}$,
which in turn corresponds to $\gen{d}+\gen{\p} = \gen{d}+\gen{\p\cap k[t]} \triangleleft \A_{utx}$.
Now $\gen{d}+\gen{\p\cap k[t]}\triangleleft \A$ is prime because
the quotient by it is $(k[u,t]/\gen{\p\cap k[t]})[x,y;\sigma,z]$,
a GWA over a domain with $z\neq 0$.
Hence by
Lemma \ref{lem:contract_gens},
$\p$ corresponds to $\gen{d}+\gen{\p\cap k[t]}\triangleleft \A$.
So $$ S_4 = \{ \gen{d} + \gen{\p\cap k[t]} \bbel \p\in\spec{k[t^\pm]} \} \subseteq T_2. $$

\stepheader{$S_5\subseteq T_2$}
We have $S_5\approx \spec{\A_{utxd}}$.
Since $\A_{utxd} = k[u^\pm,t^\pm,d^\pm][x^\pm;\sigma]$ is a quantum torus, 
Lemma \ref{lem:spec_qmtorus} and Proposition \ref{prop:Z_skew}
give $S_5\approx k[t^\pm,d^\pm]$.
Let $\p\in \spec{k[t^\pm,d^\pm]}$, and let us follow it back to $S_5$:
$\p$ corresponds to its extension $\gen{\p} = \gen{\p\cap k[t,d]} \triangleleft \A_{utxd}$.
Now $\gen{\p\cap k[t,d]}\triangleleft \A$ is prime because
the quotient by it is $(k[u,t,d]/\gen{\p\cap k[t,d]})[x,y;\sigma, z]$,
a GWA over a domain with $z\neq 0$.
Hence by
Lemma \ref{lem:contract_gens},
$\p$ corresponds to $\gen{\p\cap k[t,d]}\triangleleft \A$.
So $$ S_5 = \{ \gen{\p\cap k[t,d]} \bbel \p\in\spec{k[t^\pm,d^\pm]} \} \subseteq T_2. $$

We have established (\ref{eqn:specaq2_partition}), proving that
$$ \spec{\A} = T_1\,\sqcup\,T_2\,\sqcup\bigsqcup_{n\geq1}T_{3n}. $$
The remainder of the proof establishes homeomorphisms of $T_1$, $T_2$, and the $T_{3n}$
to spectra of commutative algebras.

\stepheader{$T_1$}
Clearly, $T_1$ is homeomorphic to the prime spectrum of $\A/\gen{u}\cong k[u_{11},u_{12},u_{21}]$
via $\p\mapsto \gen{u}+\gen{\p}$.

\stepheader{$T_2$}
Note that the ring extension $k[u,t,d]^\sigma = k[t,d]\subseteq k[u,t,d]$
satisfies the condition (\ref{eqn:gwa_centralspec1}).
It also satisfies the condition (\ref{eqn:gwa_centralspec1.6}),
due to Proposition \ref{prop:gwa_centralspec}.
We may therefore apply Lemma \ref{lem:gwa_centralspec} to conclude that $T_2\approx \spec{k[t,d]}$,
with $\p\in\spec{k[t,d]}$ corresponding to $\gen{\p}\triangleleft\A$.

\stepheader{$T_{3n}$}
Let $n\geq 1$.
We have $T_{3n}\approx \spec{A_{(n)}/\gen{x^{n}}}$.
By
Corollary \ref{cor:specaq2_specAut},
we in turn have $\spec{A_{(n)}/\gen{x^n}} \approx \spec{k[t^\pm]}$.
Let $\p\in\spec{k[t^\pm]}$,
and let us follow it back to $T_{3n}$.
In Corollary \ref{cor:specaq2_specAut}, $\p$ corresponds to
$$\bigoplus_{m\in\Z} (\gen{\pinm{n}{m}}+\gen{\p})v_m ,$$
which is
\begin{equation} \label{eqn:specaq2_specAut3}
\gen{\pinm{n}{m}v_m \bbel -n\leq m \leq n} + \gen{\p}
\ \triangleleft\  A_{(n)} = \A_{ut}/\gen{r_n}.
\end{equation}
Applying Lemma \ref{lem:contract_gens} is not as trivial in this situation,
so we will check the needed hypotheses carefully.
Write $\p\cap k[t]$ as $\gen{p}\in \spec{k[t]}$,
where $p$ is either zero or it is some irreducible polynomial in $t$ that is not divisible by $t$.
Let $A=\A/\gen{r_n}$ and let $R=k[u,t]$.
Then $A$ is a GWA $ R[x,y;\sigma,z_n]$,
with the $z_n$ given in (\ref{eqn:specaq2_zndef}).
Note that the extension of the ideal $\gen{r_n}\triangleleft \A$ to
$\A_{ut}$ is $\gen{r_n}\triangleleft \A_{ut}$.
So by Proposition \ref{prop:loc_quot}, $\A_{ut}/\gen{r_n}$
is the localization of $A$
at $\S:=\{u^it^j\bbel i,j\in\N\}\subseteq A$.
Define
$$\G:=\{\pinm{n}{m}v_m \bbel -n\leq m \leq n\}\cup\{p\}\ \subseteq\ A.$$
Let $P=\bigoplus_{m\in\Z} \gen{\pinm{n}{m},p}v_m \ \triangleleft\  A.$
Check that this is an ideal of $\A$ by using
(\ref{eqn:s_obs_4}) and (\ref{eqn:specaq2_sigma_z})
to verify that the conditions of Proposition  \ref{prop:gwa_homideals2} are met.
$P$ is generated by $\G$ as a right ideal of $A$;
this takes care of
hypothesis \ref{itm:contract_gens_2sided} of Lemma \ref{lem:contract_gens}.
Hypothesis \ref{itm:contract_gens_prime}
requires some work to verify. First we need:
\claim{
For $m\in\Z$,
\begin{equation}\label{eqn:UFDpeup}
\gen{\pinm{n}{m},p}_R = \bigcap_{j\in \Jnm{n}{m}} \gen{\snj{n}{j},p}_R.
\end{equation}
}{
Assume that $-n< m< n$, otherwise there is nothing to prove
(take an empty intersection to be $R$).
If $p=0$, then (\ref{eqn:UFDpeup})
follows from the fact that $R$ is a UFD and $\pinm{n}{m}$ is a product of the
non-associate irreducibles  $\snj{n}{j}\in R$.
Assume that $p\neq 0$, so $p$ is an irreducible polynomial in $t$ that is not divisible by $t$.
For convenience of notation, let $s_1,\ldots, s_r$ be the elements of $\{\snj{n}{j}\bbel j\in \Jnm{n}{m}\}$.
We will show that
\begin{equation}\label{eqn:UFDpeup2}
\gen{s_1s_2\cdots s_r,p} = \gen{s_1,p}\cap\gen{s_2\cdots s_r,p}
\end{equation}
and then (\ref{eqn:UFDpeup})
will follow by repeating the same principle with induction.

The inclusion $\subseteq$ of (\ref{eqn:UFDpeup2}) is obvious. For $\supseteq$,
suppose that
$\alpha s_1 + \gamma p = \beta s_2\cdots s_r + \delta p$, where $\alpha,\beta,\gamma,\delta\in R$.
Then $(\gamma-\delta)p \in \gen{s_1,s_2\cdots s_r}$.
We can see that $p$ is regular mod $\gen{s_1,s_2\cdots s_r}$
by using an isomorphism $R/\gen{s_1} \cong k[t]$ that fixes $t$:
the image of $p$ under $R\rightarrow R/\gen{s_1}\cong k[t]$ is itself,
and the image of $\gen{s_2\cdots s_r}$ is $\gen{t^{r-1}}$
(since $q$ is not a root of unity).
Hence we have $(\gamma-\delta)\in\gen{s_1,s_2 \cdots s_r}$.
Write it as $(\gamma-\delta) = \epsilon s_1+\zeta s_2\cdots s_r$,
for some $\epsilon,\zeta\in R$.
Then
\begin{align*}
(\alpha+\epsilon p)s_1
&= (\alpha s_1) + (\epsilon s_1) p
= (\beta s_2\cdots s_r + \delta p -\gamma p) + (\gamma - \delta -\zeta s_2\cdots s_r)p\\
&= (\beta - \zeta p ) s_2\cdots s_r.
\end{align*}
Since $s_1,s_2,\ldots,s_r$ are non-associate irreducibles,
it follows that $(\beta-\zeta p)=\eta s_1$ for some $\eta\in R$.
Finally,
\begin{align*}
\alpha s_1+\gamma p
&= \beta s_2 \cdots s_r + \delta p
= (\eta s_1 + \zeta p)s_2\cdots s_r +\delta p\\
&= \eta s_1s_2 \cdots s_r + (\zeta s_2 \cdots s_r + \delta)p,
\end{align*}
proving (\ref{eqn:UFDpeup2}).
}
Now we can verify hypothesis \ref{itm:contract_gens_prime};
the ideal (\ref{eqn:specaq2_specAut3}) is already known to be prime and:
\claim{$(A/P)_A$ is $\S$-torsionfree.}
{
It suffices to check that $(R/\gen{\pinm{n}{m},p})_R$
is $\S$-torsionfree for all $m\in\Z$,
for if
$$\left( \sum_{m\in\Z}a_mv_m \right) u^it^j \in P,$$
then
$q^{2mi}a_mu^it^j \in \gen{\pinm{n}{m},p}$
for each $m\in\Z$.
By (\ref{eqn:UFDpeup}),
the problem further reduces to checking that $(R/\gen{\snj{n}{j},p})_R$ is $\S$-torsionfree for each $m\in\Z$ and $j\in\Jnm{n}{m}$.
$R/\gen{\snj{n}{j}}$ is isomorphic to $k[t]$ by an isomorphism that fixes $t$,
so $R/\gen{\snj{n}{j},p}\cong k[t]/\gen{p}$ is a domain and in particular $\S$-torsionfree.
}
For hypothesis \ref{itm:contract_gens_ore}, the nontrivial case to check is $g=\pinm{n}{m}v_m$ and $s=u^i$.
In this case we have
$$
gs=s(q^{2mi}g)\in g\S\cap sP.
$$
Therefore $P\triangleleft \A/\gen{r_n}$ is the contraction of (\ref{eqn:specaq2_specAut3}).
Pulling back to $\A$, we conclude that $\p\in\spec{k[t^\pm]}$
corresponds to
$$
\gen{\pinm{n}{m}v_m \bbel -n\leq m \leq n} + \gen{r_n} + \gen{\p\cap k[t]}\ \triangleleft\ \A.
$$
\done

\cor{
The algebra $\A$ is a noetherian UFD. (See \cite{UFDs} for the definition of noetherian UFD).
}
\begin{pf}
Having just listed all the prime ideals of $\A$, we simply check off the needed conditions:
\vspace{-1em}
\begin{itemize}
\item $\A$ is a noetherian domain.
\item Every nonzero prime ideal of $\A$ contains a nonzero principal prime ideal.
(Here a \emph{principal} ideal is one generated by a single normal element).
\emph{Proof:} For $T_1$ and $T_2$ this is obvious.
For $P\in T_{3n}$, $n\geq 1$, note that $P$ contains $\gen{r_n}\in T_2$.
\item Height one primes of $\A$ are completely prime.
\emph{Proof:}
Since $\gen{r_n}$ is properly contained in any $P\in T_{3n}$ for $n\geq 1$,
the primes in $T_{3n}$ are not height one.
We check that all the other primes are completely prime.
Suppose $P\in T_2$.
Then $P$ is generated in the commutative coefficient ring $k[u,t,d]$ of the GWA
$\A=k[u,t,d][x,y;\sigma,z]$ and it does not contain $z$,  so
Proposition \ref{prop:gwa_quot} shows that $\A/P$
is a GWA over a domain, and hence a domain.
For $P=\gen{u}+\gen{\p}\in T_1$, $\A/P$ is $k[u_{11},u_{12},u_{21}]/\p$, which is a domain.
\end{itemize}
\end{pf}
\done

Since $\A$ is noetherian, every closed subset of $\spec{\A}$ is a finite union of
irreducible closed subsets.
The topology of $\spec{\A}$ is therefore known if all inclusions of prime ideals are known.
We address in the following proposition 
those inclusions that  are not already expressed in Theorem $\ref{thm:specaq2}$.

\def\q{\mathfrak{q}}
\prop{\label{prop:specaq2_incl}
The inclusions among the prime ideals of $\A$ are as follows:
\vspace{-1em}\begin{enumerate}
\item\label{itm0:specaq2_incl}
Inclusions coming from the homeomorphisms
$T_1    \approx \spec{k[u_{11},u_{12},u_{21}]}$,
$T_2    \approx \spec{ k[t,d] }$, and
$T_{3n} \approx \spec{ k[t^\pm] }$ for $n\geq 1$.
\item \label{itm1:specaq2_incl}
Let $P\in T_1$. No prime in $T_2$ contains $P$,
and no prime in $T_{3n}$ contains $P$ for any $n$.
\item \label{itm2:specaq2_incl}
Let $P\in T_2$, say $P=\gen{\p}$ with $\p\in\spec{k[t,d]}$.
\begin{enumerate}
\item \label{itm2a:specaq2_incl}
The set of $Q\in T_1$ that contain $P$ is
$$
\{ \gen{u}+\gen{\q}
\bbel
\q\in\spec{k[u_{11},u_{12},u_{21}]}
\text{ and }
\p\subseteq \phi^{-1}(\q)
\},
$$
where $\phi$ is the
homomorphism $\phi:k[t,d]\rightarrow k[u_{11},u_{12},u_{21}]$
that sends $t$ to $u_{11}$ and $d$ to $u_{12}u_{21}$.
\item \label{itm2b:specaq2_incl}
Let $n\geq 1$.
The set of $Q\in T_{3n}$ that contain $P$ is
$$
\{
\gen{\pinm{n}{m}v_m \bbel -n\leq m \leq n} 
+ \gen{r_n} + \gen{\q\cap k[t]}
\ \bbel\ 
\q\in\spec{k[t^\pm]}
\text{ and }
\p\subseteq \eta_n^{-1}(\q)
\},
$$
where $\eta_n$ is the
homomorphism $\eta_n:k[t^\pm,d]\rightarrow k[t^\pm]$
that sends $t$ to $t$ and $d$ to $\frac{-q^{2n}}{(q^{2n}+1)^2}t^2$.
\end{enumerate}
\item \label{itm3:specaq2_incl}
Let $n\geq 1$ and let $P\in T_{3n}$, say
$$P=  
\gen{\pinm{n}{m}v_m \bbel -n\leq m \leq n} 
+ \gen{r_n} + \gen{\p\cap k[t]}$$
with $\p\in\spec{k[t^\pm]}$.
If $\p=0$, then
the only $Q\in T_1$ containing $P$ is
$$ \gen{u_{11},u_{22},u_{21},u_{12}}. $$
If $\p\neq 0$, then
no prime in $T_1$ or $T_2$ contains $P$, and no prime in $T_{3n'}$
contains $P$ for any $n' \neq n$.
\end{enumerate}
}
\begin{pf}
The inclusions of assertion \ref{itm0:specaq2_incl} are addressed by the homeomorphisms in Theorem \ref{thm:specaq2}.

\renewcommand{\stepheader}[1]{\textbf{#1:}}
\stepheader{\ref{itm1:specaq2_incl}}
If $P\in T_1$, then $u\in P$.
If $Q\in T_2$ 
then $Q_0$ 
(using the notation of Definition \ref{defn:gwa_idealseq})
is generated in $k[u,t,d]$ by elements of $k[t,d]$,
so $Q$ cannot contain $u$ and therefore cannot contain $P$.
If $Q\in T_{3n}$, then by definition $Q$ cannot contain $u$ and therefore cannot contain $P$.

\stepheader{\ref{itm2a:specaq2_incl}}
Assume the setup of assertion \ref{itm2a:specaq2_incl}.
Suppose that $Q\in T_1$, and write it
as $\gen{u}+\gen{\q}$
with $\q\in\spec{k[u_{11},u_{12},u_{21}]}$.
Then $P\subseteq Q$ if and only if $\gen{u}+\gen{\p}\subseteq Q$,
which holds if and only if $(\gen{u}+\gen{\p})/\gen{u}\subseteq Q/\gen{u}$
holds in $\A/\gen{u}$.
The following composite
is the homomorphism $\phi$ that we defined:
$$\begin{array}{ccccccc}
k[t,d]&\hookrightarrow &\A  &\twoheadrightarrow &\A/\gen{u} &\cong           &k[u_{11},u_{12},u_{21}]\\
\p    &                &Q   &\mapsto            &Q/\gen{u}  &\leftrightarrow &\q.
\end{array}$$
We see that $P\subseteq Q$ if and only if $\phi(\p)\subseteq\q$.
This holds if and only if $\p\subseteq\phi^{-1}(\q)$, so assertion \ref{itm2a:specaq2_incl} is proven.
\begin{excludeThis}
Hence the minimal $Q\in T_1$ containing $P$ correspond to the minimal
primes $\q$ over $\gen{\phi(\p)}$.
The following claim then proves assertion \ref{itm2a:specaq2_incl}.
\claim{The set of minimal primes over $\gen{\phi(\p)}$
is the fiber $\Phi^{-1}(\p)$ of $\Phi$ over $\p$.}{
Let $\q\in\spec{k[u_{11},u_{12},u_{21}]}$ be a minimal prime over $\p$.
If $\p=0$, then $\q=0$ and we have $\phi^{-1}(\q)=\p$.
If $\p$ is maximal, then we also have $\phi^{-1}(\q)=\p$.
So it remains to consider the case where $\p$ has height $1$.
It is then principal, and so is $\gen{\phi(\p)}$.
By Krull's principal ideal theorem, $\q$ must also be principal.
Write $\p$ as $\gen{f}$ and $\q$ as $\gen{g}$ for some
$f\in k[t,d]$ and $g \in k[u_{11},u_{12},u_{21}]$.
Write $g$ out as
$$ g= \sum_{i,j,l\geq 0} c_{ijl} u_{11}^i u_{12}^j u_{21}^l. $$
Since $g$ is a nonzero nonunit, $c_{i_0j_0l_0}\neq 0$
for some $i_0,j_0,l_0$ that are not all $0$.
Assume that either $i_0$ or $l_0$ is nonzero.
Since $k$ is infinite,
there is some $\alpha\in k^\times$ so that
\begin{equation}\label{eqn1:specaq2_incl} \sum_{j\geq 0} c_{i_0jl_0}\alpha^{l_0-j}\neq 0. \end{equation}
Let $\psi:k[u_{11},u_{12},u_{21}]\rightarrow k[t,d]$
be a homomorphism sending $u_{11}$ to $t$, $u_{12}$ to $\alpha^{-1}$, and $u_{21}$ to $\alpha d$.
Then $\psi\circ \phi = \text{id}$ and (\ref{eqn1:specaq2_incl}) is the coefficient of
$t^{i_0}d^{l_0}$ in $\psi(g)$.
Hence $\psi(g)$ is not a unit.
From $g|\phi(f)$, we have $\psi(g)|f$. 
Further, $\psi(g)$ must be an associate of $f$ because $\psi(g)$ is not a unit and $f$ is irreducible.
This proves that $\phi^{-1}(\q)=\p$:
$$ \phi^{-1}(\q) = \phi^{-1}(\gen{g}) \subseteq \gen{\psi(g)} = \gen{f} = \p. $$
If $i_0$ and $l_0$ were both $0$ above, then $j_0\neq 0$ and we could have obtained the same result by choosing
$\alpha\in k^\times$ such that $\sum_{l\geq0} c_{i_0j_0l}\alpha^{j_0-l}\neq 0$
and letting $\psi$ send $u_{11}$ to $t$, $u_{12}$ to $\alpha d$ and $u_{21}$ to $\alpha^{-1}$.

It remains to address the converse, that $\Phi(\q)=\p$ implies that
$\q$ is minimal over $\gen{\phi(\p)}$.
\note{cripes this part is not even true!}
}
\end{excludeThis}

\stepheader{\ref{itm2b:specaq2_incl}}
Assume the setup of assertion \ref{itm2b:specaq2_incl}.
Suppose that $Q\in T_{3n}$, and write it as
$$ 
\gen{\pinm{n}{m}v_m \bbel -n\leq m \leq n} 
+ \gen{r_n} + \gen{\q\cap k[t]}
$$
with $\q\in\spec{k[t^\pm]}$.
Then
$$ Q = \bigoplus_{m\in\Z} ( 
\gen{\pinm{n}{m},r_n} + \gen{\q\cap k[t]}
 ) v_m ; $$
the inclusion $\supseteq$ is clear and
the inclusion $\subseteq$ follows from the fact that the right hand side is an ideal of $\A$,
which can be verified by using
(\ref{eqn:s_obs_4}) and (\ref{eqn:specaq2_sigma_z})
to check that the conditions of Proposition  \ref{prop:gwa_homideals2} are met.
In particular, $Q_0=
\gen{\pinm{n}{0},r_n} + \gen{\q\cap k[t]} $.
Now assertion \ref{itm2b:specaq2_incl} is proven as follows:
\begin{align}
P\subseteq Q
&\ \Leftrightarrow\ \p\subseteq Q_0 = \gen{\pinm{n}{0},r_n}+\gen{\q\cap k[t]} \nonumber \\
&\ \Leftrightarrow\ \p\subseteq \gen{r_n}_{k[t,d]}+\gen{\q\cap k[t]}   \label{eqn1:specaq2_incl} \\
&\ \Leftrightarrow\ \p\subseteq \gen{r_n}_{k[t^\pm,d]}+\gen{\q}  \label{eqn2:specaq2_incl}  \\
&\ \Leftrightarrow\ ( \gen{r_n}_{k[t^\pm,d]} + \gen{\p}_{k[t^\pm,d]})/\gen{r_n}_{k[t^\pm,d]} 
   \subseteq (\gen{r_n}_{k[t^\pm,d]}+\gen{\q}) / \gen{r_n}_{k[t^\pm,d]}  \nonumber \\
&\ \Leftrightarrow\ \eta_n(\p)\subseteq \q  \label{eqn3:specaq2_incl}  \\
&\ \Leftrightarrow\ \p\subseteq \eta_n^{-1}(\q).  \nonumber
\end{align}
Line (\ref{eqn1:specaq2_incl})
is due to the fact that $Q_0\cap k[t,d] = \gen{r_n}_{k[t,d]}+\gen{\q\cap k[t]}$.
Line (\ref{eqn2:specaq2_incl})
is due to the fact that
$k[t,d]$ mod the ideal
$\gen{r_n}_{k[t,d]}+\gen{\q\cap k[t]}$ 
is $t$-torsionfree.
Line (\ref{eqn3:specaq2_incl})
is due to the fact that $\eta_n$ is
the following composite:
$$\begin{array}{ccccc}
k[t^\pm,d]&\twoheadrightarrow &k[t^\pm,d]/\gen{r_n} &\cong               &k[t^\pm]\\
          &                   &\gen{\q}+\gen{r_n}   &\leftrightarrow     &\q.
\end{array}$$

\stepheader{\ref{itm3:specaq2_incl}}
Assume the setup of assertion \ref{itm3:specaq2_incl}.
Let $Q\in T_1$ such that $P\subseteq Q$,
say $Q=\gen{u_{22}}+\gen{\q}$ with $\q\in\spec{k[u_{11},u_{12},u_{21}]}$.
Then $Q$ contains a power of $x$ and a power of $y$,
so $\q$ contains $u_{21}$ and $u_{12}$.
$Q$ also contains $r_n$, which is equivalent to $q^{2n}u_{11}$ modulo $\gen{u_{22},u_{12},u_{21}}$.
So $Q$ contains, and therefore equals, the maximal ideal $\gen{u_{11},u_{22},u_{21},u_{12}}$.
The containment $P\subseteq \gen{u_{11},u_{22},u_{21},u_{12}}$ clearly holds if $\p=0$.
But if $\p$ is nonzero, then it contains some polynomial in $t$ with nonzero constant term, which is not in
$\gen{u_{11},u_{22},u_{21},u_{12}}$.
Thus, nothing in $T_1$ contains $P$ when $\p\neq 0$.

If $Q\in T_2$, then $Q=\bigoplus_{m\in\Z}Q_0v_m$ does not contain any power of $x$.
So nothing in $T_2$ contains $P$.

Now suppose that $Q\in T_{3n'}$ with $n'\neq n$, and suppose for the sake of contradiction that $P\subseteq Q$.
Then $Q$ contains $r_n$ and $r_{n'}$.
Since $n\neq n'$, it follows that $t,d\in Q$.
Write $Q$ as
$$Q=  
\gen{\pinm{n'}{m}v_m \bbel -n'\leq m \leq n'} 
+ \gen{r_{n'}} + \gen{\q\cap k[t]}$$
with $\q\in\spec{k[t^\pm]}$.
Since $t\in Q$, we must have $\q=0$.
We have a contradiction:
\begin{equation*}
d\in Q_0 = \gen{\pinm{n'}{0},r_{n'}} \triangleleft k[u,t,d].
\tag*{\done}
\end{equation*}
\end{pf}

\prop{The algebra $\A$ does not have normal separation (see \cite[Ch 12]{NNR} for the definition).}
\begin{pf}
Let $P=\gen{\pinm{1}{0},x,y,r_1}$ and let $Q=\gen{r_1}$, both prime ideals of $\A$.
We will show that no element of $P\setminus Q$ is normal modulo $Q$.
Note that $k[u,t,d]/\gen{r_1}\cong k[u,t]$, and let $R=k[u,t]$.
Using Proposition \ref{prop:gwa_quot}, $\A/Q$ is isomorphic to
$$ W:=R[x,y;\sigma,z=-q^{-4}\snj{1}{1}\snj{1}{2}], $$
and $P/Q$ becomes \def\bP{\overline{P}}
$$ \bP := \gen{\pinm{1}{0}, x, y} = \bigoplus_{m>0} Ry^m \oplus \gen{\snj{1}{1}}_R \oplus \bigoplus_{m>0}Rx^m. $$
By Proposition \ref{prop:gwa_normalelements2},
the nonzero normal elements of $W$ are the $\sigma$-eigenvectors in $R$.
Thus, they are all of the form $u^if(t)$ for some polynomial $f(t)$ and some $i\in\N$.
But $\bP{}$ cannot contain such elements, since $\bP{}_0=\gen{\snj{1}{1}}_R$.
\end{pf}
\done

\prop{The algebra $\A$ is not catenary.}
\begin{pf}
Let $n\geq 1$.
The information in Proposition \ref{prop:specaq2_incl}
implies that the following two chains of primes are saturated:
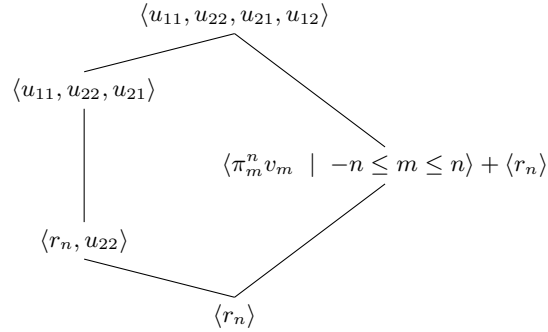
\newcommand{\de}[1]{{\the\dimexpr #1\relax}}
\def\w{2cm}    
\def\h{1cm}    
\def\c{7pt}  
\begin{center}
\begin{tikzpicture}
\draw (\de{\w},\de{\c}) -- (0,\de{\h-\c}) ;
\draw (0,\de{\h+\c}) -- (0,\de{ \h*3-\c}) ;
\draw (0,\de{ \h*3+\c}) -- (\de{\w},\de{ \h*4-\c}) ;
\draw (\de{\w},\de{\c}) -- (\de{ \w*2},\de{ \h*2-\c}) ;
\draw (\de{ \w*2},\de{ \h*2+\c}) -- (\de{\w},\de{ \h*4-\c}) ;
\node at (\de{\w},0) {$\gen{r_n}$};
\node at (0,\de{\h}) {$\gen{r_n,u_{22}}$};
\node at (0,\de{\h*3}) {$\gen{u_{11},u_{22},u_{21}}$};
\node at (\de{\w*2},\de{\h*2}) {$\gen{\pinm{n}{m}v_m\bbel -n\leq m\leq n}+\gen{r_n}$};
\node at (\de{\w*1},\de{\h*4}) {$\gen{u_{11},u_{22},u_{21},u_{12}}$};
\end{tikzpicture}
\end{center}
\end{pf}
\done

One reason to compute the prime spectrum of an algebra
is to make progress towards the lofty goal of knowing its complete representation theory.
The idea is to make progress by trying to know the algebra's \emph{primitive ideals},
those ideals that arise as annihilators of irreducible representations.
Since primitive ideals are prime, one approach is to determine the prime spectrum of the algebra
and then attempt to locate the primitives living in it.
The Dixmier-Moeglin equivalence, when it holds, provides a topological criterion for picking out primitives from the spectrum;
see \cite[II.7-II.8]{KGnB} for definitions.

\thm{
The algebra $\A$ satisfies the Dixmier-Moeglin equivalence, and its primitive ideals are as follows:
\vspace{-1em}
\begin{itemize}
\item The primitive ideals in $T_1$ are $\gen{u}+\gen{\p}$ for $\p\in\maxspec{k[u_{11},u_{12},u_{21}]}$.
\item The primitive ideals in $T_2$ are $\gen{\p}$ for $\p\in\maxspec{ k[t,d] }$.
\item The primitive ideals in $T_{3n}$ are $\gen{\pinm{n}{m}v_m \bbel -n\leq m \leq n} + \gen{r_n} + \gen{\p\cap k[t]}$
for $n\geq 1$ and $\p\in\maxspec{ k[t^\pm] }$.
\end{itemize}
}
\begin{pf}
We first observe that $\A$ satisfies the Nullstellensatz over $k$.
For this we can use \cite[II.7.17]{KGnB}, which applies because $\A$ is an iterated skew polynomial algebra
$$ \A \cong k[u_{11}][u_{22}][u_{12};\tau][u_{21};\tau',\delta'] $$
for a suitable choice of $\tau,\tau',\delta'$.
It then follows from \cite[II.7.15]{KGnB} that the following implications hold for all prime ideals of $\A$:
$$
\text{locally closed} \quad\Longrightarrow\quad \text{primitive}\quad \Longrightarrow\quad \text{rational}.
$$
To establish the Dixmier-Moeglin equivalence for $\A$,
it remains to close the loop and show that rational primes are locally closed.
We shall deal separately with the three different types of primes identified in Theorem \ref{thm:specaq2}.

\stepheader{$T_1$}
Suppose that $P\in T_1$, say $P=\gen{u}+\gen{\p}$ with $\p\in\spec{k[u_{11},u_{12},u_{21}]}$.
Then $\A/P \cong k[u_{11},u_{12},u_{21}]/\p$.
It follows that $P$ is rational if and only if $\p$ is a maximal ideal of $k[u_{11},u_{12},u_{21}]$.
In this case $P$ will be maximal and therefore locally closed. 
Thus, rational primes in $T_1$ are locally closed.

\newcommand{\Fract}[1]{\operatorname{Fract}(#1)}

\stepheader{$T_2$}
Suppose that $P\in T_2$, say $P=\gen{\p}$ with $\p\in\spec{k[t,d]}$.
Then, using Proposition \ref{prop:gwa_quot}, $\A/P$ is a GWA $R[x,y;\sigma,z]$, where $R:=k[u,t,d]/\gen{\p}$.
Since $z=d+q^{-2}tu-q^{-4}u^2$ is regular in $R$,
Proposition \ref{prop:gwa_laurent2} 
tells us that $R[x,y;\sigma,z]$ embeds
into the skew Laurent polynomial algebra $R[x^\pm;\sigma]$.
Let $K$ denote the fraction field of $R$.
The skew Laurent polynomial algebra $R[x^\pm;\sigma]$
embeds into the skew Laurent \emph{series} algebra $K((x^\pm;\sigma))$.
(We are abusing notation and writing $\sigma$ for the induced automorphism of $K$.)
Since the skew Laurent series algebra is a division ring,
we obtain an induced embedding of the Goldie quotient ring $\Fract{\A/P}$ into it:
$$ \Fract{\A/P} \hookrightarrow K((x^\pm;\sigma)). $$
For something to be in the center of $\Fract{\A/P} \cong \Fract{ R[x,y;\sigma,z] }$, it must at least commute
with $R$ and $x$.
This is sufficient to place it in the center of $K((x^\pm;\sigma))$, so
\begin{equation}\label{eqn:DM1} Z(\Fract{\A/P}) \cong Z(K((x^\pm;\sigma)))\cap\Fract{\A/P}. \end{equation}
According to Proposition \ref{prop:Z_skew}, the center of $K((x^\pm;\sigma))$ is the fixed subfield $K^\sigma$.
Since $K$ is wholly contained in $Z(\Fract{\A/P}) \cong Z(\Fract{R[x,y;\sigma,z]})$, (\ref{eqn:DM1}) becomes
$$ Z(\Fract{\A/P}) \cong K^\sigma. $$
Now to compute $K^\sigma$.
Since
$$
R=k[u,t,d]/\gen{\p} \cong (k[t,d]/\p)[u],
$$
$K$ is the rational function field $L(u)$, where $L$ is the fraction field of $k[t,d]/\p$.
\claim{$K^\sigma=L$.}{
Observe that $\sigma$ fixes $L$ and sends $u$ to $q^2 u$.
Consider any nonzero $f/g\in K^\sigma = L(u)^\sigma$, where $f,g\in L[u]$ are coprime.
We have $\sigma(f)g=f\sigma(g)$. Since $f$ and $g$ are coprime, it follows that $f\mid\sigma(f)$.
Similarly, since $\sigma(f)$ and $\sigma(g)$ are coprime, $\sigma(f)\mid f$.
It follows that $\sigma(f)=\alpha f$ for some $\alpha\in L$.
From $\sigma(f)g=f\sigma(g)$ it follows that also $\sigma(g)=\alpha g$.
We have an eigenspace decomposition for the action of $\sigma$ as an $L$-linear operator on $L[u]$;
it is $\bigoplus_{i\geq 0} Lu^i$, with
distinct eigenvalues since $q$ is not a root of unity.
Since $f$ and $g$ are $\sigma$-eigenvectors with the same eigenvalue $\alpha$,
there is some $i\geq 0$ such that $f=f_0u^i$ and $g=g_0u^i$, where $f_0,g_0\in L$.
Thus, $f/g=f_0/g_0\in L$.
}
We have found that
$$ Z(\Fract{\A/P}) \cong L. $$
The fraction field $L$ of $k[t,d]/\p$ is algebraic over $k$ if and only if $\p\triangleleft k[t,d]$ is maximal.
Thus, $P$ is rational if and only if $\p$ is maximal. 

Now assume that $P$ is rational and hence that $\p$ is maximal.
For any $\q\in\spec{k[t^\pm]}$ and any $n\geq 1$, define
$$Q_{\q,n} := \gen{\pinm{n}{m}v_m \bbel -n\leq m \leq n} + \gen{r_n} + \gen{\q\cap k[t]} \in T_{3n}.$$
If no $Q_{\q,n}$ contains $P$,
then
by using Proposition \ref{prop:specaq2_incl} we can see that
$\{P\} = V(P)\cap(\spec{\A}\setminus V(u))$,
so $P$ is locally closed.
Suppose, on the other hand, that $Q_{\q,n}$ contains $P$ for some $\q\in\spec{k[t^\pm]}$ and $n\geq 1$.
We claim that this can occur for at most one $n$.
\claim{If $Q_{\q,n}$ and $Q_{\q',n'}$ both contain $P$, then $n=n'$.}{
According to assertion \ref{itm2b:specaq2_incl} of Proposition \ref{prop:specaq2_incl},
we have
\begin{align*}
&\p\subseteq \eta_n^{-1}(\q)
\text{\ \ and\ \ }
\p\subseteq \eta_{n'}^{-1}(\q'),
\end{align*}
where $\eta_n,\eta_{n'}$ are the homomorphisms defined there.
Since $\p$ generates a maximal ideal of $k[t^\pm,d]$, this forces
$$ \eta_n^{-1}(\q) = \eta_{n'}^{-1}(\q'). $$
We have
$$
d+\frac{q^{2n}}{(q^{2n}+1)^2}t^2,\ 
d+\frac{q^{2n'}}{(q^{2n'}+1)^2}t^2\in
\eta_n^{-1}(\q) = \eta_{n'}^{-1}(\q'),
$$
so
$$
\left(\frac{q^{2n}}{(q^{2n}+1)^2} - 
\frac{q^{2n'}}{(q^{2n'}+1)^2}\right) t^2 \in \eta_n^{-1}(\q) = \eta_{n'}^{-1}(\q').
$$
Since we cannot have $\eta_n^{-1}(\q) = k[t^\pm,d]$,
the quantity in parentheses must vanish.
This leads to the equation
$$
0 = q^{2n}(q^{2n'}+1)^2 - q^{2n'}(q^{2n}+1)^2
 =  (q^{n'}-q^n)(q^{n'}+q^n)(q^{n+n'}-1)(q^{n+n'}+1).
$$
Since $q$ is not a root of unity, it follows that $n=n'$.
}
Hence $\{P\} = V(P)\cap (\spec{\A}\setminus (V(u)\cup V(x^n)))$,
and again we see that $P$ is locally closed.
Thus we have shown that all rational primes in $T_2$ are locally closed.

\stepheader{$T_{3n}$}
Suppose that $n\geq1$ and $P\in T_{3n}$, say $P=\gen{\pinm{n}{m}v_m \bbel -n\leq m \leq n} + \gen{r_n} + \gen{\p\cap k[t]}$
with $\p\in\spec{k[t^\pm]}$.
We will show that if $P$ is rational, then $\p$ must be maximal.
Assume that $P$ is rational, but $\p$ is not maximal (i.e. $\p=0$).
Then $t\in Z(\Fract{\A/P})$ is algebraic over $k$, so 
for some nonzero polynomial $f(T)\in k[T]$ we must have $f(t)=0\in Z(\Fract{\A/P})$.
For the element $t$ of $\A$, this means that $f(t)\in P$.
We can describe $P$ explicitly in terms of its homogeneous components:
$$
P=\bigoplus_{m\in\Z}\gen{\pinm{n}{m},r_n}v_m.
$$
The inclusion $\supseteq$ is obvious, and the equality can be verified
by checking that the conditions of Proposition \ref{prop:gwa_homideals2} are met by the right hand side
and it indeed defines a two-sided ideal.
So the fact that $f(t)\in P$ can be refined to $f(t)\in \gen{\pinm{n}{0},r_n}_{k[u,t,d]}$.
Pushing this fact into $k[u,t,d]/\gen{r_n}\cong k[u,t]$ gives
$$
f(t)\in \gen{\pinm{n}{0}}_{k[u,t]}
$$
which is clearly false. 

Thus we have shown that when $P$ is rational, $\p\triangleleft k[t^\pm]$ must be maximal.
Using Proposition \ref{prop:specaq2_incl}, we see that in this case $\{P\}=V(P)$.
So all rational primes in $T_{3n}$ are locally closed.

We have now shown that all rational prime ideals of $\A$ are locally closed, and 
we conclude that $\A$ satisfies the Dixmier-Moeglin equivalence.
Further, we have pinpointed which primes are rational in $T_1$ and $T_2$.
As for $T_{3n}$, we have found for $P=\gen{\pinm{n}{m}v_m \bbel -n\leq m \leq n} + \gen{r_n} + \gen{\p\cap k[t]}$ that
$$
\text{$P$ rational}\quad\Longrightarrow\quad\text{$\p$ maximal}\quad\Longrightarrow\quad\text{$P$ locally closed}.
$$
Putting this information together and applying the Dixmier-Moeglin equivalence, 
we conclude that the primitive ideals of $\A$ are as stated in the theorem.
\end{pf}
\done

\section{Appendix}

\def\S{\mathcal{S}}

There are a few aspects of noncommutative localization that make an appearance throughout this work
and that rely on the noetherian hypothesis.
For the reader's convenience, we lay them out here.
Proposition \ref{prop:loc_quot} says that
localization ``commutes'' with factoring out an ideal, and it is a standard fact.
Theorem \ref{thm:spec_loc} says that the usual correspondence of prime ideals along a localization
is a homeomorphism, also a standard fact.
Finally, Lemma \ref{lem:contract_gens}
provides a way to describe the pullback of a prime ideal along a localization
by using a ``nice'' generating set.

\prop{\label{prop:loc_quot}
Let $\S$ be a right denominator set in a right noetherian ring $R$.
Let $I$ be an ideal of $R$, with extension $I^e$ to $R\S^{-1}$.
Then:
\vspace{-1em}
\begin{itemize}
\item $I^e$ is an ideal of $R\S^{-1}$.
\item $\bar{\S}:=\{s+I\bbel s\in \S\}$ is a denominator set of $R/I$.
\item The canonical homomorphism $\phi:R/I\rightarrow (R\S^{-1})/I^e$
gives a right ring of fractions for $R/I$ with respect to $\bar{\S}$.
That is, there is an isomorphism 
$
\bar{\phi}:
(R/I)\bar{\S}^{-1}
\cong
(R\S^{-1})/I^e
$
making the following diagram commute:
$$ \begin{tikzcd}[ampersand replacement=\&]
R \arrow[d,two heads] \arrow[r,"\text{loc}"] \& R\S^{-1} \arrow[r,two heads,"\text{quo}"] \& R\S^{-1}/I^e   \\
R/I \arrow[rr,"\text{loc}"] \arrow[urr, "\exists ! \phi"] \& \& (R/I)\bar{\S}^{-1} \arrow[u,"\cong"', "\bar{\phi}"]
\end{tikzcd}$$
\end{itemize}
}
\begin{thesisOnly}
\begin{pf}
By \cite[Theorem 10.18a]{NNR}, $I$ extends to an ideal of $R\S^{-1}$.
It is clear that $\bar{\S}$ is a right Ore set of $R/I$.
It is automatically right reversible due to the noetherian hypothesis; see \cite[Proposition 10.7]{NNR}.
We will use the universal property that characterizes rings of fractions (e.g. see \cite[Proposition 10.4]{NNR}).
The homomorphism $\phi$ is uniquely defined because $I$ is in the kernel of the upper row of the diagram.
Since $\phi$ maps $\bar{\S}$ to a collection of units,
$\bar{\phi}$ is uniquely defined such that the diagram commutes.
It is surjective because $rs^{-1}+I^e$ is the image of $(r+I)(s+I)^{-1}$ for $r\in R$ and $s\in\S$.
For injectivity, suppose that $rs^{-1}\in I^e$. Then $r1^{-1}\in I^e$,
so $r\in I^{ec}$, the contraction of $I^e$ to $R$.
By \cite[Theorem 10.15b]{NNR}, this implies that $rs'\in I$ for some $s'\in\S$.
That is, $0=(r+I)1^{-1}\in (R/I)\bar{\S}^{-1}$.
\end{pf}
\done
\end{thesisOnly}

\thm{
\label{thm:spec_loc}
Let $\S$ be a right denominator set in a right noetherian ring $R$.
Then contraction and extension of prime ideals are inverse homeomorphisms:
\begin{equation}\label{eqn:spec_loc}
\spec{R\S^{-1}}\ \approx\ \{ Q\in\spec{R} \bbel Q\cap \S=\emptyset \}.
\end{equation}
}
\begin{thesisOnly}
\begin{pf}
The bijection of sets in (\ref{eqn:spec_loc}) is given to us by \cite[Theorem 10.20]{NNR},
so we only need to show that closed sets are preserved.
The closed subsets of $\{ Q\in\spec{R} \bbel Q\cap \S=\emptyset \}$
are
$$V_\S(J):=\{Q\in\spec{R}\bbel Q\supseteq J\ \text{and}\ Q\cap \S=\emptyset\},$$
for ideals $J$ of $R$.
Observe that (\ref{eqn:spec_loc}) preserves finite unions and inclusions of prime ideals.
Since $R$ is right noetherian,
every ideal of $R$ has finitely many minimal primes over it,
and the same goes for $R\S^{-1}$.
It follows that the topological spaces in (\ref{eqn:spec_loc}) are noetherian,
so every closed set is a finite union of irreducible closed sets.
Hence it suffices to show that (\ref{eqn:spec_loc}) preserves irreducible closed sets.
Irreducible closed sets of $\spec{R\S^{-1}}$ are of the form $V(P)$ for primes $P\triangleleft R\S^{-1}$.
Irreducible closed sets of $\{ Q\in\spec{R} \bbel Q\cap \S=\emptyset \}$ are of the form $V_\S(P)$ for primes
$P\triangleleft R$ such that $P\cap\S=\emptyset$.
Since (\ref{eqn:spec_loc}) preserves inclusions of prime ideals,
it sends any $V_\S(P)$ to $V(P\S^{-1})$.
\end{pf}
\done
\end{thesisOnly}

\def\G{\mathcal{G}}
\def\S{\mathcal{S}}
\lem{ \label{lem:contract_gens}
Let $R$ be a right noetherian ring,
$\S\subseteq R$ a right denominator set,
and $\phi:R\rightarrow R\S^{-1}$ the localization map.
Let $\G\subseteq R$ and assume the following:
\vspace{-1em}
\begin{enumerate}
\item \label{itm:contract_gens_2sided}
The \emph{right}
ideal $P$ generated by $\G$ is a two-sided ideal of $R$.
\item \label{itm:contract_gens_prime}
Either $P$ is a prime ideal of $R$ disjoint from $\S$, or $\gen{\phi(\G)}$ is a prime ideal of $R\S^{-1}$ and $(R/P)_R$ is $\S$-torsionfree.
\item \label{itm:contract_gens_ore}
For all $g\in\G$ and $s\in\S$,
$$ g \S\cap s P \neq \emptyset .$$
\end{enumerate}
Then
$$ P = \phi^{-1}(\gen{\phi(\G)}). $$
That is, the ideal of $R\S^{-1}$ generated by $\phi(\G)$ contracts to the ideal of $R$ generated by $\G$.
}
\begin{pf}
Assumption \ref{itm:contract_gens_ore}
guarantees that the right ideal of $R\S^{-1}$ generated by $\phi(\G)$ is a two-sided ideal.
Let superscripts ``e'' and ``c'' denote extension and contraction of ideals along $\phi$.
Observe that
\begin{align}
\gen{\phi(\G)}
&=\{\sum_{i=1}^n\phi(g_i) \phi(r_i) \phi(s_i)^{-1}
\bbel n\in\N,\ r_i\in R,\ s_i\in \S,\ g_i\in\G \text{\ for $1\leq i \leq n$}\}
\nonumber \\
&=\{\sum_{i=1}^n\phi(g_i r_i) \phi(s)^{-1}
\bbel s\in\S,\ n\in\N,\  r_i\in R,\ g_i\in\G \text{\ for $1\leq i \leq n$}\}
\label{eqn:contract_gens_explain1} \\
&=\{\phi(a) \phi(s)^{-1}
\bbel s\in\S,\ a\in \gen{\G} \}
 = P^e.
\nonumber
\end{align}
In line (\ref{eqn:contract_gens_explain1}),
we used the fact that it is possible to get a ``common right denominator''
for a finite list of right fractions; see \cite[Lemma 10.2]{NNR}.
Now assumption \ref{itm:contract_gens_prime} implies that $P$ is a prime ideal
of $R$ disjoint from $\S$, either trivially or by \cite[Theorem 10.18b]{NNR}.
To finish, we use the correspondence between prime ideals disjoint from $\S$ and prime ideals of $R\S^{-1}$:
\begin{align*}
P
&= P^{ec}
=\gen{\phi(\G)}^c = \phi^{-1}(\gen{\phi(\G)}).
\tag*{\done}
\end{align*}
\end{pf}

Note that assumption \ref{itm:contract_gens_ore} of Lemma \ref{lem:contract_gens}
holds trivially whenever $\G$ or $\S$ is central.

\begin{thesisOnly}
\subsection{Chinese Remainder Theorem}

This section serves to clarify the way in which the Chinese Remainder Theorem (CRT)
gets applied in section \ref{sec:primespwrx}.

Let $R$ be a commutative ring with pairwise comaximal ideals $I_1,\ldots,I_n$,
and let $\Pi:=I_1\cdots I_n$ denote their product.
The CRT says that the homomorphism
\begin{equation} \label{eqn:CRT} R/\Pi \rightarrow R/I_1 \times \cdots \times R/I_n, \end{equation}
whose components are induced by canonical projections, is an isomorphism.
This implies that there is a bijective correspondence between
ideals $J$ of $R$ that contain $\Pi$, and tuples $(J_1,\ldots,J_n)$ of ideals of $R$
such that $J_i\supseteq I_i$ for all $i$.
Let us describe the correspondence explicitly.

There are $e_1,\ldots e_n\in R$
such that $e_i \equiv 1$ mod $I_i$ and $e_i \equiv 0$ mod $I_j$ for $j\neq i$.
These are just the pairwise orthogonal idempotents (when taken modulo $\Pi$)
corresponding to the ring decomposition in (\ref{eqn:CRT});
they also satisfy
$$ e_1 + \cdots + e_n \equiv 1 \ \ \text{mod $\Pi$}. $$
$(R/\Pi)(e_i+\Pi)$ is the copy of the ring $R/I_i$ that is (non-unitally)
contained in $R/\Pi$ via (\ref{eqn:CRT}).
Explicitly, the correspondence of elements is given by
\begin{equation}\label{eqn:CRT2}(re_i+\Pi)\ \leftrightarrow\ (r+I_i).\end{equation}

Let $J$ be an ideal of $R$ that contains $\Pi$.
Projecting $J/\Pi$ to the copy $(R/\Pi)(e_i+\Pi)$ of $R/I_i$
yields $(J/\Pi)(e_i+\Pi)$, which corresponds to $(J+I_i)/I_i \triangleleft R/I_i$ via (\ref{eqn:CRT2}).
Going back in the other direction, suppose that $J_1,\ldots,J_n$ are ideals of $R$
such that $J_i\supseteq I_i$ for all $i$.
Then $J_i/I_i$ is an ideal of $R/I_1 \times \cdots \times R/I_n$,
and it corresponds to $(J_ie_i+\Pi)/\Pi\triangleleft R/\Pi$ via (\ref{eqn:CRT2}).
The sum of these over $i$ is the ideal of $R/\Pi$ corresponding to the tuple $(J_1,\ldots,J_n)$.
Thus we have the following explicit description of how ideals are carried across (\ref{eqn:CRT}):

\prop{\label{prop:CRT_ideals}
Let $R$ be a commutative ring with pairwise comaximal ideals $I_1,\ldots,I_n$,
and let $\Pi:=I_1\cdots I_n$ denote their product.
Let $e_1,\ldots e_n\in R$ be such that $e_i \equiv 1$ mod $I_i$ and $e_i \equiv 0$ mod $I_j$ for $j\neq i$.
There is a bijective correspondence between ideals $J$ of $R$ that contain $\Pi$, and tuples $(J_1,\ldots,J_n)$ of ideals of $R$
such that $J_i\supseteq I_i$ for all $i$, and it is given by:
$$
\begin{array}{lll}
J & \mapsto & (J+I_1,\ldots,J+I_n)\\
J_1e_1+\cdots+J_ne_n+\Pi & \mapsfrom & (J_1,\ldots,J_n).
\end{array}
$$
}

Another use we have for the CRT is primary decomposition for modules annihilated by a product of maximal ideals:

\def\m{\mathfrak{m}}
\def\a{\mathfrak{a}}
\prop{
\label{prop:CRT_primary}
Let $R$ be a commutative ring, $X$ an $R$-module annihilated by
$\a=\m_1^{i_1} \cdots \m_s^{i_s}$,
where $\m_1,\ldots,\m_s$ are distinct maximal ideals of $R$
and $i_1,\ldots,i_s\in\Z_{>0}$.
Then $X=X_1\oplus\cdots\oplus X_s$ where $X_j:=\ann_X(\m_j^{i_j})$.
}
\begin{pf}
By the CRT,
$ R/\a \cong (R/\m_1^{i_1}) \times \cdots \times (R/\m_s^{i_s}). $
This gives a complete set $e_1,\ldots,e_s\in R/\a$ of pairwise
orthogonal idempotents such that $\ann_R(e_j)=\m_j^{i_j}$ for $1\leq j \leq s$.
Note that $X$ can be viewed as a module over $R/\a$.
It is easily checked that $X=e_1X\oplus \cdots e_s X$ as an $(R/\a)$-module, 
and hence also as an $R$-module.
Since $\m_j^{i_j}e_jX=0$, the inclusion $e_jX\subseteq \ann_X(\m_j^{i_j})=:X_j$ holds
for all $1\leq j\leq s$.
For the reverse inclusion, suppose that $\m_j^{i_j}x=0$
for a given $x\in X$ and $1\leq j\leq s$.
When $l\neq j$, $\ann_Re_lx$ contains both $\m_j^{i_j}$ and $\m_l^{i_l}$,
and is therefore equal to $R$. Hence $e_lx=0$ for $l\neq j$,
and $x=e_jx\in X_j$.
\end{pf}
\done
\end{thesisOnly}

Department of Mathematics,
University of California,
Santa Barbara,
CA 93106,
USA\\\ \\
ebrahim@math.ucsb.edu

\end{document}